\documentclass[12pt, a4paper]{article}
\usepackage[margin=1in]{geometry}
\usepackage{comment}

\RequirePackage[OT1]{fontenc}
\RequirePackage{amsthm,amsmath}
\RequirePackage[numbers]{natbib}

\usepackage[english]{babel}
\usepackage{graphicx}
\usepackage{float}
\usepackage{amsfonts}
\usepackage{enumerate}
\usepackage{amsmath}
\usepackage{comment}
\usepackage{authblk}
\usepackage{mathtools}
\usepackage{tikz}

\numberwithin{equation}{section}
\theoremstyle{plain}
\newtheorem{thm}{Theorem}[section]
\newtheorem{prop}[thm]{Proposition}
\newtheorem{cor}[thm]{Corollary}
\newtheorem{lem}[thm]{Lemma}
\theoremstyle{definition}
\newtheorem{dfn}[thm]{Definition}

\newtheorem{rem}[thm]{Remark}
\theoremstyle{remark}

\newcommand{\probSpace}{(\Omega,\mathcal{F},\{\mathcal{F}_{t}\}_{t\geq 0},\mathbb{P})}

\newcommand{\Vcand}{\widetilde{V}}

\graphicspath{{../IMAGES/}}
\def\B{\mathbb B}
\def\R{\mathbb R}
\def\N{\mathbb N}
\def\Z{\mathbb Z}

\def\E{\mathbb E}
\def\P{\mathbb P}
\def\Q{\mathbb Q}
\def\sha{{\cal A}}

\def\shf{{\cal F}}

\def\shl{{\cal L}}
\def\shp{{\cal P}}

\def\shs{{\cal S}}

\def\1{\mathds{1}}
\def\ra{\rightarrow}

\def\ve{\varepsilon}

\def\s{\sigma}
\def\t{\theta}
\def\d{\delta}
\def\T{\Theta}
\def\piq{\frac{\pi}{4}}

\def\pitq{\frac{3\pi}{4}}
\def\picq{\frac{5\pi}{4}}
\def\pim{\frac{\pi}{2}}
\def\pitm{\frac{3\pi}{2}}
\def\pisq{\frac{7\pi}{4}}

\def\sa{A^\star}

\def\um{\frac{1}{2}}

\newcommand{\ind}{1}

\usepackage{geometry}
\usepackage{todonotes}

\begin{document}

\begin{center}
{\Large \bf  Stochastic Control Problems Motivated by Sailboat Trajectory Optimization}
\vskip 16pt

\bf Carlo Ciccarella\footnote{Institut de Mathématiques, Ecole Polytechnique Fédérale de Lausanne, Station 8, CH-1015, Lausanne, Switzerland. 
}
\hskip .5cm Robert C.~Dalang$^1$ \hskip .5cm
Laura Vinckenbosch$^{1,}$\footnote{School of Engineering and Management Vaud (HEIG-VD), HES-SO University of Applied Sciences and Art Western Switzerland, CH-1401, Yverdon-les-Bains, Switzerland.

The research of R.C.~Dalang is partly supported by the Swiss National Foundation for Scientific Research.}

\end{center}





\begin{abstract}
We develop a mathematical model for sailboat navigation that captures the essential features of the problem and can provide insights that might not be available otherwise. In our model, the motion of the sailboat, which would travel at speed $v>0$ in a constant wind, is the solution of a system of two stochastic differential equations driven by a Brownian motion on a circle with speed $\sigma > 0$. We formulate two stochastic control/optimal switching problems, in which the objective is to reach a circular upwind target of radius $\eta \geq 0$ as quickly as possible. In the first problem, there is a tacking cost $c > 0$, so this is an impulse control problems, while in the second problem,  we assume that $c=0$ and singular controls are needed. We establish the viability of both models (assuming that $\eta > 0$ in the second model), that is, their value functions are finite, and we obtain bounds on these value functions related to the parameters of the problem. In the second problem, since the  state equation for the optimally controlled motion has discontinuous coefficients and is driven by a degenerate diffusion, standard results on existence and uniqueness of strong solutions do not apply: we provide a proof via the Yamada-Watanabe argument.
\end{abstract}

{\small
\medskip

}

\section{Introduction}
The progress of a sailboat navigating upwind depends both on the navigation strategy chosen by the skipper and on the variations in wind speed and wind direction, much in the same way that the wealth of an investor depends both on the investor's investment strategy and on the fluctuations of the stock market. Indeed, wind variations are well-modeled by random processes, and on certain time scales, even by diffusion processes, including the fundamental example of Brownian motion (see \cite{DalangDumas2015,DUM,LV}). 

The objective of this paper is to develop a mathematical model that is similar in spirit to the Black and Scholes model \cite{black-scholes1974} in mathematical finance: this famous model is not based on data, is not realistic for describing the fluctuations of any particular stock and is described using only two free parameters (the {\em mean return rate} and the {\em volatility}), and yet it 
provides key insights into stock behavior and investment strategies, in addition to being the source of many interesting mathematical questions. Our intention here is to develop a model for upwind sailboat navigation, 
designed to capture the essential features of upwind sailboat navigation rather than empirical details, and which can provide insights that might not be available otherwise. For instance, it should help answer questions concerning how strategies and sailboat travel times depend on wind variability. In addition, just as it is the case for the Black and Scholes model, our model can be refined and extended so as to provide a good model for any specific sailboat navigating under realistic wind conditions.

This model should also be mathematically attractive and a source of mathematical questions that we hope many mathematicians will find interesting. Therefore, this paper contains a presentation of the mathematical model that we are proposing, with an explanation of the choices that we have made, a derivation of the basic system of stochastic differential equations that govern sailboat motion, and some first mathematical results that aim to show the potential of this model. Indeed, the model is quite rich and, depending on the choice of four parameters, leads to numerous questions in optimal stochastic control, including impulse control and singular stochastic control problems (in the terminology of \cite{OS}).

This approach is quite different from that taken in most existing papers on yacht routing problems. These contain mainly numerical studies, beginning with those of A.B.~Philpott and coauthors \cite{Phil1}--\cite{Phil4}. In particular, \cite{Phil3} presents two models for optimizing yacht routes. The first one is intended for short distance races and the wind is modeled as a Markov chain; the second one treats ocean races and the wind is modeled using long-term weather forecasts. The authors use numerical methods to determine tacking strategies which minimize the expected travel time. In \cite{DalangDumas2015,DUM}, the authors use statistical analysis of wind fluctuations to determine the transition probabilities of a Markov chain, then discretize the race field, the relevant wind directions and the boat's motion in order determine the fastest route via stochastic optimization methods. The paper \cite{EF} proceeds by solving numerically a Markov decision problem.

More recently,  see for example~\cite{FF2019,socs,rasc}, the yacht-routing problem has been formulated using a hybrid control approach and solved numerically with semi-Lagrangian iterative schemes. Other works in the sailing race literature consider the more complex question of winning a race against an opponent, in which getting to the target {\em first} is the main objective, rather than trying to reach the target as quickly as possible (see \cite{ CFF, Rossel2013, Tagliaferri2014}). This two-boat problem brings into play the mathematics of game theory and will not be considered here. 
\medskip

\noindent{\em Main results of the paper}
\medskip

We construct a model for wind behavior and sailboat motion in which the motion of the sailboat is the solution of a system of two stochastic differential equations (SDEs) driven by a Brownian motion on the unit circle. These equations are particularly compact when written in polar coordinates (see \eqref{eds_proc_controle_polar}). However, the equation for the radial component is degenerate, since it does not contain a diffusion term, and the coefficients become infinite as the radius decreases to $0$. 

We then formulate two optimal control problems. In both, the objective is to reach a circular target of radius $\eta \geq 0$ as quickly as possible, when this target is located upwind from the boat. Indeed, in this case, the skipper must make a sequence of important decisions, each of which consists in switching headings ({\em tacking}) so that the wind enters the sails from the left instead of the right, and vice-versa. In the first problem, we assume that the loss of time due to this change in direction is $c > 0$, while in the second problem, this loss of time is set to zero ($c=0$ may provide a good approximation for long-distance races, in which the loss of time due to tackings is negligible relative to the total travel time). The Hamilton-Jacobi-Bellman equations for these two problems are derived.

The first question of interest is the viability of these problems, that is, is their value function finite? We answer this question affirmatively for both problems, by analyzing a particular impulse strategy and showing that with this strategy, only a finite number of tacks are needed and the expected time to reach the target is finite in both problems. For the second problem, under the additional assumption that $\eta > 0$, we obtain a sharper upper bound on the value function that corresponds to a candidate optimal strategy. Under this strategy, the boat will perform an infinite number of tacks (as in the optimal strategy described in \cite{MSV90}; since the boat's motion under this strategy solves a system of SDEs with discontinuous but monotone coefficients, we use the Yamada-Watanabe argument \cite[Chapter 5, Section D]{ks} to establish existence and uniqueness of a strong solution. The optimality of this strategy is proved in the companion paper \cite{CC1}.

There is a well-developed literature for stochastic control problems of the types discussed in this paper, including \cite{Fleming75,FlemingSoner06,Krylov77,OS,pham,touzi,YongZhou99}. Most frequently, these problems involve both terminal costs and running costs. The book \cite[Chapter 6]{OS} contains several examples of impulse control problems; another is considered in \cite{OneD}. Our problem falls into the class of {\em optimal switching problems,} for which a theory is developed in the book of H.~Pham \cite{pham}; however, our problem is not standard, in part because it concerns minimizing the expected time, rather than an integral functional as in \cite{pham}.   Examples of singular control problems are given in \cite[Chapter 5]{OS} and \cite[Chapter 3]{pham}; others can be found in \cite{DS,Ma87,MSV90}. Contrary to \cite{Ma87,MSV90}, we do not use multiparameter processes in our formulations of the two control problems, since the formulation via a system of stochastic differential equations is quite natural.

\medskip

\noindent{\em Organization of the paper}
\medskip

In Section \ref{Sec:2}, we discuss the motion of a sailboat and the behavior of wind, then we present our model and derive the equations of motion. In Section \ref{Sec:3}, we formulate the two control problems and present a verification theorem. In Section \ref{Sec:4}, we derive bounds on the value functions in both problems. Finally, in Section \ref{Sec:5}, for $\eta > 0$, we analyze our candidate optimal strategy, we prove that this strategy is admissible, and we use it to provide a better bound on the value function in the first problem.

\bigskip

\noindent{\sc Acknowledgement.} This article is based on parts of the Ph.D.~theses of Laura Vinckenbosch and Carlo Ciccarella, in particular on \cite[Chapter 5]{LV} and \cite[Chapter 3]{CC}, written under the supervision of Robert C.~Dalang. A discussion of a near-final version of this paper with Claude Godr\`eche (Université Paris-Saclay) led to improvements in the interpretation of the mathematical results and motivated the writing of  Appendix A.

The authors thank an anonymous referee for suggesting in particular the inclusion of Theorems \ref{rd06_29t1}, \ref{rd07_01t1} and Appendix B.

\section{Modeling the motion of a sailboat}\label{Sec:2}

In this section, we begin with some general considerations about the motion of a sailboat, we recall some basic sailing terminology  and we give a simplified description of the behavior of wind. 
Then we describe the model and its assumptions and we derive the system of SDEs that describes the position process of the moving yacht,  both in Cartesian coordinates and in polar coordinates. 

\subsection{Basic features of sailboat motion and sailing terminology}\label{sec2.0}

The speed of a given sailboat depends on many factors, including in particular the wind speed, the wind direction, and the angle between the boat's heading and the wind direction. The function that relates this angle to the boat's speed is rather complex (see Section \ref{sec2.1}) and specific to each particular sailboat. There are also many other factors that we will not be discussing here, such as the effects of waves, currents, choice of sails, etc.

Another central feature of upwind sailboat navigation is that a sailboat {\em cannot} sail directly into the wind, and therefore, a sailboat that seeks to travel upwind must navigate at some non-zero angle with respect to the wind direction. And one of the most critical choices is to decide when to {\em tack}, that is, to change course in such a way that, if the wind was coming from one side of the boat just before the tack, then it will come from the other side just after the tack. Since this maneuver slows down the yacht, the skipper generally seeks to tack infrequently and at suitably chosen times.

We assume that the boat seeks to reach as quickly as possible a circular target (say, a buoy, or a small island) located upwind from the boat, with nonnegative radius $\eta \geq 0$. If $\eta = 0$, then the target is just a single point. 
\medskip

\noindent{\em Tacking and jibing}
\medskip

\begin{figure}[h]
\begin{center}
\includegraphics[height=4cm]{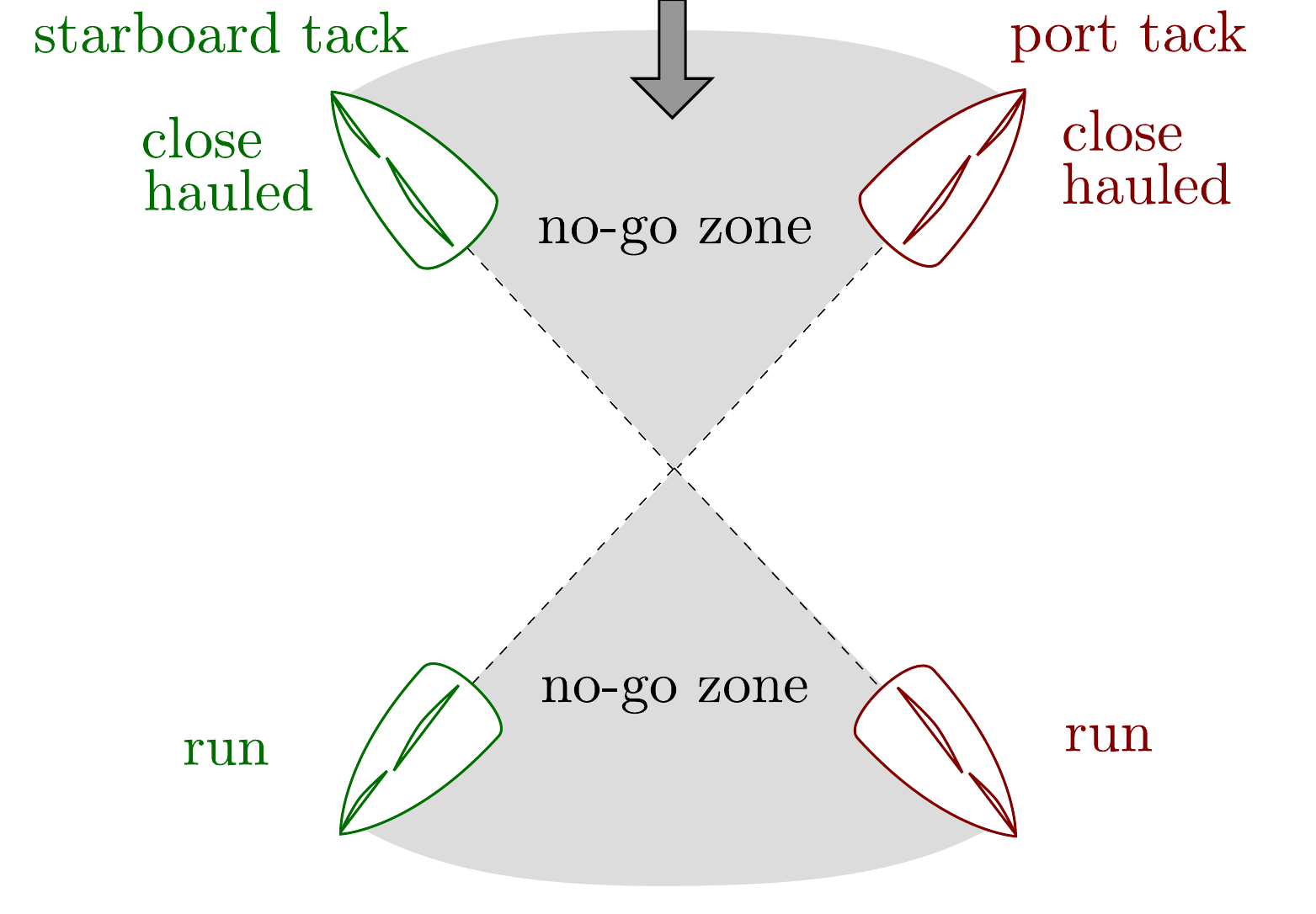}
\includegraphics[height=5cm, width=5cm,
keepaspectratio]{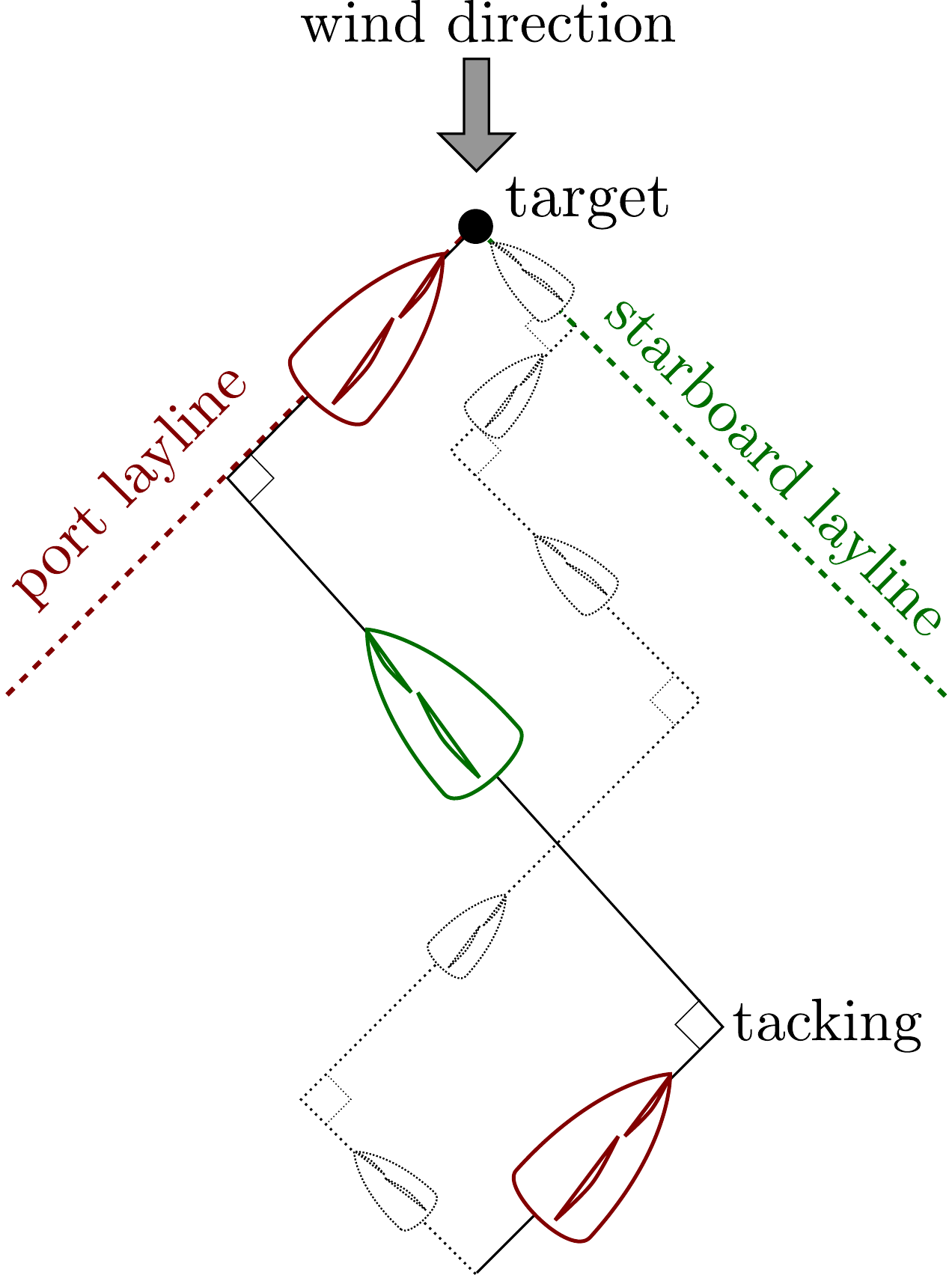}\hspace{4ex}
\end{center}
\caption{The picture on the left-hand side shows the no go zone. The picture on the right-hand side shows a tacking strategy for a boat beating upwind in a constant wind. Another boat starting at the same point and following any other  zigzag path parallel to the laylines would travel the same distance up to the target (dotted trajectory).} 
\label{Fig_points_of_sail}
\end{figure}

If  the angle between the boat's heading and the direction of the wind is too sharp, then the boat cannot move forward in this direction. The set of all  angles for which the boat cannot advance is called the \emph{no-go zone} (see Fig.~\ref{Fig_points_of_sail}). When the wind enters the sails from the right (resp.~left) side of the yacht, the yacht is sailing on \emph{starboard tack} (resp.~\emph{port tack}). A boat sailing upwind (resp.~downwind) can switch from starboard to port tack, or vice versa, by \emph{tacking} (resp.~{\em jibing}). While tacking, the yacht enters the no-go zone and hence for a few seconds, it is no longer propelled by the wind (this effect is much less pronounced when jibing). The yacht slows down and its inertia  allows it to turn to the other tack.

\medskip

\noindent{\em Sailing to a target point}
\medskip

When the boat seeks to approach a particular target point, say a buoy, then if the target is not in the boat's no-go zone, the boat can simply head directly towards it and this determines the tack on which to sail. 

However, if the target is in the yacht's no-go zone, then the yacht has to follow a zigzag path (see Fig.~\ref{Fig_points_of_sail}). This is called \emph{beating upwind}  (resp.~\emph{running downwind}): the yacht alternates between headings on starboard and on port tack.
\medskip

\noindent{\em Wind behavior}
\medskip

Both the wind direction and the wind speed vary continuously over time, and both influence the motion of the boat. However, empirical observations  indicate that at least on certain time scales, both fluctuate in an unpredictable manner, but fluctuations in wind direction are more important than those in wind speed (see \cite{DalangDumas2015}). Therefore, in our simplified model, we will assume that the wind speed is in fact constant, and only the wind direction changes over time. 

In statistical studies that led to the models developed in \cite{DalangDumas2015,DUM}, it was noticed that the assumption of centered Gaussian increments fit the available data reasonably well, and only occasionally could a better fit be obtained with non-centered increments. Therefore, in our basic model, the wind direction will be a Brownian motion on a circle with standard deviation per unit time equal to $\sigma > 0$ and without drift. Of course, this description is a major simplification, since it does not take into account the spatial variations of the wind, nor effects such as gusts and spatial correlations. But this description was in fact sufficient to be used during an actual match race (see \cite{DalangDumas2015}). 

\subsection{Model assumptions}\label{sec2.1}

\begin{figure}[h]
\begin{center}
\includegraphics[height=5cm, width=5cm,
keepaspectratio]{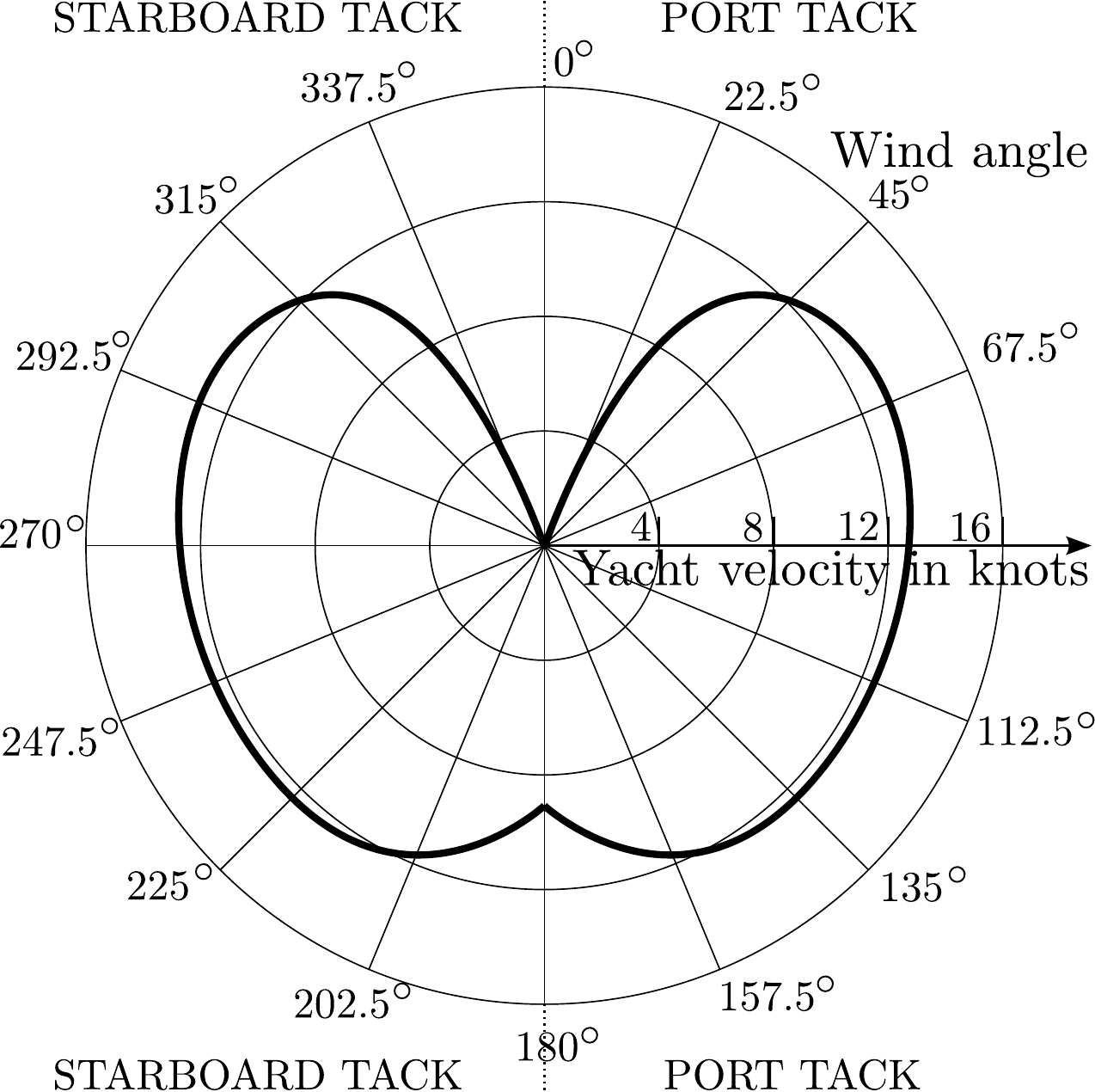}\hspace{4ex}
\includegraphics[height=5cm, width=5cm,
keepaspectratio]{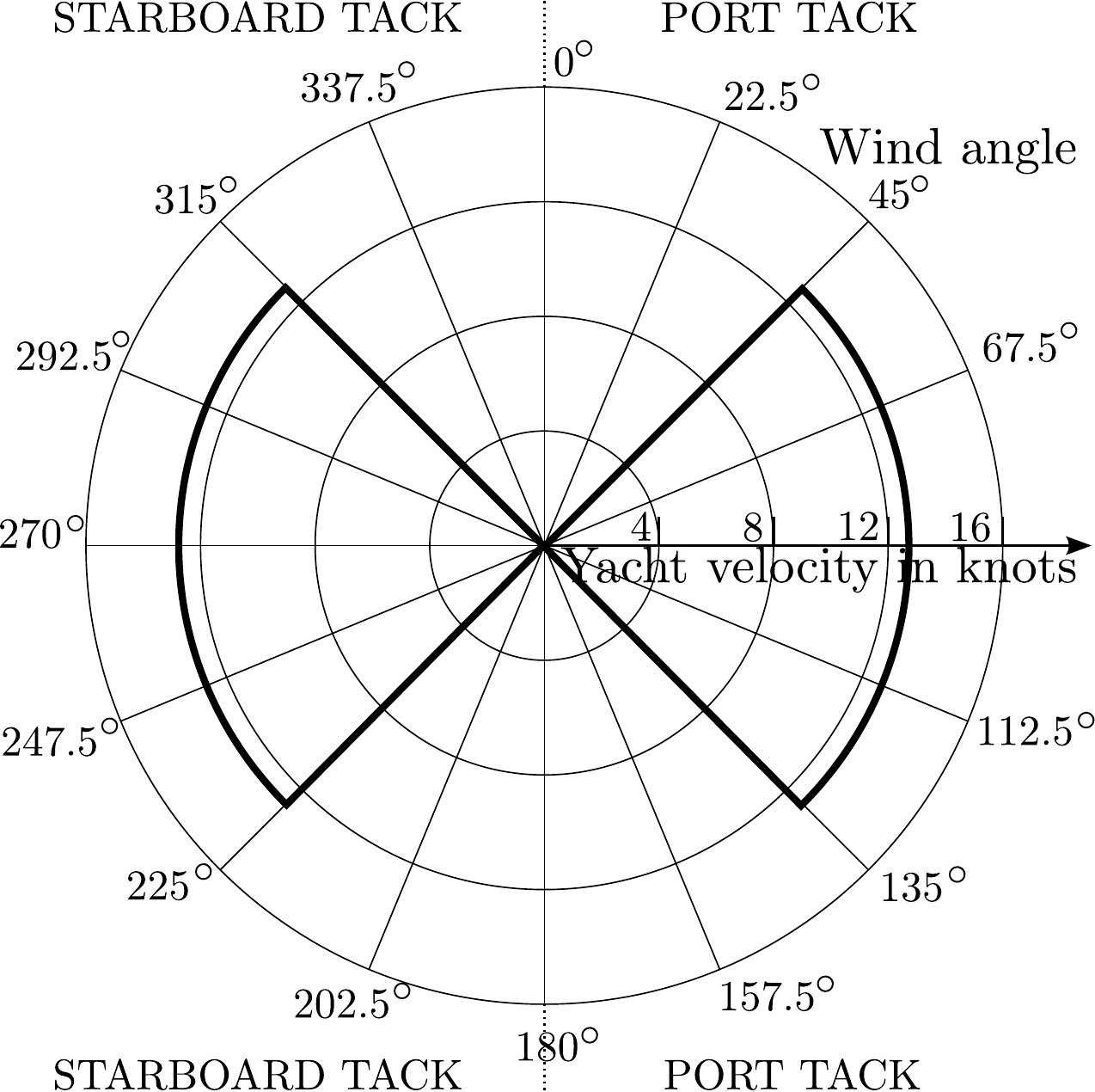}
\caption{Example of a polar diagram for a 20 knot wind speed (on the left-hand side) and our simplified {\em butterfly} model (on the right-hand side) where the optimal angle for an upwind (resp.~downwind) route is fixed at $45^o$ (resp.~$135^o$).}
\label{Fig:polar}
\end{center}
\end{figure}

For a given wind speed, the velocity of the yacht is a function of the angle between the wind direction and the yacht's heading. This function is illustrated by a {\em polar diagram} (see the left-hand side of Fig.~\ref{Fig:polar}). 
\medskip

\noindent{\em The optimal sailing angle}
\medskip

 When sailing upwind (resp.~downwind), there is one specific sailing angle that moves the boat most quickly towards the target. This is the angle in the polar diagram that maximizes the projection of the velocity vector onto the vertical axis (in Fig.~\ref{Fig:polar} left, this is  approximately $\alpha_0 = 40^o$ (resp. $\alpha_0 = 140^o$)). If the boat wants to reach as quickly as possible a target located in the no-go zone, then it should always navigate at the optimal angle $\alpha_0$ (there are certain situations in which the boat may prefer a slightly different angle, known as the Wally effect, for which we refer to \cite{DalangDumas2015}, but we do not consider this option in our model).

In our model (see Fig.~\ref{Fig:polar}, right), we consider that the yacht's behavior is symmetric regardless of whether it is beating upwind or running downwind. For simplicity, the optimal angle for an upwind (resp.~downwind) route is set to $\alpha_0 := \piq$ (resp.~$\pitq$). In between these angles, we assume that the polar curve is a butterfly-shaped curve (bounded by two quarter-circles) instead of the heart-shaped curve shown in Fig.~\ref{Fig:polar}, left. In particular, since we assume in our model that the wind speed is constant,  the yacht will always sail at a constant speed $v>0$.
\medskip

\noindent{\em The cost of tacking}
\medskip

When the yacht is beating upwind (resp.~running downwind), it will sail at angle $\alpha_0$ (resp.~$\pi - \alpha_0$) on one tack, then tack (resp.~jibe), continue at angle $\alpha_0$ (resp.~$\pi - \alpha_0$) on the other tack and then repeat this procedure until the target is reached (see Fig.~\ref{Fig_points_of_sail}). Many possible tacking strategies can be employed. Notice that if the wind direction is constant, then the number of tacks does not affect the total distance traveled. However, each tack (resp.~jibe) makes the yacht slow down and hence lose time. 
We model this loss of time by a fixed ``tacking cost" $c \geq 0$. The limiting case $c=0$ is relevant when the total travel time is very large (say several days) relative to the tacking cost (a few seconds) and it also leads to interesting mathematical questions.

   Our model as just described contains only four parameters: the wind variability $\sigma > 0$, the boat speed $v>0$, the tacking cost $c \geq 0$ and the target's radius $\eta \geq 0$.
\medskip

\noindent{\em The geographic reference frame}
\medskip

The {\em geographic reference frame} is attached to fixed geographic directions: the $\xi_1$-axis is oriented North to South and the $\xi_2$-axis is oriented West to East. We put the target at the origin and we let $\beta_t$ denote the wind direction relative to the North to South direction, measured counter-clockwise: $\beta_t = 0$ (respectively $\beta_t = \tfrac{\pi}{2}$) corresponds to the wind coming from the North (respectively the West).

The position of the boat at time $t$ is $(\xi^1_t, \xi^2_t)$. When $\xi^1_t > 0$ and $\vert \xi^2_t\vert < \xi^1_t$, the boat will beat upwind to as to move closer to the target. If it is on port tack (respectively starboard tack), then its velocity vector will be $v\, (- \cos(\beta_t - \piq), -\sin(\beta_t - \piq))$ (respectively $v\, (- \cos(\beta_t + \piq), -\sin(\beta_t + \piq))$). 
\medskip

\noindent{\em Laylines}
\medskip

The {\em laylines} are the lines $L_1$ and $L_2$ that pass through the origin/target with respective angles $\alpha_0$ and $- \alpha_0$. Because $\alpha_0$ has been set to $\piq$, these lines are orthogonal with respective equations $\xi_2 = \xi_1$ and $\xi_2 = - \xi_1$. For upwind sailing (when $\vert \xi_2 \vert < \xi_1$), $L_1$ is the {\em starboard layline} and $L_2$ is the {\em port layline}. 

\medskip

\noindent{\em A rotating reference frame}
\medskip

It turns out to be useful to attach the reference frame to the wind direction rather than to fixed geographic directions. 
In this new and rotating reference frame (centered at the origin/target), the $x$- (respectively $y$-) coordinate is measured along the port (respectively starboard) layline. The polar coordinates of the boat in this rotating frame are denoted $(r, \t)$. The radial coordinate $r$ of the boat represents its distance from the origin, while the angular coordinate $\t$ is $\piq + \alpha$, where $\alpha$ is the angle between the wind direction and the line connecting the boat to the target (equivalently, this is the angle between the port layline and this line). More precisely, $r = \sqrt{x^2 + y^2}$, and $\t = \arctan(y / x)$ if $x>0$, $\t = \pi + \arctan(y/x)$ if $x < 0$.


For example, a motionless boat would appear to move on a circle because of the changes in wind direction; and the speed vector of a boat moving upwind at speed $v$ and on port tack at angle $\alpha_0 = 45^o$ with respect to the wind will be horizontal and of the form $(-v, 0)$.

\begin{figure}[h]
\begin{center}
\includegraphics[height=7cm, width=7cm, keepaspectratio]{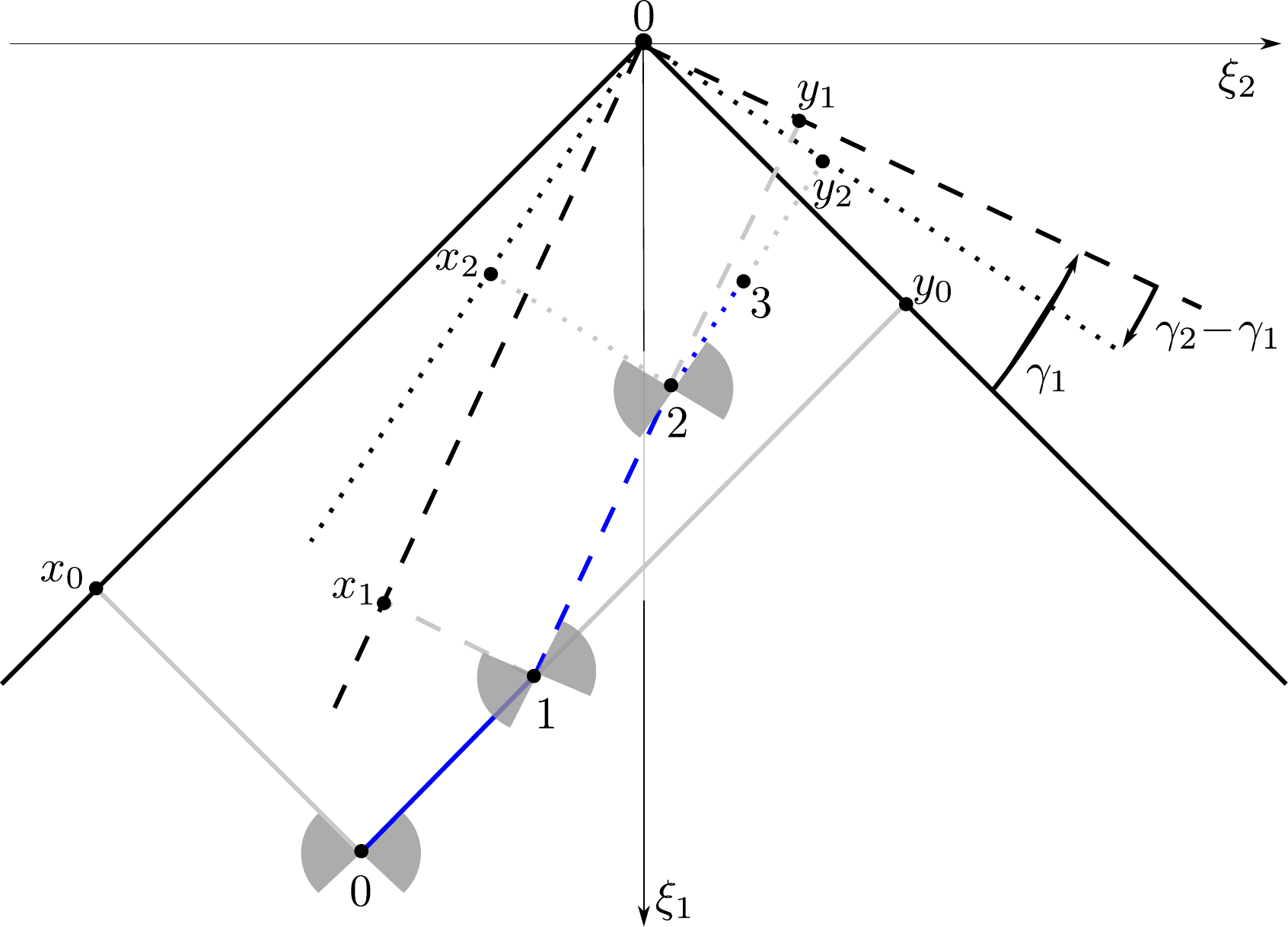}
\hspace{4ex}
\includegraphics[height=7cm, width=7cm, keepaspectratio]{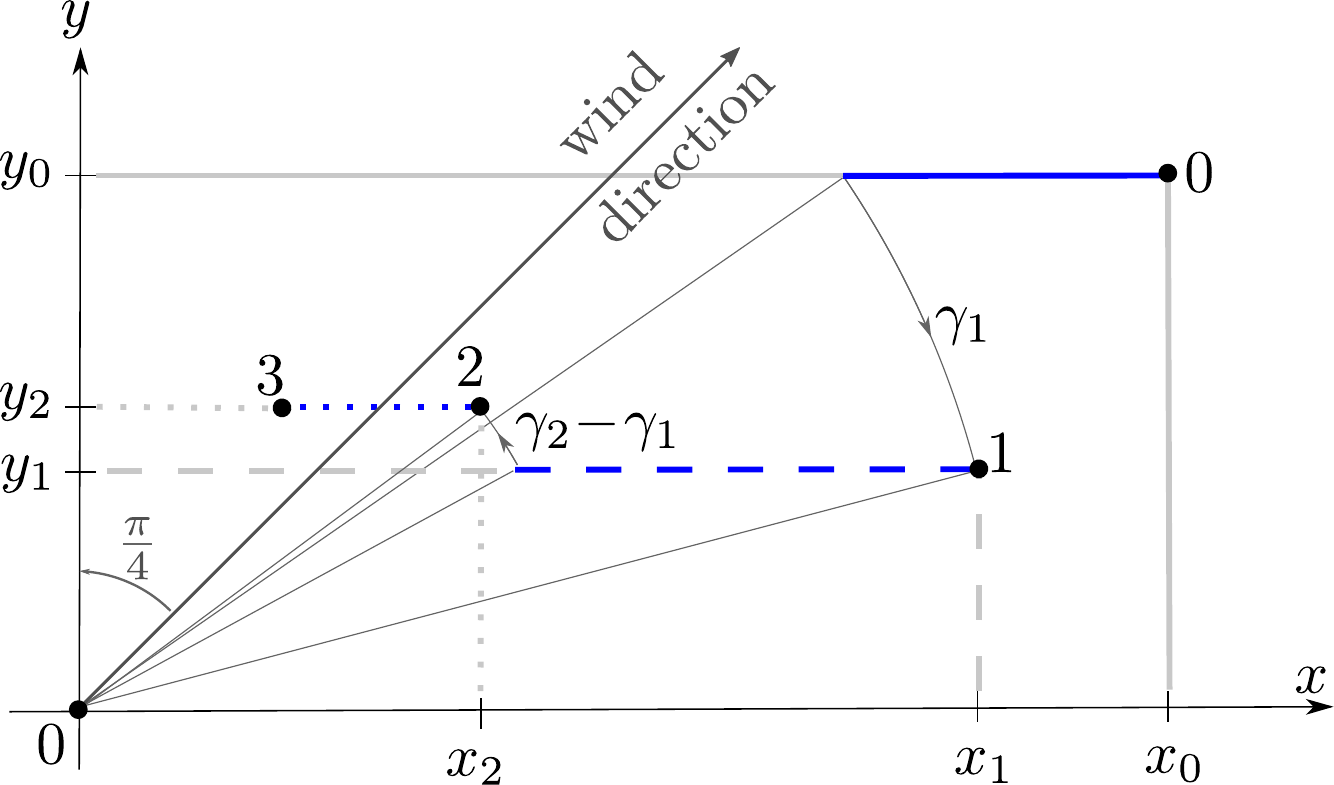}
\end{center}
\caption{The picture on the left-hand side uses the geographic reference frame to show the trajectory of a boat starting at position $0$ and sailing on port tack. The wind direction $\beta_t$ in the geographic frame is given by $\beta_t =  \gamma_1 \, 1_{[\tau_1, \tau_2[}(t) + \gamma_2\, 1_{[\tau_2, \tau_3[}(t)$, where $0 < \tau_1 < \tau_2 < \tau_3$ and $0 < \gamma_2 < \gamma_1 < \piq$. We simplify the notation by using  $x_i := x_{\tau_i}$ and $y_i := y_{\tau_i}$, $i=0, 1, 2$.  The butterfly-shaped regions correspond to the feasible directions when sailing on port tack (resp.~starboard tack). The boat follows the solid line to move from position $0$ to position $1$ between times $0$ and $\tau_1$ while the wind is coming from the North, then, when the wind direction changes to $\beta_{\tau_1} = \gamma_1$, the boat follows the dashed line to move from position $1$ to position $2$ between times $\tau_1$ and $\tau_2$, and when the wind direction changes to $\gamma_2$ (which implies a change of $\gamma_2 - \gamma_1$ from the previous direction), the boat follows the dotted line to move from position $2$ to position $3$ between times $\tau_2$ and $\tau_3$.
The picture on the right-hand side shows the corresponding trajectory of the boat in the rotating reference frame attached to the wind direction.
\label{Fig:beating:upwind}
}
\end{figure}

 To gain intuition about the problem in the rotating reference frame, consider first the case where the wind is initially coming from the North in the geographic frame, and then jumps to an angle $\gamma_ 1 \in \, ]0, \piq[$ at time $\tau_1$, then to an angle $\gamma_2 \in \,  ]0, \gamma_1[$ at time $\tau_2$, so that $\beta_t = \gamma_1 \, 1_{[\tau_1, \tau_2[}(t) + \gamma_2\, 1_{[\tau_2, \tau_3[}(t)$, where $0 < \tau_1 < \tau_2 < \tau_3$. Assume that the boat sails on port tack during the entire time interval $[0, \tau_3[$. Then during the time-interval $[0, \tau_1[$, 
 the boat's velocity vector in the geographic frame is $v\, (-\sqrt{2}/2, \sqrt{2}/2)$. Now suppose that $x_{\tau_1} > 0$ and $y_{\tau_1} > 0$. At time $\tau_1$, the wind rotates clockwise to angle $\gamma_1$, and during the time-interval $[\tau_1, \tau_2[$, the boat's 
 velocity vector in the geographic frame is $v\, (-\cos(\gamma_1 - \piq), - \sin(\gamma_1 - \piq))$. Suppose that $x_{\tau_2} > y_{\tau_2} > 0$. At time $\tau_2$, the wind rotates counterclockwise to angle $\gamma_2$, and during the time-interval $[\tau_2, \tau_3[$, the boat's 
 velocity vector in the geographic frame is $v\, (-\cos(\gamma_2 - \piq), - \sin(\gamma_2 - \piq))$. This path is shown on the left-hand side of Fig.~\ref{Fig:beating:upwind}.

In the rotating reference frame, between time $0$ and $\tau_1$, the boat's speed vector is $(-v, 0)$. At time $\tau_1$, the boat's position is subject to a counterclockwise rotation of amplitude $\gamma_1$, and during $[\tau_1, \tau_2[$, the boat's speed vector is again $(-v, 0)$. At time $\tau_2$, the boat's position is subject to a clockwise rotation of amplitude $\gamma_1 - \gamma_2$, and during $[\tau_2, \tau_3[$, the boat's speed vector is yet again $(-v, 0)$. This path is shown on the right-hand side of Fig.~\ref{Fig:beating:upwind}. Observe that the rotating reference frame provides a much simpler picture of the boat's motion.

\medskip

\noindent{\em Reducing the considered headings}
\medskip

\begin{figure}
\begin{center}
\includegraphics[height=5cm, width=5cm,
keepaspectratio]{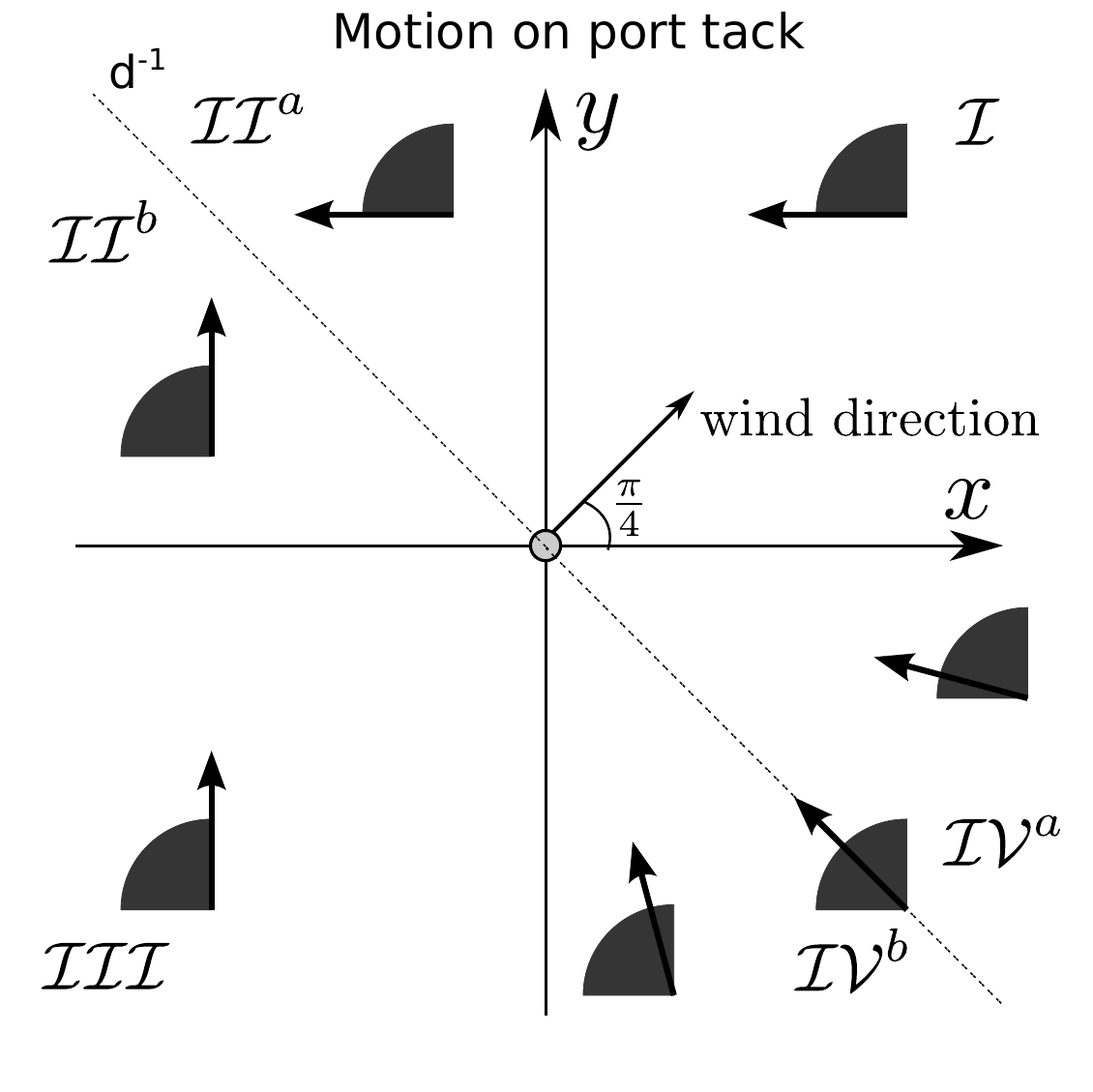}
\includegraphics[height=5cm, width=5cm,
keepaspectratio]{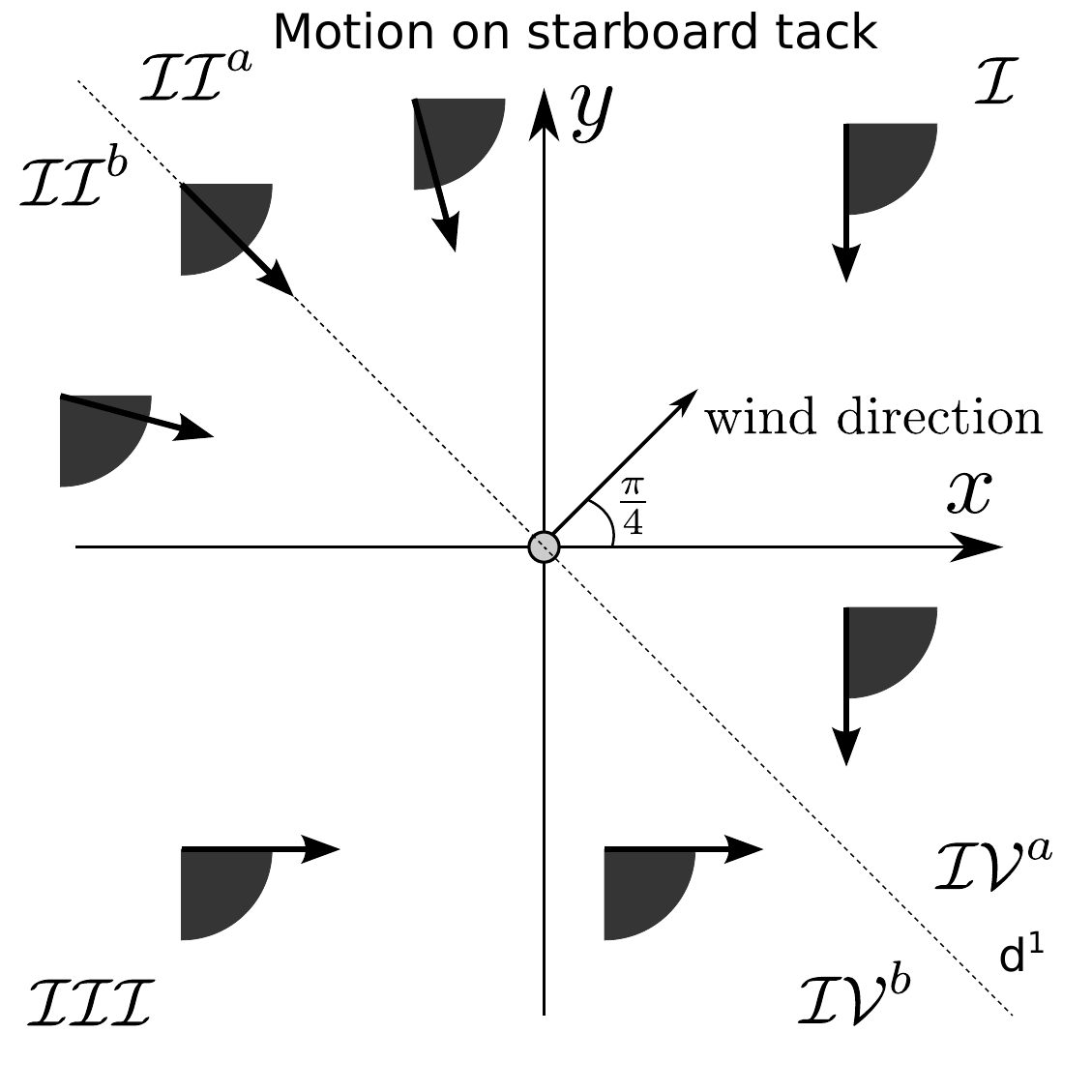}
\caption{The sailboat's motion on port tack on the left hand side and on starboard tack on the right hand side. The reference frame rotates with the wind.}
\label{Fig:4}
\end{center}
\end{figure}

For each position of the boat and each tack, the boat will usually choose the heading which moves it closest to the target. Therefore, for each position of the boat, we will only consider two headings, represented by the arrows in Fig.~\ref{Fig:4}, one for port tack and the other for starboard tack. Notice that in regions $II^a$, $II^b$, $IV^a$ and $IV^b$ shown in Fig.~\ref{Fig:4}, the boat is in fact moving {\em away} from the target; however, the boat {\em must} enter these regions and continue for at least some time there in order to reach the target with a finite number of tacks.

In the rotating reference frame, regardless of the wind direction, the motion of the boat on starboard (resp. port) tack is shown on the right- (resp.~left-) hand side of Fig.~\ref{Fig:4}. In order to prevent the boat from moving too far away from the target, we require that the yacht tack (or jibe) on the half line $d^1$ (resp $d^{-1}$) shown in Fig.~\ref{Fig:4}.


\subsection{Equations of motion}
The purpose of this section is to write the position of the yacht as a 2-dimensional It\^{o} process, as it sails in a wind whose direction is given by a Brownian motion on a circle. 

Let $\eta \geq 0$ and consider a circular target buoy $\B_\eta$ centered at the origin, which is a closed disc of radius $\eta$ if $\eta > 0$, and just a target point at the origin if $\eta = 0$.

In order to describe the motion in our model, for $\ell \in \{+1, -1\}$, we set
\begin{align}
d^\ell&=\{(x,y)\in\R^2: y = -x\mbox{ and } \ell\, x>0\} 
\label{d_p_m}
\end{align}
(see Fig.~\ref{Fig:4}). We assume that the lines $d^1$ and $d^{-1}$ are {\em switching boundaries} which cannot be crossed: when the boat encounters $d^1$ (resp.~$d^{-1}$) and is on starboard (resp.~port) tack, then the boat {\em must} tack 
(otherwise, the boat would move far away from the buoy even though it could choose to move towards it).
We also split the race field  $ \R^2 \setminus \B_\eta$ into four zones (see Fig.~\ref{Fig:4}):
\begin{align*}
\textsl{I}&=\{(x,y)\in\R^2:\, x>0 \mbox{ and } y>0\} \cap \B_\eta^c , \qquad \textsl{II}=(\textsl{II}^a\cup\textsl{II}^b\cup d^{-1}) \cap \B_\eta^c,\\
\textsl{III}&=\{(x,y)\in\R^2:\, x<0 \mbox{ and } y<0\} \cap \B_\eta^c,\qquad  \textsl{IV}=(\textsl{IV}^a\cup \textsl{IV}^b\cup d^1) \cap \B_\eta^c,
\end{align*}
where
\begin{align*}
\textsl{II}^a&=\{(x,y)\in\R^2:\, x<0 \mbox{ and } y> -x\},\qquad \textsl{II}^b=\{(x,y)\in\R^2:\, y>0 \mbox{ and } y<-x\},\\
\textsl{IV}^a&=\{(x,y)\in\R^2:\, y<0 \mbox{ and } y> -x\},\qquad \textsl{IV}^b=\{(x,y)\in\R^2:\, x>0 \mbox{ and } y<-x\}.
\end{align*}

Let $\sigma>0$ be a given real number and $(B_t)_{t\geq 0}$ be a standard real-valued Brownian motion defined on a filtered probability space $\probSpace$. In our model, $\sigma B_t$ (modulo $2\pi$) represents the wind direction at time $t$. We recall that the boat's speed is constant in our model.

Recall (see \cite[Example 5.14]{Oeksendal}) that if $\vec{m}$ and $\vec{\Sigma}$ are the two functions from $\R^2$ to $\R^2$ defined by
\begin{align*}
\vec{m}(\vec{x})=-\tfrac{\sigma^2}{2}\,\,\vec{x}\qquad\text{ and }\qquad\vec{\Sigma}(\vec{x})=\vec{\Sigma}(x,y)=\left(\begin{array}{r} -\sigma y \\ \sigma x \end{array}\right),
\end{align*}
then the unique solution $\vec{X} = (\vec{X}_t)_{t \geq 0}$ of the system of SDEs
\begin{align*}
d\vec{X}_t =\vec{m}(\vec{X}_t)\,dt+\vec{\Sigma}(\vec{X}_t)\,dB_t, \qquad \vec{X}_0=\vec{x}=(x,y),
\end{align*}
is a Brownian motion with speed $\sigma$ on the circle with radius $\|\vec{x}\| = (x^2 + y^2)^{1/2}$, given by
\begin{align*}
\vec{X}_t=\left(\cos(\sigma B_t)x-\sin(\sigma B_t)y\,,\,\sin(\sigma B_t)x+\cos(\sigma B_t)y\right).
\end{align*}
In the rotating frame, this process $(\vec{X}_t)$ defines the position of a motionless yacht.

We now define the velocity vectors that correspond to each considered sailing direction of the moving yacht. In agreement with Fig.~\ref{Fig:4},  set
\begin{align}\label{motions_def}
  \vec{v}_u =(0,v), \quad \vec{v}_d = (0,-v), \quad \vec{v}_\ell = (-v,0), \quad  \vec{v}_r =(v,0), 
\end{align}
and for $\vec{x}\in\R^2$, set $\vec{v}_\rho(\vec{x})=-v\frac{\vec{x}}{\|\vec{x}\|}$ ($u$, $d$, $\ell$, $r$ and $\rho$ stand respectively for upward, downward, leftward, rightward and radial motion).
We add these velocity vectors to the drift $\vec{m}(\vec{x})$ to obtain the drift term $\vec{\mu}(\vec{x},1)$ (resp.~$\vec{\mu}(\vec{x},-1)$) for  the starboard  (resp.~port) tack: 
\begin{align}\notag
\vec{\mu}(\vec{x},+1)&=\vec{v}_d\, \ind_{\textsl{I}\, \cup\, \textsl{IV}^a }(\vec{x})+\vec{v}_r\, \ind_{\textsl{III}\, \cup\, \textsl{IV}^b}(\vec{x})+\vec{v}_\rho(\vec{x})\, \ind_{\textsl{II}}(\vec{x})+\vec{\mu}_c(\vec{x}),\\
\vec{\mu}(\vec{x},-1)&=\vec{v}_\ell\,  \ind_{\textsl{I}\, \cup\, \textsl{II}^a}(\vec{x})+\vec{v}_u\,  \ind_{\textsl{III}\, \cup\, \textsl{II}^b}(\vec{x})+\vec{v}_\rho(\vec{x})\, \ind_{\textsl{IV}}(\vec{x})+\vec{\mu}_c(\vec{x}).
\label{mu_def_vec}
\end{align}
The position process of a yacht sailing on starboard tack, respectively on port tack, is now the solution of the following system of SDEs:
\begin{align}\label{eq_star_1}
   \text{starboard:}\qquad\qquad&d\vec{X}_t =\vec{\mu}(\vec{X}_t, +1)\,dt+\vec{\Sigma}(\vec{X}_t)\,dB_t, \qquad \vec{X}_0=\vec{x}_0, \\[6pt]
   \text{port:}\qquad\qquad\qquad\ &d\vec{X}_t =\vec{\mu}(\vec{X}_t,-1)\,dt+\vec{\Sigma}(\vec{X}_t)\,dB_t, \qquad \vec{X}_0=\vec{x}_0.
\label{eq_port}
\end{align}
On starboard tack and in regions $I$, $III$ and $IV$, and on port tack and in regions $I$, $II$ and $III$, these equations are of the form
\begin{align}\label{eq_sol_position_process}
d\vec{X}_t =(\vec{\mu}_c(\vec{X}_t)+\vec{v})\,dt+ \vec{\Sigma}(\vec{X}_t)\,dB_t, \qquad \vec{X}_0=\vec{x}, \quad t \geq 0, 
\end{align}
where $\vec{v}$ is one of the four vectors in \eqref{motions_def}.

Proposition \ref{prop_brw_SDE_solut} below exhibits the solution of these SDEs for any constant velocity vector $\vec{v} = (v_1, v_2)$, and in particular for the four velocity vectors defined in \eqref{motions_def}, while Proposition \ref{prop_brw_SDE_solut_radial} below gives a solution for the radial motion (starboard tack in region $II$ or port tack in region $IV$).

\begin{prop}\label{prop_brw_SDE_solut}
Let $\vec{x}=(x, y)\in\R^2 \setminus \B_\eta$ and $\vec{v} = (v_1,v_2)\in\R^2$ be the initial position of the boat and the velocity vector selected by the boat, respectively. The process $(\vec{X}_t=(X_t, Y_t))_{t\geq 0}$, taking values in $\R^2$ and defined by

\begin{equation}\label{eq_position_process}
\vec{X}_t=\left(\begin{array}{c}
\displaystyle\cos(\sigma B_t)x-\sin(\sigma B_t)y + v_1\int_0^t\cos(\sigma B_t-\sigma B_s)ds - v_2\int_0^t\sin(\sigma B_t-\sigma B_s)ds\\
\displaystyle\sin(\sigma B_t)x+\cos(\sigma B_t)y + v_1\int_0^t\sin(\sigma B_t-\sigma B_s)ds + v_2\int_0^t\cos(\sigma B_t-\sigma B_s)ds
\end{array}
\right)
\end{equation}
is the unique solution of \eqref{eq_sol_position_process}.
\end{prop}

\begin{proof}
Uniqueness follows directly from \cite[Thm.~5.2.1]{Oeksendal}, since the drift in \eqref{eq_sol_position_process}, as well as the diffusion coefficient there, is affine. Thus, it only remains to verify that the process $\vec{X}$ defined by \eqref{eq_position_process} is indeed a solution of \eqref{eq_sol_position_process}, by using It\^{o}'s formula. For details see \cite[Prop. 5.1]{LV}.
\end{proof}

\begin{prop}\label{prop_brw_SDE_solut_radial}
Let  $\eta > 0$ and $\vec{x}=(x, y)\in\R^2\setminus \B_\eta$. The process $(\vec{X}_t,\,t\in\left[0,(\|\vec{x}\| - \eta)/v \right])$ taking values in $\R^2 \cap (\text{\rm int }\B_\eta)^c$ defined by
\begin{equation}\label{eq_process_radial}
\vec{X}_t=\left(\begin{array}{c}
X_t\\
Y_t
\end{array}
\right)
=\left(\begin{array}{c}
\left(\|\vec{x}\|-vt\right)\cos\left(\sigma B_t+\theta\right)\\
\left(\|\vec{x}\|-vt\right)\sin\left(\sigma B_t+\theta\right)
\end{array}
\right)
\end{equation}
where $\theta\in\R$ is such that $\cos(\theta)=\frac{x}{\|\vec{x}\|}$ and $\sin(\theta)=\frac{y}{\|\vec{x}\|}$, is the unique solution of
\begin{equation}\label{eq_sde_system_radial}
d\vec{X}_t =(\vec{\mu}_c(\vec{X}_t)+\vec{\mu}_\rho({\vec{X}_t}))\,dt+ \vec{\Sigma}(\vec{X}_t)\,dB_t, \qquad \vec{X}_0=\vec{x}, \quad t\in\left[0,\|\vec{x}\| / v\right].
\end{equation}
\end{prop}

\begin{proof}
Observe that the drift coefficient in \eqref{eq_sde_system_radial} is  Lipschitz continuous in $\R^2 \cap \B_\eta^c$, hence the existence and uniqueness of the solution holds on $t\in \left[0, (\|\vec{x}\| - \eta)/v \right]$. 
It remains to verify that the process $\vec{X}$ defined by \eqref{eq_process_radial} is indeed a solution of \eqref{eq_sde_system_radial}, using It\^o's formula. For details see \cite[Prop. 5.2]{LV}.
\end{proof}

Notice that Proposition \ref{prop_brw_SDE_solut} and \eqref{eq_sol_position_process} give us the system of SDEs and the corresponding solutions for the following four motions: 
(1) {\em downward motion} ($v_1=0, v_2=-v$) when the boat is on starboard tack in zone $\textsl{I}\cup \textsl{IV}^a$; (2)  {\em upward motion} ($v_1=0, v_2=v$) when the boat is on port tack in zone $\textsl{II}^b\cup \textsl{III}$; {\em leftward motion}  ($v_1=-v, v_2=0$) when the boat is on port tack in zone $\textsl{I}\cup\textsl{II}^a$; (4) {\em rightward motion}  ($v_1=v, v_2=0$) when the boat is on starboard tack in zone $\textsl{III}\cup \textsl{IV}^b$.
Therefore, Propositions \ref{prop_brw_SDE_solut} and \ref{prop_brw_SDE_solut_radial} provide explicit formulas to describe the solution of \eqref{eq_star_1}   and \eqref{eq_port}.

\begin{rem}\label{strong_cart}
Denote $(X^1_t,Y^1_t)$ (resp. $(X^{-1}_t,Y^{-1}_t)$)  the process defined by \eqref{eq_star_1} (resp. \eqref{eq_port}).
Define $\tau^a_f:= \inf \{ t: (X^a_t,Y^a_t) \in d^{a} \cup \B_\eta\}$, $a\in \{-1,+1 \} $, where the half-lines $d^{1}$  and $d^{-1}$ are defined in \eqref{d_p_m}.
Equation \eqref{eq_star_1} (resp.~\eqref{eq_port}) has Lipschitz continuous coefficients in $\R^2\setminus d^a$, 
therefore by a standard localization argument, there exists a strong solution until time $\tau^a_f$, 
for any initial condition $\bar{x}_0=(x_0,y_0)\in \R^2 \setminus d^a$. 
\end{rem}
\medskip

\noindent{\em Motion in polar coordinates}
\medskip

\begin{figure}
\begin{center}
\includegraphics[width=5cm]{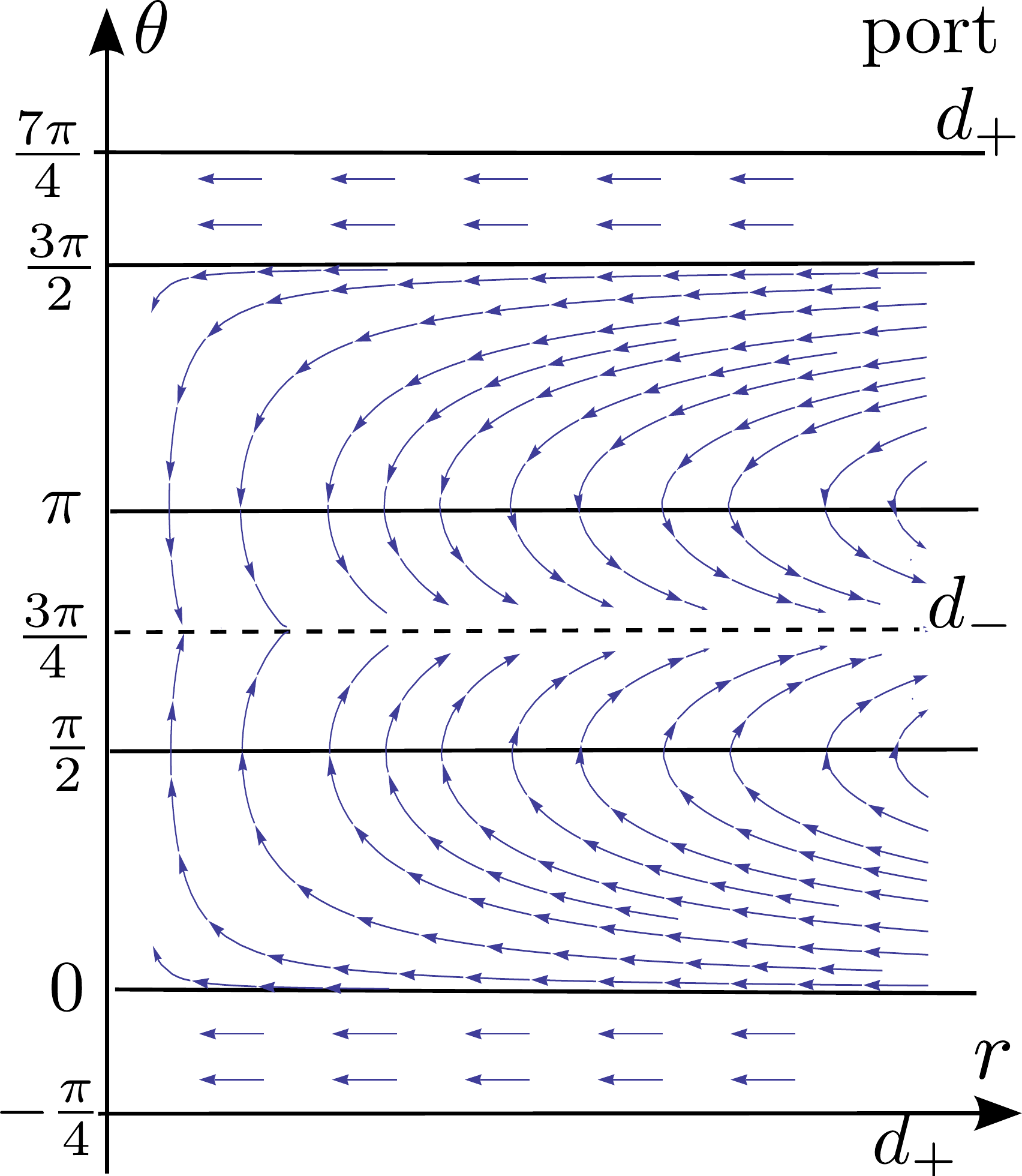}\hspace{5ex}\includegraphics[width=5cm]{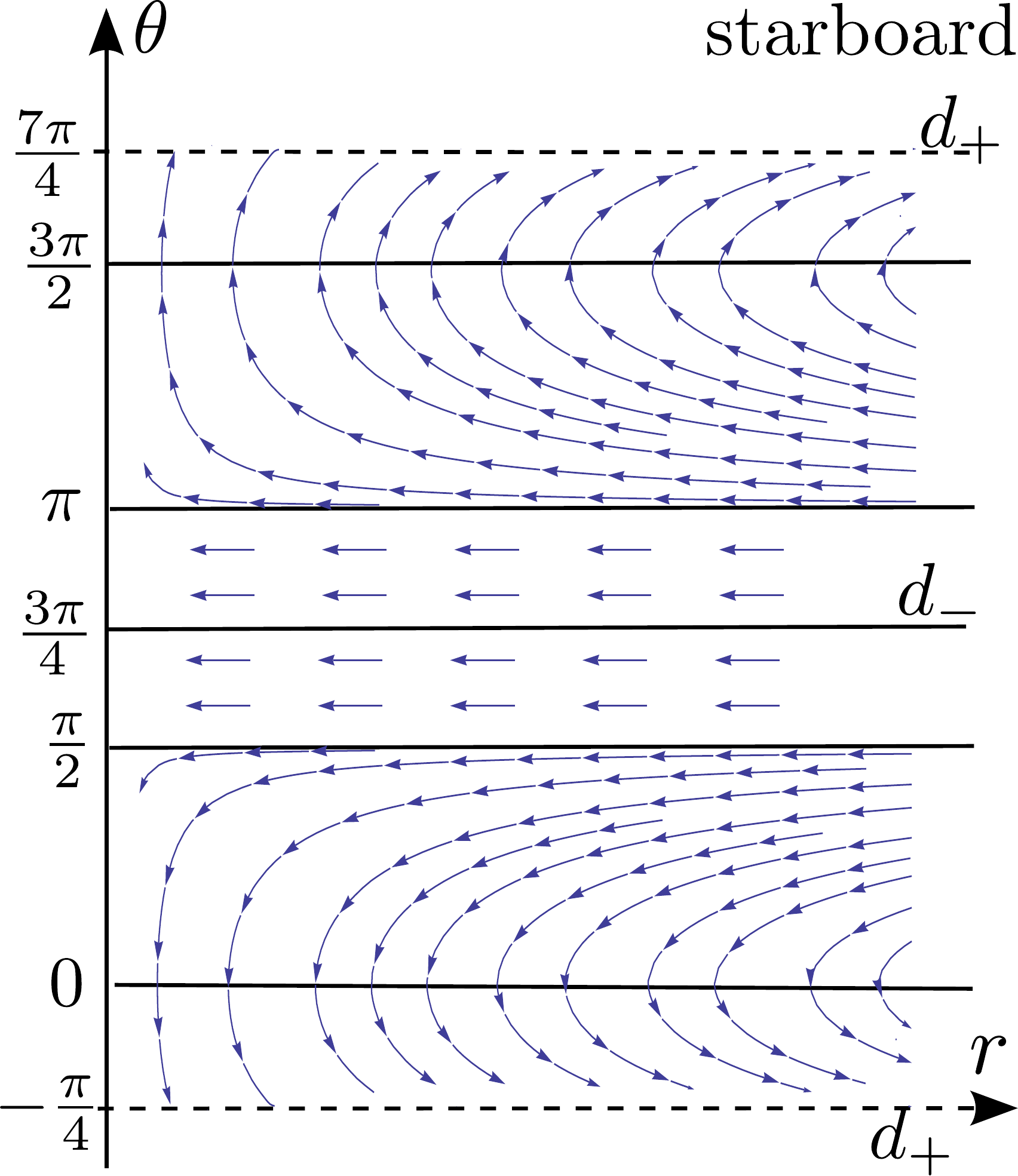}
\end{center}
\caption{Illustration,  in polar coordinates, of the motion of the boat in a constant wind ($\sigma$ set to $0$ in \eqref{eds_proc_controle_polar}), which corresponds to the flow induced by the drift $(\mu_1(\t, a), \mu_2(r,\t, a))$. On the left-hand (respectively right-hand) side, the yacht is on port tack (respectively starboard tack). In our fluctuating wind, in addition to the motion along the flow, the yacht moves as a Brownian motion in a direction parallel to the $\theta$-axis, with a speed that depends on the diffusion coefficient $\sigma$. Note that there is no Brownian component in the radial direction.}
\label{Fig_brw_sailing_rules_polar}
\end{figure}

It turns out that the equations of motion become more compact when written in polar coordinates. For definiteness, we consider $2 \pi$-periodic functions defined on $\left[-\frac{\pi}{4},\frac{7\pi}{4}\right[$ and then extend them to $\R$ by $2 \pi$-periodicity. Consider, for instance, the process describing a yacht sailing upwind on starboard tack: this correponds to $v_1=0$ and $v_2=-v$ in \eqref{eq_sol_position_process} and  \eqref{eq_position_process}. Therefore, the position process $(\vec{X}_t = (X_t, Y_t))$ satisfies the system of SDEs
\begin{equation*}
\left\{
\begin{array}{rll}
dX_t &=-\frac{\sigma^2}{2}\,X_t\,dt-\sigma Y_t\,dB_t,                   &\quad X_0=x,\\
dY_t &=\left(-\frac{\sigma^2}{2}\,Y_t-v\right)\,dt+\sigma X_t\,dB_t,    &\quad Y_0=y.
\end{array}
\right.
\end{equation*}
Applying the multidimensional It\^{o}'s formula \cite[Section 4.2]{Oeksendal} to 
\begin{equation*}
    R_t=\sqrt{X_t^2+Y_t^2},  \qquad \Theta_t=\arctan(Y_t/X_t), 
\end{equation*}
we find that
\begin{align*}
R_t &=R_0 +\int_0^t\frac{X_s}{R_s}\left(-\tfrac{\sigma^2}{2}\,X_s\,ds -\sigma Y_s\,dB_s \right)
+\int_0^t\frac{Y_s}{R_s}\left(\left(-\tfrac{\sigma^2}{2}\,Y_s-v\right)\,ds +\sigma X_s\,dB_s\right) \\
 &  \quad  +\tfrac{1}{2}\int_0^t \frac{R_s^2-X_s^2}{R_s^3}Y_s^2 \sigma^2\,ds + \int_0^t \frac{- X_sY_s}{R_s^3}(-X_sY_s)\sigma^2\,ds
    +\tfrac{1}{2}\int_0^t\frac{R_s^2-Y_s^2}{R_s^3}X_s^2\sigma^2\,ds\\
    &=R_0 -\int_0^tv\frac{Y_s}{R_s}\,ds = R_0 -\int_0^tv\sin(\Theta_s)\,ds
\end{align*}
after several simplifications, and similarly,
\begin{align*}
\Theta_t &=\Theta_0 +\int_0^t-\frac{Y_s}{R_s^2}\left(-\tfrac{\sigma^2}{2}\,X_s\,ds -\sigma Y_s\,dB_s \right)
+\int_0^t\frac{X_s}{R_s^2}\left(\left(-\tfrac{\sigma^2}{2}\,Y_s-v\right)\,ds +\sigma X_s\,dB_s\right) \\
&\quad +\tfrac{1}{2}\int_0^t \frac{2X_sY_s}{R_s^4}Y_s^2 \sigma^2\,ds + \int_0^t \frac{-X_s^2+Y_s^2}{R_s^4}(-X_sY_s)\sigma^2\,ds
 +\tfrac{1}{2}\int_0^t\frac{-2X_sY_s}{R_s^4}X_s^2\sigma^2\,ds\\
    &= \Theta_0 -\int_0^t\frac{v\cos(\Theta_s)}{R_s}\,ds+\int_0^t\sigma\,dB_s.
\end{align*}
In other words, given a starting position $(x, y) \in \R^2\setminus \B_\eta$, the position process $(R_t, \Theta_t)$ for a yacht sailing upwind on starboard tack satisfies
\begin{equation*}
\left\{
\begin{array}{rll}
dR_t        &=-v\sin(\Theta_t)\,dt,                                 &\hspace{5ex} R_0=\|\vec{x}\|,\\
d\Theta_t   &=-\dfrac{v\cos(\Theta_t)}{R_t}\,dt+\sigma\,dB_t,       &\hspace{5ex}\Theta_0=\arctan(y/x).
\end{array}
\right.
\end{equation*}
We do the same for the other four sailing directions and we find that for $a\in \{+1, -1 \}$, \emph{the position process $(R^a_t, \Theta^a_t)_{t\geq0}$ of the yacht in polar coordinates}, is the solution of
\begin{equation}\label{eds_proc_controle_polar}
\left\{\begin{aligned}
dR^a_t      &= 1_{\{R^a_t > \eta\}}\, \mu_1(\Theta^a_t,a)\,dt, \\                                             
d\Theta^a_t &= 1_{\{R^a_t > \eta\}} \left[\mu_2(R^a_t, \Theta^a_t,a)\,dt+ \sigma\,dB_t\right],                          
\end{aligned}\right.
\end{equation}
where $\vec{\mu}(r,\theta,a)=(\mu_1(\theta,a),\mu_2(r,\theta,a))$ is defined, for the motion of the boat on \textit{starboard tack} and $\theta \in [-\tfrac{\pi}{4},\tfrac{7\pi}{4}[$,  by
\begin{equation}
 \mu_1(\t, +1)= \left\{\begin{array}{l}
-v \sin(\theta), 
\\
-v, 
\\
\ \ v\cos(\theta), 
\end{array}
\right.
\quad \mu_2(r,\t, +1)= \left\{\begin{array}{cl}
- \frac{v \cos(\theta)}{r}, &\qquad \theta \in [-\frac{\pi}{4}, \frac{\pi}{2} [,\\
0, &\qquad \theta \in [\frac{\pi}{2}, \pi [,\\
- \frac{v\sin(\theta)}{r}, &\qquad \theta \in [\pi, \frac{7 \pi}{4} [,
\end{array}
\right.
\label{eq_drift_ra}
\end{equation}
while for the motion of the boat on \textit{port tack} and $\t \in [-\piq,\pisq[$, they are defined by
\begin{equation}
\mu_1(\t,-1)= \left\{\begin{array}{l}
-v, 
\\
-v \cos(\theta), 
\\
\ \ v\sin(\theta), 
\end{array}
\right.
\ \mu_2(r,\t,-1)= \left\{\begin{array}{cl}
0, &\qquad \theta \in  [-\piq,0 [ \,\cup\, [ \pitm, \pisq [,\\
\frac{v \sin(\theta)}{r}, &\qquad \theta \in [0, \frac{3\pi}{4} [,\\
\frac{v\cos(\theta)}{r}, &\qquad \theta \in [\frac{3\pi}{4}, \frac{3 \pi}{2} [.
\end{array}
\right.
\label{eq_drift_pa}
\end{equation}
Therefore, after extending these definitions to $\t \in \R$ by $2\pi$-periodicity, the mapping $(\theta,a) \mapsto \mu_1(\theta,a)$ (resp.~$(r,\theta,a) \mapsto \mu_2(r,\theta,a)$) is defined on  $\R \times \{ +1,-1\}$ with values in  $[-v,v]$ (resp.~on $\R_+^*\times \R\times \{+1,-1 \}$ with values in $\R$).
See Fig.~\ref{Fig_brw_sailing_rules_polar} for an illustration of the drift $(\mu_1(\t, a), \mu_2(r,\t, a))$ in polar coordinates. Notice that $R^a_t \geq \eta$ and $\T^a_t \in \R$.

\begin{rem}\label{strong_pol}
Even though the coefficients $\mu_2(r,\theta,a)$ may diverge as $r \downarrow 0$,   the system \eqref{eds_proc_controle_polar} admits a strong solution to up to time $\tau^a_f$ (defined in Remark \ref{strong_cart}) because of the corresponding result in Cartesian coordinates (see  Remark \ref{strong_cart}). At time $\tau^a_f$, if $\B_\eta$ has not been reached, then $a$ is replaced by $-a$ (meaning the boat tacks), and the solution can be continued and this procedure repeated, either until $\B_\eta$ is reached or for all time if this never occurs.
\end{rem}

\section{Formulating two control problems}\label{Sec:3}
In this section, we introduce two control problems, one with $c= 0$ and the other with $c>0$ and we present the Hamilton-Jacobi-Bellman (HJB) equation for each. 

Fix $(r, \t) \in\, ]\eta, \infty[ \times \R$ and suppose that at time $t \geq 0$, the boat is on tack $A_t \in \{+1, -1\}$, where $A = (A_t)_{t \geq 0}$ is a process chosen by the skipper (a strategy). In view of  \eqref{eds_proc_controle_polar}, the position process $(R^A_t , \Theta^A_t)$ of the boat at times $t \geq 0$ is given by the solution of the system of SDEs
\begin{equation}\label{state1}
\left\{\begin{array}{l l}
dR^A_t = 1_{\{R_t^A > \eta \}}\, \mu_1(\Theta^A_t, A_t) dt, & R^A_0 = r, \\
d\Theta^A_t = 1_{\{R_t^A > \eta \}} \left[\mu_2(R^A_t, \Theta^A_t,A_t) dt +   \sigma  dB_t \right], & \Theta^A_0 = \t,
\end{array}
\right.
\end{equation}
and  the time $\tau^{A}_\eta$ to reach the target buoy $\B_\eta$ is
 \begin{equation}\label{tau_n}
   \tau^{A}_\eta:= \inf\left\{t: R^{A}_t \leq  \eta  \right\}.
 \end{equation}
 Of course, the process $(A_t)_{t\geq 0}$ will have to satisfy certain properties, so as to be an \textit{admissible strategy}: these properties depend on whether $c=0$ or $c>0$  and are stated in Definitions \ref{Def:adm_con} and \ref{Def:3.3} below, respectively.
 
The position of the yacht at time $t$ when driven by $A$ is defined as the $\P_{r,\t, a}$-strong solution $(R_t^A,\T_t^A)_{t \in [0,   \tau^{A}_\eta]}$ of \eqref{state1},
where $\P_{r,\t,a}$ denotes a family of probability measures on path space such that
\begin{equation}\label{p_def}
\P_{r,\t,a}\left\{ (R^A_0,\T^A_0, A_{0-}) = (r, \t,a) \right\} = 1.
\end{equation}
The initial tack $a$ is the yacht's tack just before the race starts, so $A_{0}$ may differ from $A_{0-}$ if the yacht tacks immediately at the start of the race.

\subsection{Case where $c=0$}\label{Subs:1.2.1}
When $c=0$, a strategy that allows an infinite number of tacks may lead the boat to $\B_\eta$ in finite time, and the initial tack $a$ of the boat has no bearing on the problem, hence we will omit it in \eqref{p_def}, and write $\P_{r,\t}$ instead of $\P_{r,\t,a}$.

\begin{dfn}\label{Def:adm_con}
 Let $\eta \geq 0$, $(r,\theta) \in S^0 := \, ]\eta, \infty[ \times \R$ and $A: \R_+ \times \Omega \mapsto \{+1, -1\}$ be a jointly measurable and $(\shf_t)$-adapted  process.  We say that $A$ is an \textit{admissible strategy} if:
 
(1)  the system \eqref{state1} admits a strong solution;

(2)  $\E_{r,\t} \left( \tau^{A}_\eta \right) < \infty$.
\end{dfn}

The set of admissible strategies under $\P_{r,\t}$ is denoted by $\sha^\eta_{r,\t}$ (sometimes the subscript $\cdot_{r,\t}$ will be omitted).  In Section \ref{Sec:4}, we shall prove  that for every $\eta \geq 0$, the set $\sha^\eta$ is non-empty (see Theorem \ref{thm_brw_exist_V}).  

The control problem $\shp(\eta)$ consists in finding an \textit{admissible strategy} that minimizes the expected time to reach the upwind buoy:  for $A\in \sha^\eta$, let
\begin{align}\label{Jeta}
J^{\eta}(r,\theta,A)= \E_{r,\t} \left(\tau^{A}_\eta \right);
\end{align}
the \textit{value function} is defined as
\begin{align}\label{Veta}
V^{\eta}(r,\theta) = \inf_{A \in \sha^\eta} J^{\eta}(r,\t,A).
\end{align}

In the above definition, $V^{\eta}: [\eta, \infty[ \times \R \to \R_+$ and $J^{\eta}: [\eta, \infty[ \times \R \times \sha^\eta \to \R_+ $,
hence the control problem consists in finding an admissible strategy $A^{\star} \in \sha^\eta$ that is \textit{optimal,} that is, $V^{\eta}(r,\theta) = J^{\eta}(r,\t,A^{\star})$.
In Section \ref{Sec:5}, we will present a candidate strategy that will turn out to be optimal for this control problem (see 
\cite{CC1}).
\medskip

\noindent{\em The HJB equation}
\medskip

In this case where $c=0$, our control problem can be interpreted in the general framework for stochastic control problems outlined for instance in \cite{Oeksendal,pham}. In this general framework, the running cost is simply the constant function $1$ (cost $1$ per unit time), and there is no terminal cost.  

For $\psi \in C^{1,2}(S^0)$ and $a \in \{1, -1\}$, let
\begin{align*}
    \shl^a \psi(r, \t) := \mu_1(\theta,a)\,\frac{\partial}{\partial r}\psi(r,\t) 
  +\mu_2(r,\theta,a)\, \frac{\partial}{\partial \theta}\psi(r,\t) + \frac{\sigma^2}{2}\, \frac{\partial^2}{\partial \theta^2}\psi(r,\t).
\end{align*}
According to \cite[Theorem 11.2.1]{Oeksendal} or \cite[Section 3.4.1]{pham}, the HJB equation for the value function $V = V^{\eta}$ is 
\begin{equation}
\min_{a\in\{-1,+1\}}
\left[
\mathcal{L}^a V(r,\theta)+1
\right]
=0,
\qquad
(r,\theta)\in S^0,
\label{eq:HJB-c-zero}
\end{equation}
with the initial condition
\[
V(\eta,\theta)=0, \qquad \theta \in \mathbb{R},
\]
and we seek a solution that is $2\pi$-periodic in the variable $\theta$. This HJB equation will be confirmed in Theorem \ref{rd07_01t1}.

As it is the case in most stochastic control problems, we do not expect the value function $V^{\eta}$ to be sufficiently regular  that the differential operator can be applied to it everywhere in $S^0$: see also the case $c>0$ below.

\subsection{Case where $c>0$}\label{rd07_15ss1}
When the tacking cost is positive, the yacht is penalized by $c$ units of time for each tack,  therefore a useful strategy cannot allow an infinite number of tacks, since these would induce an infinite total for the tacking penalties. In particular, the strategy should be piecewise constant in the following sense.

\begin{dfn}
Let $A: (\{0-\} \cup \R_+)\times \Omega \ra \R^n$ be an $(\shf_t)_{t\geq 0}$ adapted and right-continuous process. If, for all $(r, \theta, a) \in \, ]\eta, \infty[ \times \R \times \{+1, -1\}$, for $\P_{r, \theta, a}$-almost all $\omega \in \Omega$  and $T <+ \infty$, there exists $\ve > 0$ such that for all $t \in [0,T]$, if $A_t \neq A_{t-}$, then  $A_{t+r}=A_t$ for $r\in [0,\ve]$,  then we say that $A$ is a \textit{piecewise constant process}.
The \textit{jump} of $A_t$ at $t \geq 0$ is defined as $\Delta A_t = A_t - A_{t-}$.
\end{dfn}

The main difference between Definition \ref{Def:3.3} below and 
Definition \ref{Def:adm_con} 
is the restriction to piecewise constant processes and the condition on integrability of the number of tacks, formulated in Item (d) below.

\begin{dfn}\label{Def:3.3}
Let $\eta \geq 0$, $c>0$, and $(r ,\t, a)\in S^> := \, ]\eta, \infty[ \times \R\times \{+1,-1 \}$. Let $A: (\{0-\} \cup \R_+) \times \Omega \mapsto \{+1, -1\}$. Then $A$  is an \textit{admissible strategy} if the following properties hold: 

(a) $A$ is a piecewise constant process;

(b) there exists  a $\P_{r,\t,a}$-strong solution to \eqref{state1}; 

(c) $\E_{r,\t,a} \left( \tau^{A}_\eta \right) < \infty$, where $\tau^{A}_\eta$ is defined in \eqref{tau_n};

(d) $\E_{r,\t,a} \left(N_{\tau^{A}_\eta }(A) \right) < \infty$, where 
$ 
N_{t}(A) :=  \sum_{s=0}^t \vert \Delta A_s \vert. 
$ 
\end{dfn}

The set of the admissible strategies  under $\P_{r,\t,a}$  is denoted $\sha^{c, \eta}_{r,\t,a}$  (sometimes, the subscript $\cdot_{r,\t,a}$ will be omitted).  In Section \ref{Sec:4}, we will show that $\sha^{c,\eta}$ is non-empty.


The optimal control problem consists in finding an \textit{admissible strategy} that minimizes the expected {\em penalized time} to reach the target, which accounts for the penalties due to the tacks. The {\em payoff} for strategy $A$ is therefore $\tau^A_\eta + c\, N_{\tau^A_\eta}(A)$, and letting
\begin{align}\label{Jetac}
   J^{c,\eta}(r,\theta,a,A)= \E_{r,\t,a} \left(\tau^A_\eta + c\, N_{\tau^A_\eta}(A)\right),
\end{align}
the \textit{value function} is defined by
\begin{align}\label{Vetac}
    V^{c,\eta}(r,\theta,a) = \inf_{A \in \sha^{c,\eta}_{r,\t,a}} J^{\eta}(r,\t,a,A).
\end{align}

In the above definition, $V^{c,\eta}$ and $J^{c,\eta}$ are defined respectively on $\R_+ \times \R \times \{-1,+1 \}$ and $\R_+ \times \R \times \{-1,+1 \} \times \sha^{c,\eta}$, and the control problem $\shp(c,\eta)$ consists in finding an  \textit{admissible strategy} $A^\star \in \sha^{c,\eta}_{r,\t,a}$ such that  $V^{c,\eta}(r,\theta,a) = J^{c,\eta}(r,\t,a,A^\star)$.

Note that for $\eta \geq 0$ and $c>0$, $\sha^{c,\eta}  \subset \sha^\eta $, where $\sha^\eta$ is introduced just after  Definition \ref{Def:adm_con}. In particular,  $V^\eta(r, \theta) \leq V^{c, \eta}(r, \theta, a)$ always holds.
\medskip

\noindent{\em HJB equation}
\medskip

In this case where $c>0$, our control problem can be interpreted in the framework of optimal switching problems considered in \cite[Chapter 6]{OS}.  The running cost is again the constant function $1$ (cost $1$ per unit time), and there is no terminal cost.  

For $\varphi \in C^{1,2}(S^>)$, let 
\begin{align*}
    \shl \varphi(r, \t, a) := \mu_1(\theta,a)\,\frac{\partial}{\partial r}\varphi(r,\t,a) 
  +\mu_2(r,\theta,a)\, \frac{\partial}{\partial \theta}\varphi(r,\t,a) + \frac{\sigma^2}{2}\, \frac{\partial^2}{\partial \theta^2}\varphi(r,\t,a).
\end{align*}
According to \cite[Chapter 6]{OS}, the HJB equation for the value function $V = V^{c,\eta}$ is
\begin{equation}
\min\left[
\mathcal{L} V(r,\theta,a)+1,\, V(r,\theta,-a)+c-V(r,\theta,a)
\right] =0, \qquad (r,\theta,a) \in S^>,
\label{eq:QVI}
\end{equation}
with the initial condition
\[
V(\eta,\theta)=0, \qquad \theta\in\mathbb{R},
\]
and we seek a solution that is $2\pi$-periodic in the variable $\theta$. This HJB equation will be confirmed in Theorem \ref{rd06_29t1}.

As in the case $c=0$, we do not expect the value function $V^{\eta}$ to be sufficiently regular that the differential operator can be applied to it everywhere in $S^>$. We have not succeeded in solving this equation, but we indicate in Appendix B the anticipated shape of the solution. 

\section{Verification theorems}

The theorem below is a so-called {\em verification theorem}: let $S^> =\, ]\eta, \infty[ \times \R \times \{-1,+1 \}$. Given a candidate value function $\Vcand:  \bar S^> \to \R_+$ and $A^{\star} \in \sha^{c,\eta}_{r,\t,a}$ that is a candidate for an optimal strategy, it provides a sufficient condition for $\Vcand$ to be the value function of the problem $\shp(c, \eta)$ and for the strategy $A^{\star}$ to be optimal. 

If $c = 0$,  notations are understood related to the control problem of Section \ref{Subs:1.2.1}, that is, $S^0 =\, ]\eta, \infty[ \times \R$,  $\sha^{0,\eta}_{r,\t,a} =\sha^{\eta}_{r,\t}$ (recall that the initial tack is irrelevant when $c=0$) and  $\shp(0, \eta) =  \shp(\eta)$.

\begin{thm}[Verification theorem]\label{Ver}
Fix $\eta\ge0$ and $c\ge0$. We let $S$ denote $S^>$ if $c > 0$ and $S^0$ if $c=0$. Suppose there is a continuous function
$\Vcand: \bar S \to\R_+$ such that
\begin{itemize}
  \item[(a)] $\Vcand(\eta,\theta,a)= 0$,  for all $(\theta, a) \in \R \times \{ -1,+1 \}$;
  \item[(b)] $\Vcand$ is continuous at $r=\eta$, uniformly in $\theta\in\R$ and $a  \in  \{ -1,+1 \}$;
  \item[(c)] there are $b,d\in\R_+$ with $\Vcand(r,\theta,a)\le b\, r+d$ for all $(r,\theta,a)\in [\eta,\infty[\times\R \times \{ -1,+1\}$.
\end{itemize}
Fix $(r,\theta,a)\in S$.  For $A\in\sha^{c,\eta}_{r,\theta,a}$ set
\begin{equation}\label{eq:N}
  H(t,r,\theta,a,n):=t+c\, n+ \Vcand(r,\t,a),\qquad
  M^{A}_t:=H\left(t,R^A_t,\T^A_t,A_t,N_t(A)\right),
\end{equation}
$t\in \{0- \} \cup \R_+$, where $(R^A,\Theta^A)$ is the  process introduced in \eqref{state1}. 
\begin{enumerate}
  \item[(1)] If, for every $A\in\sha^{c,\eta}_{r,\theta,a}$, the process
  $\bigl(M^{A}_{t \wedge \tau^A_\eta}\bigr)_{t\ge0}$ is a $\P_{r,\theta,a}$-submartingale,
  then
  \[
    \Vcand(r,\theta,a)\le V^{c,\eta}(r,\theta,a);
  \]
  \item[(2)] if, in addition, there is $A^\star\in\sha^{c,\eta}_{r,\theta,a}$ for which
  $\bigl(M^{A^*}_{t \wedge \tau^{A^*}_\eta}\bigr)_{t\ge0}$ is a $\P_{r,\theta,a}$-martingale, then
  \[
    \Vcand(r,\theta,a)= V^{c,\eta}(r,\theta,a)\qquad\text{and }\quad A^\star\text{ is optimal.}
  \]
\end{enumerate}
\end{thm}

\begin{proof}
(1.) Let $A\in\sha^{c,\eta}_{r,\theta,a}$ be an admissible strategy. Since $t\wedge\tau^A_\eta$ is a
bounded stopping time, the submartingale property gives
\begin{align}\notag
  \Vcand(r,\theta,a) &=  M^{A}_{0-} \le\E_{r,\theta,a}\bigl[M^{A}_{t \wedge \tau^A_\eta}\bigr] \\
   &=\E_{r,\theta,a}\Bigl[t\wedge\tau^A_\eta+c\, N_{t\wedge\tau^A_\eta}(A)
  +\Vcand\bigl(R^A_{t\wedge\tau^A_\eta},\Theta^A_{t\wedge\tau^A_\eta}, A_{t\wedge\tau^A_\eta}\bigr)\Bigr].
\label{34}
\end{align}
Let $t\to\infty$. As $A$ is admissible, $\tau^A_\eta<\infty$ $\P_{r,\theta,a}$-a.s., so
$t\wedge\tau^A_\eta\uparrow\tau^A_\eta$ and
$N_{t\wedge\tau^A_\eta}(A)\uparrow N_{\tau^A_\eta}(A)$; by monotone convergence,
\begin{equation}\label{34a}
  \E_{r,\theta,a}\bigl[t\wedge\tau^A_\eta\bigr]\uparrow\E_{r,\theta,a}\bigl[\tau^A_\eta\bigr],
  \qquad
  \E_{r,\theta,a}\bigl[N_{t\wedge\tau^A_\eta}(A)\bigr]\uparrow\E_{r,\theta,a}\bigl[N_{\tau^A_\eta}(A)\bigr].
\end{equation}
For the last term, because $\tau^{A}_\eta < \infty$, $\P_{r,\theta,a}$-a.s., we have $R^A_{t\wedge \tau^{A}_\eta} \to R^A_{\tau^{A}_\eta} = \eta$ as $t\rightarrow  \infty$, $\P_{r,\theta, a}$-a.s. By  assumption (b), and since $\Vcand(\eta,\t,a) = 0$ by (a), we have $\Vcand(R^A_{t\wedge \tau^{A}_\eta},\Theta^A_{t\wedge\tau^{A}_\eta}, A_{t\wedge\tau^A_\eta}) \to 0$ as $t\rightarrow \infty$.

By \eqref{eq_drift_ra} and \eqref{eq_drift_pa}, $| \mu_1(\theta, a) | \leq v$, for all $\theta \in \R$ and $a \in \{+1, -1\}$. It follows from the first equation in \eqref{state1} that
\begin{equation}\label{rd07_01e1}
|R^A_{t\wedge \tau^{A}_\eta} -r| \leq (t\wedge \tau^{A}_\eta) v \leq v\tau^{A}_\eta, \qquad\text{therefore }\quad 
 R^A_{t\wedge \tau^{A}_\eta} \leq v\, \tau^{A}_\eta + r .
\end{equation}
By assumption (c),
\begin{align*}
\Vcand(R^A_{t\wedge \tau^{A}_\eta},\Theta^A_{t\wedge\tau^{A}_\eta},A_{\tau^{A}_\eta}) \leq b R_{t\wedge \tau^{A}_\eta} + d \leq bv\, \tau^{A}_\eta + br+d .
\end{align*}
Since $\E_{r,\theta,a} \left(\tau^{A}_\eta \right) < \infty$, we can apply the dominated convergence theorem to see that
\begin{equation}\label{36}
\E_{r,\theta,a} \left(\Vcand(R^A_{t\wedge \tau^{A}_\eta},\Theta^A_{t\wedge\tau^{A}_\eta},A_{\tau^{A}_\eta})\right) \rightarrow  0 \qquad \text{as } t\rightarrow \infty.
\end{equation}

From \eqref{34}--\eqref{36}, and using that, by construction, $c\, N_{\tau^A_\eta}(A)$ is exactly the total cost of tackings, we conclude that
$ 
\Vcand(r,\theta,a) \leq \E_{r,\theta,a}\left(\tau^{A}_\eta  +c\, N_{\tau^A_\eta}(A)\right) = J^{\eta}(r,\theta, a, A) .
$ 
Since this holds for all $A\in \sha^{c,\eta}_{r,\theta,a}$,  by taking the infimum over $A$, we obtain
\begin{align}\label{37}
\Vcand(r,\theta,a) \leq V^{c,\eta}(r,\theta,a),
\end{align}
and (1) is proved.

(2) Under the strategy $A^\star$, by the martingale property assumed here, \eqref{34} holds with the inequality there replaced by an equality, hence
\begin{equation*}
\Vcand(r,\theta,a)= \E_{r,\theta,a} \left(t \wedge \tau^{A^\star}_\eta + +c\, N_{t\wedge\tau^{A^\star}_\eta}(A^\star) + \Vcand(R^{A^\star}_{t \wedge \tau^{A^\star}_\eta},\Theta^{A^\star}_{t \wedge \tau^{A^\star}_\eta}, A^\star_{t \wedge \tau^{A^\star}_\eta} )\right).
\end{equation*}
By the same argument as in part (1), this converges to $\E_{r,\theta,a}\left( \tau^{A^\star}_\eta + c\, N_{t\wedge\tau^{A^\star}_\eta}(A^\star)\right)$ as $t \to \infty$.
This implies that
\begin{align*}
\Vcand(r,\theta,a) = \E_{r,\theta,a}\left(\tau^{A^\star}_\eta + c\, N_{t\wedge\tau^{A^\star}_\eta}(A^\star) \right) \geq \inf_{A\in \sha^\eta} \E_{r,\theta,a}\left(\tau^{A}_\eta  + c\, N_{t\wedge\tau^{A}_\eta}(A)\right)= V^{c,\eta}(r ,\theta,a),
\end{align*}
and by \eqref{37}, the conclusion in (2) follows.
\end{proof}

We can use the previous theorem to give verification theorems that are based on the Hamilton-Jacobi-Bellman equations \eqref{eq:QVI} when $c>0$, or \eqref{eq:HJB-c-zero} when $c=0$. We begin with the case $c>0$. 

For $\varphi \in C^{1,2}(S^>)$, let 
\begin{align*}
    \shl \varphi(r, \t, a) := \mu_1(\theta,a)\,\frac{\partial}{\partial r}\varphi(r,\t,a) 
  +\mu_2(r,\theta,a)\, \frac{\partial}{\partial \theta}\varphi(r,\t,a) + \frac{\sigma^2}{2}\, \frac{\partial^2}{\partial \theta^2}\varphi(r,\t,a).
\end{align*}

\begin{thm}[HJB verification theorem, case $c > 0$] Fix  $\eta \geq 0$, $c > 0$.

\noindent (a) Suppose that there is a function $\Vcand \in C(\bar S^>) \cap C^1(S^>)$ such that:

   (i) For fixed $r > \eta$ and $a \in \{+1, -1\}$, $\theta \mapsto \Vcand(r, \t,a)$ is periodic with period $2\pi$;

   (ii) $\Vcand(r,\t,a) \leq \Vcand(r,\t,-a) +c$, for all $(r,\t,a) \in \bar S^>$;
   
   (iii) define 
\begin{align*}
   C = \{(r,\t,a) \in S^>: \Vcand(r,\t,a) < \Vcand(r,\t,-a) + c \}\qquad \text{(the no-tacking region).} 
\end{align*}
We assume that $D:= \partial C$ is a finite union of Lipschitz curves in $S^>$, $\Vcand \in C^2(S^> \setminus D)$, and $\Vcand$ has locally bounded first and second derivatives near $D$;
   
   (iv)   For $(r,\t,a) \in S^> \setminus  D$, 
\begin{align}\label{rd06_29e1}
    \shl \Vcand(r,\t,a) + 1 \geq 0;
\end{align}

   (v) for all $(r,\t,a) \in S^>$, for all $A \in \sha^{c,\eta}_{r,\t,a}$, 
\begin{align*}   
   \E_{r,\t,a}\left[\int_0^{\tau^A_\eta} 1_{D}(R^A_s,\Theta^A_s, A_s)\, ds \right] = 0.
\end{align*}
   
   (vi) Conditions (a), (b) and (c) of Theorem \ref{Ver} are satisfied.
   
\noindent Then for all $(r,\t,a) \in \bar S^>$, $\Vcand(r,\t,a) \leq V^{c,\eta}(r,\theta,a) $.
\medskip

\noindent(b) Suppose, in addition to (i)--(vi), that:

    (vii)  For $(r,\t,a) \in C$, $\shl \Vcand(r,\t,a) + 1 = 0$;
    
     (viii) there is an admissible strategy $A^*$ such that for all $(r,\t,a) \in S^>$, $\P_{r,\t,a)}$-a.s., for all $s \in [0, \tau^{A^*}_\eta]$, $A^*_s \neq A^*_{s-}$ if and only if $(R^{A^*}_s, \Theta^{A^*}_s, A^*_s) \in S^> \setminus C$.
     
     (ix) Define $S^>_\pm = \{(s,\t,a) \in S^>: a = \pm 1\}$ and $C_\pm = C \cap S^>_\pm$. Then $(s,\t,a) \in S^>_+ \setminus C_+$ implies $(s,\t,-a) \in C_-$; and $(s,\t,a) \in S^>_- \setminus C_-$ implies $(s,\t,-a) \in C_+$.

\noindent Then for all $(r,\t,a) \in \bar S^>$, $\Vcand(r,\t,a) = V^{c,\eta}(r,\theta,a)$ and the strategy $A^*$ is optimal.
\label{rd06_29t1}
\end{thm}

\begin{proof}
(a) Let $A$ be an admissible strategy. By Hypothesis (iii) and \cite[Appendix D, Theorem D.1]{Oeksendal}, there is a sequence $(\Vcand_j) \subset C(\bar S^>) \cap  C^{2}(S^>)$ such that  $\Vcand_j \to \Vcand $ (respectively $\shl \Vcand_j \to \shl \Vcand$) uniformly on compact subsets of $\bar S^>$ (respectively uniformly on compact subsets of $S^> \setminus D$), and the sequence $(\shl \Vcand_j)$ is locally bounded on $S^>$.  By It\^o's formula,
\begin{align*} 
   \Vcand_j(R^A_t,\Theta^A_t, A_t) &= \Vcand_j(R^A_0,\Theta^A_0, A_{0-}) + \int_0^t \frac{\partial\Vcand_j}{\partial r} (R^A_s,\Theta^A_s, A_s) \, dR^A_s \\  
   &\qquad + \int_0^t  \frac{\partial\Vcand_j}{\partial \t} (R^A_s,\Theta^A_s, A_s) \, d\Theta^A_s + \frac12 \int_0^t  \frac{\partial^2 \Vcand_j}{\partial \theta^2}(R^A_s,\Theta^A_s, A_s)\, d\langle \Theta^A_\cdot \rangle_s \\ \notag
   &\qquad + \sum_{s \in [0, t]} \left(\Vcand_j(R^A_s,\Theta^A_s, A_s) - \Vcand_j(R^A_s,\Theta^A_s, A_{s-}) \right). 
\end{align*}
By Hypothesis (v), we have
  \begin{align}  \notag
  \Vcand_j(R^A_t,\Theta^A_t, A_t)  &= \Vcand_j(R^A_0,\Theta^A_0, A_{0-}) + \int_0^t  \shl \Vcand_j(R^A_s,\Theta^A_s, A_s) \, 1_{S^> \setminus D}(R^A_s,\Theta^A_s, A_s) \, ds \\ \notag
   &\qquad +  \sigma \int_0^t  \frac{\partial\Vcand_j}{\partial \t} (R^A_s,\Theta^A_s, A_s)\, 1_{S^> \setminus D}(R^A_s,\Theta^A_s, A_s) \, dB_s \\
   &\qquad + \sum_{s \in [0, t]} \left(\Vcand_j(R^A_s,\Theta^A_s, A_s) - \Vcand_j(R^A_s,\Theta^A_s, A_{s-}) \right).
\label{rd06_29e2}
\end{align}
Define $H$ and $M^A_t$ as in \eqref{eq:N}, and let
\begin{align*}
   H_j(t,r,\theta,a,n):=t+c\, n+ \Vcand_j(r,\t,a),\qquad
  M^{A}_{j,t}:=H_j\left(t,R^A_t,\T^A_t,A_t,N_t(A)\right).
\end{align*}
By \eqref{rd06_29e2},
  \begin{align}  \notag
M^{A}_{j,t}
 &= \Vcand_j(R^A_0,\Theta^A_0, A_{0-}) + \int_0^t  (1 + \shl \Vcand_j(R^A_s,\Theta^A_s, A_s)) \, 1_{S^> \setminus D}(R^A_s,\Theta^A_s, A_s) \, ds \\ \notag
    &\qquad +  \sigma \int_0^t  \frac{\partial\Vcand_j}{\partial \t} (R^A_s,\Theta^A_s, A_s)\, 1_{S^> \setminus D}(R^A_s,\Theta^A_s, A_s)  \, dB_s \\
    & \qquad + c\, N_t(A) +  \sum_{s \in [0, t]} \left(\Vcand_j(R^A_s,\Theta^A_s, A_s) - \Vcand_j(R^A_s,\Theta^A_s, A_{s-}) \right).
\label{rd06_30e1}
\end{align}

   We have observed in \eqref{rd07_01e1} that $R^A_t \leq vt + r$, therefore the radial component $R^A_s$ is uniformly bounded (in $\Omega$). By (i) and the ``locally bounded'' and ``uniform convergence on compact sets'' properties of $\shl\Vcand_j$ and $\Vcand_j$, we can let $j \to \infty$ in \eqref{rd06_30e1}  to find that
  \begin{align}  \notag
M^{A}_{t \wedge \tau^A}
 &= \Vcand_j(R^A_0,\Theta^A_0, A_{0-}) + \int_0^{t \wedge \tau^A}  (1 + \shl \Vcand(R^A_s,\Theta^A_s, A_s)) \, 1_{S^> \setminus D}(R^A_s,\Theta^A_s, A_s) \, ds \\ \notag
    &\qquad +  \sigma \int_0^{t \wedge \tau^A}  \frac{\partial\Vcand}{\partial \t} (R^A_s,\Theta^A_s, A_s)\, 1_{S^> \setminus D}(R^A_s,\Theta^A_s, A_s) \, dB_s \\
    & \qquad + c\, N_{t \wedge \tau^A}(A) +  \sum_{s \in [0, t \wedge \tau^A]} \left(\Vcand(R^A_s,\Theta^A_s, A_s) - \Vcand(R^A_s,\Theta^A_s, A_{s-}) \right).
\label{rd06_30e2}
\end{align}   
   
By Hypothesis (ii),  the sum of the last two terms in \eqref{rd06_30e2} is nonnegative, and by (iv), the integrand in the $ds$-integral is also nonnegative. Therefore, $(M^{A}_{t \wedge \tau^A})$ is a submartingale, and by (vi), the conclusion in part (a) follows from Theorem \ref{Ver} (1).

(b) We consider the equality \eqref{rd06_30e2} with $A$ replaced by $A^*$ from (viii). By (viii),
each term in the sum over $s \in [0, t]$ vanishes unless $(R^{A^*}_s, \Theta^{A^*}_s, A^*_s) \in \bar S^> \setminus C$, in which case the corresponding term is equal to $-c$, by the definition of $C$ and by Hypothesis (ii). Therefore, the sum of the last two terms in \eqref{rd06_30e2} vanishes. 

By (vii), the integrand in the $ds$-integral in \eqref{rd06_30e2} vanishes when $(R^{A^*}_s, \Theta^{A^*}_s, A^*_s) \in C \cup \partial C$. We now check that $\P_{r,\t,a}$-a.s., $(R^{A^*}_s, \Theta^{A^*}_s, A^*_s)$ belongs to the open set $C$ except at the discrete set of times where $A^*_s \neq A^*_{s-}$. 

Assume that $(r,\t,a) \in C$. By (viii), $(R^{A^*}_s, \Theta^{A^*}_s, A^*_s)$ remains in $C$ during $[0, \tau_1[$, where $\tau_1> 0$ is the first time $\tau_1$ where $(R^{A^*}_s, \Theta^{A^*}_s, A^*_s) \in D$. At this time, $A^*_{s-}$ switches to $A^*_{s} = - A^*_{s-}$. By (ix), $(R^{A^*}_{\tau_1}, \Theta^{A^*}_{\tau_1}, A^*_{\tau_1}) \in C$. Again by (viii), $(R^{A^*}_s, \Theta^{A^*}_s, A^*_s)$ remains in $C$ during $[\tau_1, \tau_2[$, where $\tau_2> \tau_1$ is the first time $\tau_2$ where $(R^{A^*}_s, \Theta^{A^*}_s, A^*_s) \in D$. The tack switches again, and the yacht is again in $C$. Proceeding inductively, we see that $(R^{A^*}_s, \Theta^{A^*}_s, A^*_s)$ belongs to $C$ except at the discrete set of times where $A^*_s \neq A^*_{s-}$

It follows that the $ds$-integral in \eqref{rd06_30e2} vanishes. Since the only remaining term is the stochastic integral, $(M^{A}_{t \wedge \tau^A})$ is a martingale, and the conclusion  in part (b) follows from Theorem \ref{Ver} (2). 
\end{proof}

\begin{thm}[HJB verification theorem, case $c = 0$] Fix  $\eta \geq 0$, $c = 0$. For $\psi \in C^{1,2}(S^0)$ and $a \in \{1, -1\}$, let
\begin{align*}
    \shl^a \psi(r, \t) := \mu_1(\theta,a)\,\frac{\partial}{\partial r}\psi(r,\t) 
  +\mu_2(r,\theta,a)\, \frac{\partial}{\partial \theta}\psi(r,\t) + \frac{\sigma^2}{2}\, \frac{\partial^2}{\partial \theta^2}\psi(r,\t).
\end{align*}

\noindent (a) Suppose that there is a function $\Vcand \in C(\bar S^0) \cap C^1(S^0)$ such that:

   (i) For fixed $r > \eta$, $\theta \mapsto \Vcand(r, \t)$ is periodic with period $2\pi$;
   
   (ii) Assume that $S^0 = C_{1} \cup C_{-1} \cup D$, where $D = \partial C_{1} = \partial C_{-1}$, and $C_{1}$ and $C_{-1}$ are disjoint open sets.
We assume that $D$ is a finite union of Lipschitz curves in $S^0$, $\Vcand \in C^2(S^0 \setminus D)$, and $\Vcand$ has locally bounded first and second derivatives near $D$;
   
   (iii)   For  $a \in \{1, -1\}$ and $(r,\t) \in C_a$, 
\begin{align}\label{rd07_01e2}
    \shl^a \Vcand(r,\t) + 1 \geq 0;
\end{align}

   (iv) for all $(r,\t) \in S^0$, for all $A \in \sha^{\eta}_{r,\t}$, 
\begin{align*}   
   \E_{r,\t,a}\left[\int_0^{\tau^A_\eta} 1_{D}(R^A_s,\Theta^A_s)\, ds \right] = 0.
\end{align*}
   
   (v) Conditions (a), (b) and (c) of Theorem \ref{Ver} are satisfied.
   
\noindent Then for all $(r,\t) \in \bar S^0$, $\Vcand(r,\t) \leq V^{\eta}(r,\theta) $.
\medskip

\noindent(b) Suppose, in addition to (i)--(v), that:

    (vi)  For $(r,\t) \in C_1$, $\shl^1 \Vcand(r,\t) + 1 = 0$, and for $(r,\t) \in C_{-1}$, $\shl^{-1} \Vcand(r,\t) + 1 = 0$;
    
     (vii) there is an admissible strategy $A^*$ such that for all $(r,\t,a) \in S^0$, $\P_{r,\t,a)}$-a.s., $A^*_s = +1$ if $(R^{A^*}_s, \Theta^{A^*}_s) \in C_1$, and $A^*_s = -1$ if $(R^{A^*}_s, \Theta^{A^*}_s) \in C_{-1}$;
   
\noindent Then for all $(r,\t) \in \bar S^0$, $\Vcand(r,\t) = V^{\eta}(r,\theta)$ and the strategy $A^*$ is optimal.
\label{rd07_01t1}
\end{thm}

\begin{proof}
The main difference with the proof of Theorem \ref{rd06_29t1} is that $\Vcand$ no longer depends on the variable $a$. Using Hypotheses (i), (ii) and (iv), \eqref{rd06_30e2} is replaced by
  \begin{align}  \notag
M^{A}_{t \wedge \tau^A}
 &= \Vcand_j(R^A_0,\Theta^A_0) + \int_0^{t \wedge \tau^A}  \left[(1 + \shl^{1} \Vcand(R^A_s,\Theta^A_s)) \, 1_{C_{1}}(R^A_s,\Theta^A_s) \right. \\ \notag
 &\qquad\qquad\qquad\qquad\qquad\qquad  \left. + (1 + \shl^{-1} \Vcand(R^A_s,\Theta^A_s)) \, 1_{C_{-1}}(R^A_s,\Theta^A_s) \right] \, ds \\ 
    &\qquad +  \sigma \int_0^{t \wedge \tau^A}  \frac{\partial\Vcand}{\partial \t} (R^A_s,\Theta^A_s)\, 1_{S^0 \setminus \partial C}(R^A_s,\Theta^A_s) \, dB_s . 
\label{rd06_30e3}
\end{align}   
The  submartingale  property of $(M^{A}_{t \wedge \tau^A})$ follows from (iii) and (v), and this proves part (a).

   For part (b), let $A^*$ be an admissible strategy as in (vii). By Hypotheses (vi) and (vii), the $ds$-integral in \eqref{rd06_30e3} (with $A$ there replaced by $A^*$) vanishes, therefore $(M^{A}_{t \wedge \tau^A})$ is a martingale and the conclusion follows as in the proof of Theorem \ref{rd06_29t1}.
\end{proof}

\section{Upper bound on the value function}\label{Sec:4}
 
 In this section, we establish that the value functions $V^\eta(r, \theta)$ and $V^{c, \eta}(r, \theta, a)$ defined in \eqref{Veta} and \eqref{Vetac}, respectively, are finite. We will do this by analyzing a particular strategy $A$ in $\sha^{c,\eta}$ which is in fact close to what a typical skipper would use when approaching the target.

The main difficulty is to prove that when using this strategy, the target will be reached after
a finite number of tacks. Indeed, in a highly variable wind ($\sigma$ large), the boat may repeatedly be pushed away from the target while in Region $IV$ and on starboard tack, or in Region $II$ and on port tack (see Figs.~\ref{Fig:4} and \ref{Fig_brw_sailing_rules_polar}, and Appendix A for examples of wind variations that prevent the boat from reaching the target). The main idea will be to show that when the wind angle is given by a Brownian motion on a circle, then such situations arise with probability $0$: no matter the variability of the wind, the boat can maneuver so as to remain within a fixed distance of the target (see Remark \ref{rd10_11r1}), and then, when the wind calms down for a sufficiently long period of time, the target will be reached with only a few tacks (Corollary \ref{cor_exist_taui_taua}) and in a bounded amount of time (Corollary \ref{cor_exist_E1_tauA}); finally, since such calm periods will eventually occur, and in fact after a number of tacks related to a geometric distribution (Proposition \ref{prop_exist_p0}), the expected penalized time to reach the target will be finite (Proposition \ref{prop_NtauA_tauA_finite}).


Theorem \ref{thm_brw_exist_V}, which provides an upper bound on $V^\eta(r, \theta)$ and $V^{c, \eta}(r, \theta, a)$ and is the main result of the section, will be then a direct consequence of this analysis.
\medskip

\noindent{\em Description of the strategy}
\medskip

Let $\alpha\in\left]0,\frac{\pi}{8}\right]$ and $r_0>\eta$. We define
\begin{align*}
D^1_1&= \, ]\eta,r_0[\,\times\left( \left[\tfrac{3\pi}{2}+\alpha, \tfrac{7\pi}{4}\right[\, \cup \left[-\tfrac{\pi}{4},-\alpha\right] \right), \qquad D^{1}_2 = \{r_0\}\times \left(\left[\tfrac{3\pi}{2},\tfrac{7\pi}{4}\right[ \cup  \left[-\tfrac{\pi}{4},0\right] \,\right);  \\ 
D^{-1}_1&=\, ]\eta,r_0[\,\times\left[\tfrac{\pi}{2}+\alpha,\pi-\alpha\right],  \qquad\qquad\qquad\quad D^{-1}_2 = \{r_0\}\times \left[\tfrac{\pi}{2},\pi\right],
\end{align*}
and 
\begin{align*}
D^1&= D^1_1 \cup D^1_2,&\quad  D^{-1}= D^{-1}_1 \cup D^{-1}_2. 
\end{align*}
Let $C^1$  (resp. $C^{-1}$) be the complement sets of $D^1$ (resp. $D^{-1}$) in $[\eta,r_0]\times\left[-\tfrac{\pi}{4},\tfrac{7\pi}{4}\right[\,$ (see Fig.~\ref{rd09_15f1}). We denote by $\bar{C}^a$ the closure of $C^a$, $a \in \{+1, -1\}$. 
\begin{figure}
\begin{center}
\begin{tikzpicture}[scale=1.5, line cap=round]
  \def\thA{-18}    
  \def\thB{-72}    

  \def\circlethin{0.3pt}
  \def\thicklines{1.2pt}
  \def\shadeop{0.14}

  \draw[->, line width=\circlethin] (-1.2,0) -- (1.2,0) node[right] {$x$};
  \draw[->, line width=\circlethin] (0,-1.2) -- (0,1.2) node[above] {$y$};

  \draw[line width=\circlethin] (0,0) circle (1);

  \path (\thA:1) coordinate (A1) (\thA+180:1) coordinate (A2);
  \path (\thB:1) coordinate (B1) (\thB+180:1) coordinate (B2);

  \fill[black, opacity=\shadeop]
    (\thB:1) arc[start angle=\thB, end angle=\thA, radius=1] -- (0,0) -- cycle;
  \fill[black, opacity=\shadeop]
    (\thB+180:1) arc[start angle=\numexpr\thB+180\relax, end angle=\numexpr\thA+180\relax, radius=1]
    -- (0,0) -- cycle;

  \draw[line width=\thicklines] (A1) -- (A2);
  \draw[line width=\thicklines] (B1) -- (B2);

  \draw[line width=\thicklines] (-90:1) arc[start angle=-90, end angle=\thB, radius=1];
  \draw[line width=\thicklines] (\thA:1) arc[start angle=\thA, end angle=0, radius=1];
  \draw[line width=\thicklines] (90:1) arc[start angle=90, end angle=\numexpr\thB+180\relax, radius=1];
  \draw[line width=\thicklines] (162:1)
        arc[start angle=\numexpr\thA+180\relax, end angle=180, radius=1];

  \node at (-45:0.7) {$_{D^1}$};
  \node at (135:0.7) {$_{D^{-1}}$};

  \fill (0,0) circle (0.015);
\end{tikzpicture}\end{center}
\caption{The tacking regions $D^1$ and $D^{-1}$. \label{rd09_15f1}}
\end{figure}

The strategy $A$ can now be described as follows. When the boat is on tack $a \in \{+1, -1\}$,  $C^a$ is the \textit{continuation region}, whereas entering $D^a$ triggers a tack. With these definitions, since $D^a \subset C^{-a}$ (see Fig.~\ref{rd09_15f1}), the boat will then remain for a positive amount of time in $C^{-a}$  after each tack from $a$ to $-a$.

\begin{rem}\label{rd10_11r1}
Notice that with the strategy just described, the boat will never touch the no-cross line $d^1$ when on starboard tack, nor $d^{-1}$ when on port tack (see \eqref{d_p_m}). Further, when the boat enters $D^a$ after starting from a position $(r, \t)$ in $C^a$ with $r \leq r_0$, it either enters $D^a_1$, in which case it is at distance less than $r_0$ from the origin, or it enters $D^a_2$, in which case it has moved {\em back} to a position at distance $r_0$ from the origin. Of course, whenever the boat is on port tack and in Region $II^a$, or on starboard tack and in Region $IV^a$ (see Fig.~\ref{Fig:4}), its distance to the origin increases; it must never-the-less enter these regions if it seeks to reach the origin with a finite number of tacks. In addition, the boat will never exit the disk $\{r \leq r_0\}$ via the subset $\{r_0 \} \times ]0, \tfrac{3\pi}{2}[$ of the boundary when on starboard tack, nor via the subset $\{r_0 \} \times ([{-}\tfrac{\pi}{4}, \tfrac{\pi}{2}[\, \cup\, ]\pi, \tfrac{7\pi}{4}[$ when on port tack, since according to \eqref{state1}, \eqref{eq_drift_ra} and \eqref{eq_drift_pa}, $dR^A_t/dt$ is negative on these sets.
\end{rem}

\noindent{\em Formal construction of the strategy $A$ and the associated position process}
\medskip

Suppose that the initial position of the boat is $(r,\t, a)\in\, ]\eta,r_0[\, \times \, ]{-\piq},\pisq[\, \times \{-1, 1\}$. We define the strategy $A = (A_t)_{t \in \R_+ \cup \{0-\}}$ as follows: $A_{0-} = a$, and
\begin{equation*}
\left\{\begin{array}{l l}
A_{0} = -a, &\qquad \text{if} \   (r,\t)\in D^a, \\
A_{0} = a, &\qquad \text{if} \   (r,\t)\in C^a .
\end{array}
\right.
\end{equation*}
We recall that the coefficients $\mu_1(\t, a)$ and $\mu_2(r, \t, a)$ are defined in \eqref{eq_drift_ra} and \eqref{eq_drift_pa}. In order to define $A_t$ for $t > 0$, define $(R^{(1)}_t,\T^{(1)}_t)$ by $ (R^{(1)}_0,\T^{(1)}_0) = (r,\t)$, and for $t > 0$,
\begin{equation}\label{rd09_15e1}
\left\{
\begin{array}{l l}
dR^{(1)}_t &= 1_{ \{R^{(1)}_t> \eta\}} \,  \mu_1(\T^{(1)}_t,A_0) dt,\\
d\T^{(1)}_t &=  1_{\{ R^{(1)}_t> \eta\}}  \left[\mu_2(R^{(1)}_t, \T^{(1)}_t,A_0) dt + \s dB_t\right].
\end{array}
\right.
\end{equation}
Observe that a strong solution  $(R^{(1)}_t,\T^{(1)}_t)_{t \geq 0}$ to \eqref{rd09_15e1} exists by Remark \ref{strong_pol}, since the process never hits the lines $d^1$ (resp. $d^{-1}$) when $A_0 =1$ (resp.~$A_0 = -1$). Let $\psi^{(0)} = 0$, 
\begin{align*}
\tau^{(1)} &:= \inf \left\{ t\geq 0: \, R_t^{(1)}  = \eta   \right\} 
        , \\
\nu^{(1)}&:= \inf \left\{ t\geq 0: \, (R^{(1)}_t,\T_t^{(1)}) \in D^{A_0}_1   \right\} \wedge \tau^{(1)},\\
\rho^{(1)}&:= \inf \left\{ t\geq 0: \, (R^{(1)}_t,\T_t^{(1)})  \in D^{A_0}_2  \right\} \wedge \tau^{(1)},\\
\psi^{(1)}&:= \nu^{(1)} \wedge \rho^{(1)}
\end{align*}
(by convention, if the set over which the infimum is taken is empty, then the infimum is set to $+\infty$). Let $A_t = A_0$ for $t \in [0, \psi^{(1)}[$ (that is, the boat does not tack during $]0,\psi^{(1)}[$). Notice that on $\nu^{(1)} < \rho^{(1)}$, the boat exits the continuation region $C^{A_0}$ via $D^{A_0}_1$, and in this case, it has progressed towards the target, whereas if $\nu^{(1)} > \rho^{(1)}$, then the boat exits $C^{A_0}$ via $D^{A_0}_2$, and in this case, the boat's radial coordinate is again $r_0$ at this time.

For $i \geq 2$, assuming that $(R^{(i-1)}_t,\T^{(i-1)}_t)$ and  $\tau_\eta^{(i-1)}$,  $\nu^{(i-1)}$, $\rho^{(i-1)}$ and $\psi^{(i-1)}$ have been defined, as well as $A_t$ for $t \in [0, \psi^{(i-1)}[$, define $(R^{(i)}, \T^{(i)}_t, \, t \geq \psi^{(i-1)})$ as follows:
 $$ 
     (R^{(i)}_{\psi^{(i-1)}},\T^{(i)}_{\psi^{(i-1)}}) = (R^{(i-1)}_{\psi^{(i-1)}},\T^{(i-1)}_{\psi^{(i-1)}}), 
$$ 
and for $t > \psi^{(i-1)}$, 
\begin{equation}\label{A_tau_def}
\left\{
\begin{array}{ll}
dR^{(i)}_t &=  \mu_1(\T^{(i)}_t, - A_{\psi^{(i-1)}}) dt,\\
d\T^{(i)}_t &=   \mu_2(R^{(i)}_t, \T^{(i)}_t, - A_{\psi^{(i-1)}}) dt + \s dB_t.
\end{array}
\right.
\end{equation}
Let
\begin{align}\notag
\tau^{(i)} &:= \inf \left\{ t\geq \psi^{(i-1)}: \, R_t^{(i)}  = \eta   \right\} 
       ,\\ \notag
\nu^{(i)}&:= \inf \left\{ t\geq  \psi^{(i-1)}: \, (R^{(i)}_t,\T_t^{(i)}) \in D^{- A_{\psi^{(i-1)}}} _1   \right\} \wedge \tau^{(i)},\\ \notag
\rho^{(i)}&:= \inf \left\{ t\geq \psi^{(i-1)}: \, (R^{(i)}_t,\T_t^{(i)}) \in D^{- A_{\psi^{(i-1)}}} _2  \right\} \wedge \tau^{(i)},\\
\psi^{(i)}&:= \nu^{(i)} \wedge \rho^{(i)}.
\label{rd09_16e2}
\end{align}
Let $A_t = -A_{\psi^{(i-1)}}$ for $t \in [\psi^{(i-1)}, \psi^{(i)}[$ (that is, the boat tacks at time $\psi^{(i-1)}$ but does not tack during $]\psi^{(i-1)}, \psi^{(i)}[$). This extends the definition of $A_t$ to $[0, \psi^{(i)}[$. Observe that a strong solution to \eqref{A_tau_def} up to time $\psi^{(i)}$ exists by Remark \ref{strong_pol}.

For $i \geq 1$, the process $(A_t)_{t \geq 0}$ satisfies $A_t = (-1)^{i-1} A_0$ on $[\psi^{(i-1)},\psi^{(i)} \wedge \tau^{(i)}[$, that is, the boat tacks at each time $\psi^{(i)}$ (when $\psi^{(i)} < \tau^{(i)}$, otherwise the boat has reached the target). 
Then, during $]\psi^{(i)}, \psi^{(i+1)}\wedge \tau^{(i+1)}_\eta[$, the boat navigates through the continuation region $C^{A_{\psi^{(i)}}}$ until reaching the region $D^{A_{\psi^{(i)}}}$ (or the target). Then the boat tacks again and continues through the region $C^{-A_{\psi^{(i)}}} = C^{A_{\psi^{(i+1)}}}$.

Notice that on $\nu^{(i)} < \rho^{(i)}$, the boat exits the continuation region $C^{-A_{\psi^{(i-1)}}}$ via $D^{-A_{\psi^{(i-1)}}}_1$, and in this case, it has progressed towards the target, whereas if $\nu^{(i)} > \rho^{(i)}$, then the boat exits $C^{-A_{\psi^{(i-1)}}}$ via $D^{-A_{\psi^{(i-1)}}}_2$, and in this case, the boat's radial coordinate is again $r_0$ at this time. In order to check that the target is reached with probability one, it will be necessary to show that this second situation can only occur a finite number of times.

Define 
\begin{equation*}
\tilde{i} = \inf\{i >0 : \psi^{(i)} = \tau^{(i)} \}. 
\end{equation*}
On the event $\{\tilde{i} < \infty\}$, the yacht reaches $\B_\eta$ in a finite amount of time. 
Observe that on this event,  $ \psi^{(i)} = \psi^{(\tilde i)}$ for $i \geq \tilde i$, so after tacking for the last time at $\psi^{(\tilde i -1)}$, the yacht reaches the target.

The  process  $(R^{A}_t,\T^{A}_t)_{t \geq 0}$ is defined by
\begin{equation}\label{RA_tau_def}
(R^{A}_t,\T^{A}_t) = \left\{
\begin{array}{ll}
(R^{(i)}_t, \T^{(i)}_t ), &\text{for  } \, i \leq \tilde i, \  t \in [\psi^{(i-1)} ,\psi^{(i)} [, \\
(\eta, \T^{(\tilde i)}_{\psi^{(\tilde{i})}} ), &\text{for  } \, t \geq  \psi^{(\tilde i)}.
\end{array}
\right.
\end{equation}

\noindent{\em Towards admissibility of $A$}
\medskip

In order to show that $A$ satisfies Items (a) and (b) of Definition \ref{Def:3.3}, we need to prove that it is piecewise continuous and that \eqref{state1} admits a strong solution for $t\in [0,\infty[$.
\begin{lem}\label{strong_sol}
The process $A =(A_t)$ defined above is piecewise constant and the system \eqref{state1} admits a strong solution $(R^{A}_t,\T^{A}_t)_{t \geq 0}$, defined in \eqref{RA_tau_def}.
\end{lem}

\begin{proof}
In order to show that $A$ is piecewise constant, fix $\omega \in \Omega$ and $T<\infty$. Suppose that $A_0 = a \in \{+1, -1\}$. In order for the strategy $A_t$ to switch from $a$ to $-a$ at time $t > 0$, the boat must enter the closed set $D^a$ while starting from a position in the open set $C^a$. This requires a strictly positive amount of time and is sufficient to ensure that $t \mapsto A_t(\omega)$ is piecewise constant.

Since $d^1$ (respectively $d^{-1}$)  is in the interior of $D^1$ (respectively $D^{-1}$), the process $((R^{(i)}_t, \T^{(i)}_t ))_{t \geq \psi^{(i-1)}}$ will not hit $d^{-A_{\psi^{i-1}}}$ before time $\psi^{(i)}$ (see Fig.~\ref{rd09_15f1}). The process $(R^{A}_t,\T^{A}_t)$  obtained by ``gluing together" the $(R^{(i)}_t, \T^{(i)}_t )$ as described in \eqref{RA_tau_def} is then the strong solution to \eqref{state1}.
\end{proof}

\begin{rem}\label{tilde_tau_A_def}
Now that the process $(R^{A}_t,\T^{A}_t)_{t \geq 0}$ is defined, it is straightforward to verify that $\tau^A_\eta = \tau^{(\tilde i - 1)}_{\eta}$ a.s., where $\tau^A_\eta$ is defined in \eqref{tau_n}. From now on, we will denote $\tau^A:= \tau^A_\eta$.
\end{rem}
\medskip

\noindent{\em Reducing to a single tack by using symmetry}
\medskip

Rather than having two different pictures for the starboard and port tacks, as in Figures \ref{Fig:4} and \ref{Fig_brw_sailing_rules_polar}, it is convenient to use the symmetries in the problem to reduce to one single picture, as we now describe.

\begin{figure}
\begin{center}
\includegraphics[width=6cm]{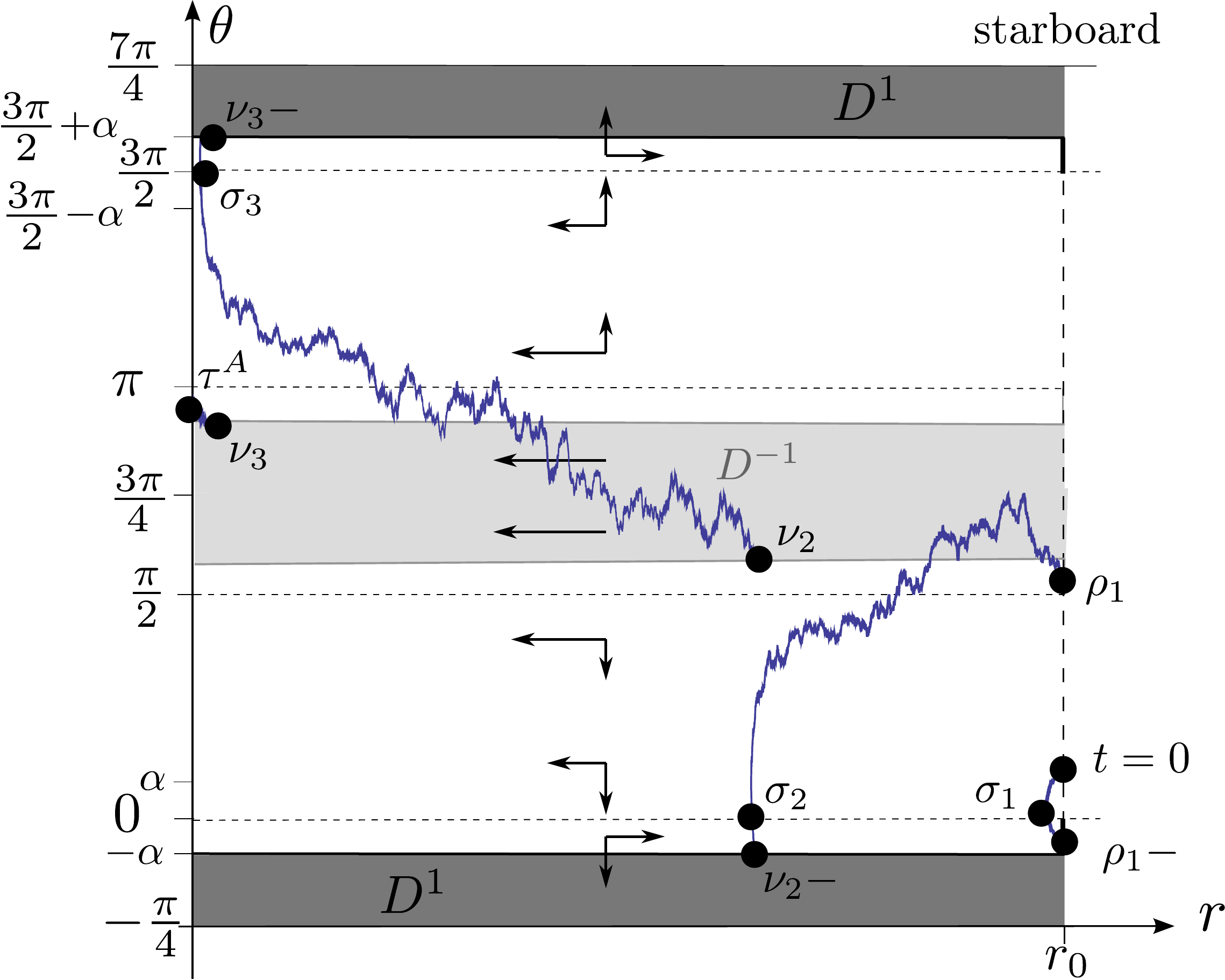}
\end{center}
\caption{Tacking and continuation regions of the strategy $A$. Length of arrows indicate the drift's  intensity in the directions of the $r$ and $\theta$ axes.  Each time the process hits $D^1$ (dark grey), $A$ prescribes to tack. We picture this by a jump to the symmetric point on $D^{-1}$ (light grey). Notice that $D^1$ contains the vertical segments $\{r_0\}\times\left[-\tfrac{\pi}{4},0\right]$ and $\{r_0\}\times\left[\tfrac{3\pi}{2},\tfrac{7\pi}{4}\right[\,$.
The bullet points show the position of the yacht at the stopping times of \eqref{eq_def_brw_A_stop_times}.}
\label{Fig_brw_exist_A}
\end{figure}

Under the symmetry $S: \R^2 \to \R^2$ defined by $S(x, y) = (y, x)$, $D^1$ transforms into $D^{-1}$ and $C^1$ transforms into $C^{-1}$ (see Fig.~\ref{rd09_15f1}). Therefore, for all $(x,y)\in\R^2$, for the strategy $A$ defined above, the payoff function \eqref{Jetac} satisfies, 
\begin{equation}\label{eq_brw_cartesian_symmetry}
J^{c,\eta}((x,y),1,A)=J^{c,\eta}((y,x),-1,A).
\end{equation}
Equivalently, in polar coordinates,  for all $(r,\theta)\in\R_+\times\left[-\tfrac{\pi}{4},\tfrac{7\pi}{4}\right[$, we have
\begin{equation}\label{eq_brw_polar_symmetry}
J^{c,\eta}((r,\theta),1,A)=J^{c,\eta}\left(\left(r,\tfrac{\pi}{2}-\theta\right),-1,A\right),
\end{equation}
where $\tfrac{\pi}{2}-\theta$ is considered modulo $2\pi$ in $\left[-\tfrac{\pi}{4},\tfrac{7\pi}{4}\right[$. Therefore, instead of considering both tacks, as in Fig.~\ref{Fig_brw_sailing_rules_polar}, we can consider only 
the starboard tack, provided each tack is replaced by a jump of the position process: 
\begin{equation}\label{rd09_16e1}
   \text{if } A_t \neq A_{t-}, \text{ then } \T^A_t = \tfrac{\pi}{2} - \T^A_{t-} \text{ modulo } 2 \pi \text{ in } \left[-\tfrac{\pi}{4},\tfrac{7\pi}{4}\right[.
\end{equation}
In particular, if $\Theta^A_{t-} = -\alpha$ (resp. $\Theta^A_{t-}= \pitm + \alpha $), then $\Theta^A_t=\frac{\pi}{2}+\alpha$ (resp $\Theta^A_t= \pi-\alpha$). From now on, we will consider that the yacht is always on starboard tack when it is in $C^1$ and that its motion is governed by \eqref{state1}, and when the boat tacks, the angle process $(\T^a_t)$ jumps according to the rule \eqref{rd09_16e1}, and then moves again according to \eqref{state1}. With this rule, $\P_{r,\theta,1}$ becomes $\P_{r,\theta}$ and the motion $(R^A_t, \T^A_t)$ is piecewise continuous. Further, $\T^A_t$ always belongs to $[-\frac{\pi}{4}, \frac{7 \pi}{4}[$. We use this description of the motion process in Fig.~\ref{Fig_brw_exist_A}.

Notice that there is also  a symmetry property about the line with equation $y = -x$:
\begin{equation}\label{eq_brw_polar_symmetry_2}
J^{c,\eta}((r,\theta),1,A)=J^{c,\eta}\left(\left(r,\tfrac{3\pi}{2}-\theta\right),1,A\right),\qquad\text{for all } (r,\theta)\in\R_+\times\left[-\tfrac{\pi}{4},\tfrac{7\pi}{4}\right[,
\end{equation}
where $\theta \mapsto \tfrac{3\pi}{2}-\theta$ is one-to-one from $\left[-\tfrac{\pi}{4},\tfrac{7\pi}{4}\right[$ into itself.

Using this new description, we now reformulate the definitions of the stopping times in \eqref{rd09_16e2}. Recall  that $\tau^A$ is defined in Remark \ref{tilde_tau_A_def}.

\begin{dfn}\label{def_stopping_times_exist}
Let $(R^A_t,\Theta^A_t)$ denote the position at time $t$ of the yacht when using the strategy $A$. We set $\sigma_0=\nu_0=\rho_0=\psi_0=0$, $E_0=\emptyset$ and for all $i\in\N^*$, we define
\begin{equation}\label{eq_def_brw_A_stop_times}
\begin{array}{ll}
\sigma_i&=\inf\left\{t\geq\psi_{i-1}:\Theta^A_t \in [-\tfrac{\pi}{4},0] \cup [\tfrac{3\pi}{2}, \tfrac{7\pi}{4}]\right\}\wedge\tau^A, \\[1ex]
\nu_i&=\inf\left\{t\geq \sigma_i:(R^A_t,\Theta^A_t) \in D^1_1\right\}\wedge\tau^A, \\[1ex]
\rho_i&=\inf\left\{t\geq \sigma_i :(R^A_t,\T^A_t) \in D^1_2\right\}\wedge\tau^A,\\
\psi_i&=\nu_i\wedge\rho_i ,\\
u_i&=\psi_i-\psi_{i-1} 
\end{array}
\end{equation} 
(see Fig.~\ref{Fig_brw_exist_A}) and
$$
   E_i=\left\{\sup_{\psi_{i-1}\leq t\leq\tau^A} \left|\sigma B_t-\sigma B_{\psi_{i-1}}\right|<\tfrac{\alpha}{2}\right\}.
$$
In particular, $E_i$ is the event ``the wind is calm (its angle varies by less than $\tfrac{\alpha}{2}$) from time $\psi_{i-1}$ until the target is reached."
\end{dfn}


The stopping times of (\ref{eq_def_brw_A_stop_times}) are illustrated in Fig.~\ref{Fig_brw_exist_A}. The times $\psi_i$ are the tacking times of the strategy $A$ (replaced here by jumps of the process $(\T^A_t)$). 
The times $u_i$ are {\em holding times}, that is, the {\em times between tacks}. Clearly, the stopping times defined in \eqref{eq_def_brw_A_stop_times} are in fact the same as those defined in \eqref{rd09_16e2}:  $\nu_i = \nu^{(i)}$, $\rho_i = \rho^{(i)}$ and $\psi_i = \psi^{(i)}$, $\P_{r, \t, a}$-a.s.
We also note that according to their definitions, $\psi_i\leq\tau^A$ and if the event $\{\psi_i=\tau^A\}$ occurs, then $E_{i+1}$ occurs trivially. 
\medskip

\noindent{\em Bounding the expected time between tacks}
\medskip


\begin{lem}\label{Lem_exist_A_ui}
For all $(r,\theta)\in\bar{C}^1$, the holding times $u_i$, $i \in \N^*$, of the strategy $A$ defined in \eqref{eq_def_brw_A_stop_times} satisfy
$$
   \E_{r,\theta}\left(u_i\right)\leq \sigma^{-2} \left(\alpha + \tfrac{3 \pi}{4} \right)^2 .
$$
Moreover, the stopping times $\psi_i$ defined in \eqref{eq_def_brw_A_stop_times} have finite expectations.
\end{lem}

\begin{proof}
Let $(r,\theta)\in\bar{C}^1$ be given numbers and let $f(x)=\left(x-\tfrac{3\pi}{4}\right)^2$. By It\^{o}'s formula, since there is no tack during the interval $]0, \psi_1[$ and therefore $t \mapsto \Theta^A_{(t\wedge\psi_1)_-}$ is a continuous diffusion,
\begin{align*}
f\left(\Theta^A_{(t\wedge\psi_1)_-}\right)-f\left(\theta\right)&=\int_0^{t\wedge\psi_1}\!\!2 \left(\Theta^A_s-\tfrac{3\pi}{4}\right)  \mu_2\left(R^A_s,\Theta^A_s,1\right) ds\\
&\qquad +2\sigma\int_0^{t\wedge\psi_1}\!\left(\Theta^A_s-\tfrac{3\pi}{4}\right)\,dB_s+\sigma^2 \left(t\wedge\psi_1\right).
\end{align*}
On one hand, $\Theta^A_s \in [-\pi/4, 7\pi/4]$ is bounded. On the other hand, for all $(r,\theta)\in \bar{C}^1$, $f\left(\theta\right)\leq\left(-\alpha-\tfrac{3\pi}{4}\right)^2$ and $ \left(\theta-\tfrac{3\pi}{4}\right) \mu_2\left(r,\theta,1\right)\geq0$ (just check the three cases in \eqref{eq_drift_ra}). Hence, taking the expectation in the previous equation and using these inequalities, we find that $\sigma^2\, \E_{r,\theta}\left(t\wedge\psi_1\right) \leq \E\left(f(\Theta^A_{(t\wedge\psi_1)_-}) \right) - f(\theta) \leq  \left(\alpha+\tfrac{3\pi}{4}\right)^2$. Therefore, letting $t \uparrow \infty$,  the monotone convergence theorem implies that
\begin{align}\label{rd09_16e3}
    \E_{r,\theta}\left(\psi_1\right)\leq \sigma^{-2} \left(\alpha+\tfrac{3\pi}{4}\right)^2=:K, \qquad\text{for all } (r, \t) \in \bar{C}^1.
 \end{align}
Observe that by \eqref{eq_def_brw_A_stop_times}, $u_1=\psi_1$. 
Furthermore, by the strong Markov property and \eqref{rd09_16e3}, for all $i\in\N^*$, we have 
$$
   \E_{r,\theta}(u_i)=\E_{r,\theta}\left(\E_{r,\theta}\left(u_i \mid \shf_{\psi_{i-1}}\right)\right)=\E_{r,\theta}\left(\E_{R^A_{\psi_{i-1}},\Theta^A_{\psi_{i-1}}}(u_1)\right)\leq K,
 $$
because $(R^A_{\psi_{i-1}},\Theta^A_{\psi_{i-1}})$ belongs to $\bar C^{1}$ by the definition of the strategy $A$. 

Finally, $\psi_i = u_1 + \cdots + u_i$. Thus, $\E_{r, \t}(\psi_i) < \infty$, for all $i \in \N$.
\end{proof}
\medskip

\noindent{\em Reducing the relevant starting positions}

\begin{rem}\label{rem_favorable_initial_conditions}
A consequence of this lemma is that, without loss of generality, we can consider only starting positions $(r,\theta)$ such that $r\leq r_0$ and $\theta\in\,\left]\tfrac{\pi}{2} - \alpha, \tfrac{\pi}{2}+\alpha\right]$. Indeed, with at most one tack and in a time with finite expectation, using \eqref{rd09_16e1}, the yacht will be at
$$(R^A_{\psi_1}, \Theta^A_{\psi_1})\in\left([\eta,r_0[\times\left\{\tfrac{\pi}{2}+\alpha,\pi-\alpha\right\}\right)\,\cup\,\left( \{r_0\}\times\left(\,\left[\tfrac{\pi}{2},\tfrac{\pi}{2}+\alpha\right]\cup[\pi-\alpha,\pi]\,\right)\right).$$
Moreover, by the second symmetry \eqref{eq_brw_polar_symmetry_2}, we can assume that
$$(R^A_{\psi_1}, \Theta^A_{\psi_1})\in\left([\eta,r_0[\times\left\{\tfrac{\pi}{2}+\alpha\right\}\right)\,\cup\,\left( \{r_0\}\times\left[\tfrac{\pi}{2},\tfrac{\pi}{2}+\alpha\right]\right)$$
and this set is clearly included in $[\eta,r_0]\times\left[\tfrac{\pi}{2} - \alpha, \tfrac{\pi}{2}+\alpha\right]$.
From now on, we will denote
\begin{equation}\label{s_def}
\shs = \shs_{\eta, r_0, \alpha}:= \, ]\eta,r_0]\times\left[\tfrac{\pi}{2} - \alpha, \tfrac{\pi}{2}+\alpha\right].
\end{equation}
In particular, when the boat starts in the region $\shs$, we will always have $M_0(A) = 0$.
\end{rem}


\noindent{\em Bounding the time to exit $\shs$}
\medskip

The next lemma states that while the boat is in the region $\shs$, it will not move away from the origin, and it also provides a deterministic bound on the time needed to exit $\shs$.

\begin{lem}\label{Lem_exist_A_zeta}
Consider the stopping time
\begin{equation}\label{eq_def_st_zeta}
\zeta=\inf\left\{t\geq0:\Theta^A_t\notin\left[\tfrac{\pi}{2}-\alpha,\,\tfrac{\pi}{2}+\alpha\right]\right\}\wedge\tau^A .
\end{equation}
For all initial positions $(r,\theta)\in \shs$,
\begin{align}\label{rd09_20e1}
   R^A_{t} < r, \quad t \in\, ]0,\zeta], 
   \quad\text{ and }\qquad \zeta\leq\frac{r-\eta}{v\cos\alpha}, \qquad\P_{r,\theta}\text{-a.s.}
\end{align}
\end{lem}

\begin{proof}
By \eqref{eq_drift_ra}, for all $0\leq t\leq\zeta$,  $\mu_1\left(\Theta^A_t,1\right)\in[-v,-v\cos\alpha]$. Hence, for all $t \geq 0$,
$$
  0 \leq  \eta\leq R^A_{t\wedge\zeta}=r+\int_0^{t\wedge\zeta}\mu_1\left(\Theta^A_s,1\right)\,ds\leq r- ( \cos\alpha) \,(t\wedge\zeta) v.
 $$
Therefore, $R^A_{t}\leq r$ for $t\in\, ]0,\zeta]$, and $\zeta\leq\frac{r-\eta}{v\cos\alpha}$.
\end{proof}
\medskip

\noindent{\em Reaching the target when the wind is calm}
\medskip

The next proposition shows that if the wind is calm from time $0$ on, that is, if $E_1$ occurs, then the time to reach $D^1$ is bounded and the yacht is never pushed back into the circle of radius $r_0$ during $[0, \psi_1]$. In order to establish this, we will make use of the equations of motion both in polar coordinates and in Cartesian coordinates.

\begin{prop}\label{Lem_exist_A_cst}
If the event $E_1$ occurs, then for all initial positions $(r,\theta)\in \shs$, the following inequalities hold $\P_{r,\theta}$-almost surely:
\begin{align}
\sigma_1-\zeta &\leq \frac{d_1\, r}{v}\,,\label{eq_Lem_cst_1}\\
  \psi_1-\sigma_1  &\leq \frac{d_2\, R^A_{\sigma_1}}{v}\,, \label{eq_Lem_cst_3}
\end{align}
where
\begin{align}
d_1&=\sin\left(\tfrac{\alpha}{2}\right)\tan\alpha+1, \qquad 
d_2=\frac{\left(\sin\left(\tfrac{\alpha}{2}\right)+\tan\alpha\right)\cos\alpha}{\cos(2\alpha)}.
\label{rd10_20e1}
\end{align}
Moreover, on $E_1$,
\begin{equation}\label{r_leq_r0}
R_t^A \leq r, \quad t\in [0,\psi_1], \qquad\text{ and } \qquad R^A_{\psi_1} < r.
\end{equation}
\end{prop}

\begin{proof}
In order to establish this statement, we will treat three cases separately. Let $\zeta$ be defined as in \eqref{eq_def_st_zeta}.
\medskip

\textbf{Case 1:} $R^A_{\zeta}=\eta$ (which is equivalent to $\zeta=\tau^A$). If $\zeta=\tau^A$, then by definition of the strategy $A$, $(R_t^A, \Theta^A_t) \in \shs$ for $t \in [0, \tau^A[$, therefore $M_{\tau^A}(A)=0$,  $\sigma_1=\nu_1=\rho_1=\psi_1=\tau^A$ and since $R^A_{\tau^A}=\eta$, the inequalities \eqref{eq_Lem_cst_1} and \eqref{eq_Lem_cst_3} are satisfied because their left-hand sides vanish, and \eqref{r_leq_r0} follows from Lemma \ref{Lem_exist_A_zeta}.
\medskip

\textbf{Case 2:} $\Theta^A_{\zeta}=\tfrac{\pi}{2}+\alpha$ and $R^A_{\zeta}> \eta$. In this case, on the event $E_1$,  $\Theta^A_t\in\left[\tfrac{\pi}{2},\pi\right]$ for all $t\geq\zeta$. Indeed, if we set
$$
   \kappa=\left\{\begin{array}{ll}
\inf\left\{t\geq\zeta:\Theta^A_t\notin\left[\tfrac{\pi}{2}\,,\pi\right]\right\},& \text{ if } \{\cdots\}\neq\emptyset,\\
+\infty,&\text{ otherwise},\\
\end{array}\right.$$
then by  \eqref{eds_proc_controle_polar}, on $\{\kappa<+\infty\}$,
\begin{align}\label{rd09_17e1}
   \Theta^A_\kappa=\Theta^A_\zeta+\int_\zeta^\kappa \mu_2\left(R^A_s,\Theta^A_s,1\right) \,ds+ \sigma B_\kappa - \sigma B_\zeta .
\end{align}
Observe that $\T^A_s \in [\tfrac{\pi}{2}, \pi]$ for $s \in [\zeta, \kappa[$, so $\mu_2\left(R^A_s,\Theta^A_s,1\right) = 0$ by \eqref{eq_drift_ra}, and since $E_1$ occurs, $\vert \sigma B_\kappa\vert < \tfrac{\alpha}{2}$ and $\vert \sigma B_\zeta \vert < \tfrac{\alpha}{2}$. It follows that on $E_1$,
$$
      \Theta^A_\kappa =  \tfrac{\pi}{2} + \alpha +  \sigma B_\kappa -  \sigma B_\zeta > \tfrac{\pi}{2}+\alpha- \tfrac{\alpha}{2}- \tfrac{\alpha}{2}=\tfrac{\pi}{2}.
$$
Similarly, on $E_1$,
$$
   \Theta^A_\kappa<\tfrac{\pi}{2}+\alpha+2\tfrac{\alpha}{2}=\tfrac{\pi}{2}+2\alpha<\pi.
$$
Since $\Theta^A_\kappa \in \{\tfrac{\pi}{2},\pi\}$ on $\{\kappa<+\infty\}$, this contradiction shows that $\P\{\kappa<+\infty\}=0$ and therefore, $\kappa\equiv+\infty$. It follows that $\Theta^A_t\in\left[\tfrac{\pi}{2},\pi\right]$ for all $t\geq\zeta$.
By \eqref{eds_proc_controle_polar},  for all $t\geq\zeta$,
\begin{equation}\label{rd09_17e2}
R^A_{t\wedge\tau^A} = R^A_\zeta+\int_\zeta^{t\wedge\tau^A} \mu_1\left(\Theta^A_s,1\right)\,ds.
\end{equation}
By \eqref{eq_drift_ra}, $\mu_1(\Theta^A_s,1) = -v$, therefore
\begin{equation}\label{eq45a}
    R^A_{t\wedge\tau^A} = R^A_\zeta-(t\wedge\tau^A-\zeta) v.
\end{equation}
Taking the limit when $t\rightarrow+\infty$, we find that $\tau^A-\zeta = \frac{R^A_\zeta - R^A_{\tau^A}}{v} = \frac{R^A_\zeta - \eta}{v}\leq\frac{r}{v}$ by Lemma \ref{Lem_exist_A_zeta}. Since $ \zeta < \infty$, it follows that $\tau^A < \infty$, and since $\kappa=+\infty$, the boat never exits $[\tfrac{\pi}{2}, \pi]$, hence it never enters $D^1$. It follows that $\psi_1 = \sigma_1=\nu_1=\rho_1=\tau^A$  and $R^A_{\sigma_1}=\eta  \leq r$. Thus, the inequality \eqref{eq_Lem_cst_1} is satisfied because $d_1 > 1$ and \eqref{eq_Lem_cst_3}  is also satisfied because $\psi_1 - \sigma_1 = 0$. 
By \eqref{eq45a}, it follows that $R_t^A \leq R_\zeta^A < r$ for $t\in [\zeta,\psi_1]$, where the last inequality follows from Lemma \ref{Lem_exist_A_zeta}, hence we have also proved  \eqref{r_leq_r0}.
\medskip

\textbf{Case 3:} $\Theta^A_{\zeta}=\tfrac{\pi}{2}-\alpha$ and $R^A_{\zeta} > \eta$. On $E_1$, we see from \eqref{eq_drift_ra} that $\mu_2(R^A_s, \T^A_s, 1) = -\tfrac{v \cos \T^A_s}{r} \leq 0$, therefore, by \eqref{rd09_17e1} and the definition of $\nu_1$ in \eqref{eq_def_brw_A_stop_times}, for $t \in [\zeta, \psi_1[$, 
\begin{align*}
  - \alpha \leq \Theta^A_t \leq \T_\zeta + \sigma\, (B_t - B_\zeta) < \tfrac{\pi}{2} - \alpha + \tfrac{\alpha}{2} + \tfrac{\alpha}{2} = \tfrac{\pi}{2},
\end{align*}
that is,
\begin{equation}\label{eq_thetaA_zeta_tau1}
\Theta^A_t\in\left[-\alpha, \tfrac{\pi}{2}\right[\,, \qquad\qquad\text{ for all } t\in[\zeta,\psi_1[\,.
\end{equation}
Therefore, by Proposition \ref{prop_brw_SDE_solut} and by \eqref{mu_def_vec}, since the yacht is in zone $I \cup IV^a$ (see Fig.~\ref{Fig:4}), the process $\vec{X}^A_t$, written in Cartesian coordinates,  satisfies for all $t\in[\zeta,\psi_1[$,
\begin{align*}
X^A_t&=\cos(\sigma B_t)X^A_\zeta-\sin(\sigma B_t)Y^A_\zeta+v\int_\zeta^t\sin(\sigma B_t-\sigma B_s)ds,\\
Y^A_t&=\sin(\sigma B_t)X^A_\zeta+\cos(\sigma B_t)Y^A_\zeta-v\int_\zeta^t\cos(\sigma B_t-\sigma B_s)ds.
\end{align*}
Since $\T^A_\zeta =  \tfrac{\pi}{2} - \alpha$, we have $Y^A_\zeta > 0$ and $X^A_\zeta > 0$. Also, on $E_1$, $|\sigma B_t|<\tfrac{\alpha}{2}$ for all $t\leq\tau^A$, and we can bound $Y^A_t$ in the second equation from above and below as follows:
\begin{align*}
&-\sin\left(\tfrac{\alpha}{2}\right)X^A_\zeta+\cos\left(\tfrac{\alpha}{2}\right)Y^A_\zeta- (t-\zeta) v  < Y^A_t  < \sin\left(\tfrac{\alpha}{2}\right)X^A_\zeta+ Y^A_\zeta- (t-\zeta) v \cos\alpha,
\end{align*}
for all $t\in[\zeta,\psi_1[\,$. Noting that $X^A_\zeta=R^A_\zeta\cos\left(\tfrac{\pi}{2}-\alpha\right)=R^A_\zeta\sin\alpha$ and  $Y^A_\zeta=R^A_\zeta\cos\alpha$, the above inequalities become
\begin{equation}\label{eq_proof_lem_cst_Ya}
\cos\left(\tfrac{3\alpha}{2}\right)\,R^A_\zeta-(t-\zeta) v < Y^A_t < \left(\sin\left(\tfrac{\alpha}{2}\right)\sin\alpha+\cos\alpha\right)R^A_\zeta- (t-\zeta) v \cos\alpha.
\end{equation}
Observe that on $E_1$, by \eqref{eq_thetaA_zeta_tau1} and under our assumptions,  $Y^A_{\sigma_1} \geq 0$ and a strict inequality corresponds to the yacht hitting the target before the line $y=0$, that is $\sigma_1 = \tau^A$. Replacing $t$ by $\sigma_1$, 
recalling the definition of $d_1$ and rearranging terms, these inequalities become
\begin{equation}\label{eq_proof_lem_cst_Ya_sigma1}
    \cos\left(\frac{3 \alpha}{2}\right) R^A_\zeta < (\sigma_1-\zeta) v <d_1\,R^A_\zeta -\frac{Y^A_{\sigma_1}}{\cos \alpha}.
\end{equation}
By \eqref{rd09_20e1}, we conclude that $\sigma_1-\zeta<\tfrac{d_1 r}{v}$ and this establishes inequality \eqref{eq_Lem_cst_1}. 

Turning to \eqref{eq_Lem_cst_3}, notice that if $\sigma_1=\tau^A$, then $\psi_1 = \nu_1=\rho_1=\tau^A$, $R^A_{\sigma_1}=\eta$ and \eqref{eq_Lem_cst_3} is trivial, so from now on, we assume that $\sigma_1<\tau^A$.

On $E_1$, since $X^A_{\sigma_1}=R^A_{\sigma_1}$ and $Y^A_{\sigma_1}=0$, we have for all $t\in[\sigma_1,\psi_1[$, by Proposition \ref{prop_brw_SDE_solut} and by \eqref{mu_def_vec}, since the yacht is in zone $I \cup IV^a$, 
\begin{align*}
X^A_t&=\cos(\sigma B_t)R^A_{\sigma_1}+v\int_{\sigma_1}^t\sin(\sigma B_t-\sigma B_s)ds,\\
Y^A_t&=\sin(\sigma B_t)R^A_{\sigma_1}-v\int_{\sigma_1}^t\cos(\sigma B_t-\sigma B_s)ds,
\end{align*}
and so $X^A_t\leq R^A_{\sigma_1}+(t-\sigma_1) v\sin\alpha$ and
\begin{equation}\label{eq_proof_lem_YA_ineg}
-\sin\left(\tfrac{\alpha}{2}\right) R^A_{\sigma_1}-(t-\sigma_1) v \leq\; Y^A_t\;\leq \sin\left(\tfrac{\alpha}{2}\right) R^A_{\sigma_1}-(t-\sigma_1) v \cos\alpha.
\end{equation}
For $t \in [\sigma_1, \psi_1[$, the boat has not yet entered $D^1$, so $X^A_t\tan\alpha+Y^A_t>0$, and we find by combining these inequalities that
$$0< X^A_t\tan\alpha+Y^A_t\leq R^A_{\sigma_1}\left(\tan\alpha+\sin\left(\tfrac{\alpha}{2}\right)\right)+(t-\sigma_1) v \left( \tan\alpha\sin\alpha-\cos\alpha\right).$$
We note that $\cos\alpha - \tan\alpha\sin\alpha > 0$ if $\cos^2 \alpha > \sin^2 \alpha$, and this is the case when $\alpha \in [0, \pi/4[$. Therefore, we deduce that
\begin{equation}\label{eq_proof_lem_t-sigma}
t-\sigma_1\leq\frac{R^A_{\sigma_1}}{v}\,\frac{\tan\alpha+\sin\left(\tfrac{\alpha}{2}\right)}{\cos\alpha-\tan\alpha\sin\alpha}=\frac{R^A_{\sigma_1}}{v}\,d_2
\end{equation}
by definition of $d_2$. Taking the limit $t\rightarrow\psi_1$, we obtain inequality \eqref{eq_Lem_cst_3}.
 
   We now prove the inequalities in \eqref{r_leq_r0}.
 By \eqref{eq_thetaA_zeta_tau1}, \eqref{eq_drift_ra} and \eqref{eds_proc_controle_polar}, on $E_1$, for all $t\in[\zeta,\psi_1[$, $R^A_t$ satisfies  
\begin{equation}\label{eq45b}
dR^A_t=-v\sin\left(\Theta^A_t\right)\,dt=-v\frac{Y^A_t}{R^A_t}\,dt.
\end{equation}
 Hence, we have
\begin{equation}\label{eq_EDS_RtA_2}
d\left(R^A_t\right)^2=2R^A_t\,dR^A_t =-2v Y^A_t\,dt,
\end{equation}
and by \eqref{eq_proof_lem_cst_Ya},
\begin{align}\nonumber
\left(R^A_t\right)^2&\leq \left(R^A_\zeta\right)^2-2v\int_\zeta^t \left(\cos\left(\tfrac{3\alpha}{2}\right)\,R^A_\zeta-(s-\zeta) v \right)\,ds\\
&=\left(R^A_\zeta\right)^2+\left((t-\zeta) v -\cos\left(\tfrac{3\alpha}{2}\right)\,R^A_\zeta\right)^2-\cos^2\left(\tfrac{3\alpha}{2}\right)\,\left(R^A_\zeta\right)^2.
\label{rd10_13e1}
\end{align}
Setting $t= \sigma_1$ and using \eqref{eq_proof_lem_cst_Ya_sigma1} to observe that the quantity inside the square is nonnegative,  \eqref{eq_Lem_cst_1} and Lemma \ref{Lem_exist_A_zeta} imply that
\begin{equation}\label{eq_415a}
\left(R^A_{\sigma_1}\right)^2\leq\left(R^A_\zeta\right)^2+\left(d_1R^A_\zeta-\cos\left(\tfrac{3\alpha}{2}\right)\,R^A_\zeta\right)^2-\cos^2\left(\tfrac{3\alpha}{2}\right)\,\left(R^A_\zeta\right)^2\leq r^2c_1^2,
\end{equation}
where $c_1 =(1+\left(d_1-\cos\left(\tfrac{3\alpha}{2}\right)\right)^2-\cos^2\left(\tfrac{3\alpha}{2}\right))^{1/2}$.
Moreover, by \eqref{eq_EDS_RtA_2}, we have that on $E_1$, for all $t\in[\sigma_1,\psi_1[$,
\begin{align*}
\left(R^A_t\right)^2&= \left(R^A_{\sigma_1}\right)^2-2v\int_{\sigma_1}^t Y^A_s\,ds,
\end{align*}
and by the lower bound in \eqref{eq_proof_lem_YA_ineg}, 
\begin{align}\notag
\left(R^A_t\right)^2 &\leq\left(R^A_{\sigma_1}\right)^2+2v\int_{\sigma_1}^t \left(\sin\left(\tfrac{\alpha}{2}\right) R^A_{\sigma_1}+(s-\sigma_1) v\right)\,ds\\ \notag
&=\left(R^A_{\sigma_1}\right)^2+\left(\sin\left(\tfrac{\alpha}{2}\right)\,R^A_{\sigma_1}+(t-\sigma_1) v \right)^2-\sin^2\left(\tfrac{\alpha}{2}\right)\,\left(R^A_{\sigma_1}\right)^2\\
&\leq \left(R^A_{\sigma_1}\right)^2\left(1+\left(\sin\left(\tfrac{\alpha}{2}\right)\,+d_2\right)^2-\sin^2\left(\tfrac{\alpha}{2}\right)\right)= \left(R^A_{\sigma_1}\right)^2c_2^2, 
\label{eq_414a}
\end{align}
where we have used \eqref{eq_proof_lem_t-sigma}, the definition of $d_2$ and have set $c_2= (1+\left(d_2+\sin\left(\tfrac{\alpha}{2}\right)\right)^2-\sin^2\left(\tfrac{\alpha}{2}\right))^{1/2}$ (notice that the quantity under the square root is nonnegative).
Since the process $R^A_t$ is continuous, this inequality remains true for $t=\psi_1$.

The first inequality in \eqref{r_leq_r0} holds for $t\in [0,\zeta]$ by Lemma \ref{Lem_exist_A_zeta}. For $t\in [\zeta,\sigma_1]$, \eqref{r_leq_r0}  follows from \eqref{eq45b} since $\sin(\T^A_t) > 0$ by \eqref{eq_thetaA_zeta_tau1}.  Finally, for $t\in [\sigma_1,\psi_1]$, \eqref{rd10_13e1}, \eqref{eq_415a} and \eqref{eq_414a}  imply that 
\begin{align}\label{rd09_18e2}
R_t^A \leq c_1 c_2 r.
\end{align}
Direct computations show that
\begin{align*}
c_1^2c_2^2=\frac{1}{2}& \left(1+\cos\alpha +2 \left((1+\cos\alpha\sec(2\alpha)) \sin\left(\tfrac{\alpha }{2}\right)+\sec(2\alpha)\sin(\alpha)\right)^2\right)\\
&\times\left(1-\cos^2\left(\tfrac{3\alpha}{2}\right) +\left(1-\cos\left(\tfrac{3\alpha}{2}\right)+\sin\left(\tfrac{\alpha}{2}\right) \tan\alpha\right)^2\right).
\end{align*}
where $\sec x =\tfrac{1}{\cos x}$ (indeed, the second big factor in $c_1^2 c_2^2$ is equal to $c_1^2$, the quantity inside the square in the first big factor of $c_1^2 c_2^2$ simplifies to $\sin(\frac{\alpha}{2}) + d_2$, therefore the first big factor in $c_1^2 c_2^2$ is equal to $2 c_2^2$). Furthermore, it can be shown that $c_1^2$ and $c_2^2$ are  increasing and positive functions of $\alpha\in\left[0,\tfrac{\pi}{4}\right[$ and $c_1^2 c_2^2 \simeq 0.72$ when $\alpha = \tfrac{\pi}{8}$. Since we assumed that $\alpha\leq\tfrac{\pi}{8}$, we conclude that $c_1c_2 \leq \tfrac{3}{4} <1$; this establishes \eqref{r_leq_r0} and concludes the proof.
\end{proof}

The corollary below shows that if at any time during the race, the wind has a calm period ($E_i$ holds for some $i$), then the yacht reaches the origin with at most one additional tack.

\begin{cor}\label{cor_exist_taui_taua}
Let $\shs$ be the region  defined in \eqref{s_def}. For all $(r,\theta)\in \shs$,  if $E_i$ occurs for some $i \in\N^*$, then $\nu_i = \psi_i\leq\tau^A=\psi_{i+1}$  and $M_{\tau^A}(A)\leq i$, $\P_{r,\theta}$-a.s. Moreover, if $\psi_i < \tau^A$, then  $\T^A_{\psi_i} = \pim + \alpha$.
\end{cor}

\begin{proof}
By definition of the strategy $A$ and according to the statement just after \eqref{s_def}, tacks occur only at the times $\psi_i$, $i \in \N^*$. Suppose that $E_i$ occurs for some $i \in \N^*$. 
By construction, after each tack, the boat can be assumed to be in $\shs$ (see Remark \ref{rem_favorable_initial_conditions}), and in particular, we can assume that the boat is in $\shs$ at time $\psi_{i-1}$.
If $\psi_i=\tau^A$, then $\tau^A = \psi_i = \nu_i = \psi_{i+1}$, so we assume that $\psi_i<\tau^A$. This implies that we are in the situation described in Case 3 in the proof of Proposition \ref{Lem_exist_A_cst}, since $\psi_i=\tau^A$ in Cases 1 and 2 there. By 
\eqref{rd09_18e2}, the observation that $c_1 c_2 \leq \tfrac{3}{4}$ in Case 3 and by the strong Markov property, we have that $R^A_{\psi_i} <r_0$, $\P_{r,\theta}$-a.s. Thus, at time $\psi_i-$, the yacht hits the line $\theta=-\alpha$, $\nu_i<\rho_i$ (by \eqref{r_leq_r0} since $r<r_0$) and $\Theta^A_{\psi_i}=\tfrac{\pi}{2}+\alpha$. Again by the strong Markov property, this configuration corresponds to Case 2 in the proof of Lemma \ref{Lem_exist_A_cst}. Therefore, on $E_i$ and in this case, $\sigma_{i+1}=\nu_{i+1}=\psi_{i+1}=\tau^A$ and so $M_{\tau^A}(A)\leq i$. \end{proof}

\begin{cor}\label{cor_exist_E1_tauA}
Let $d_1$ and $d_2$ be as defined in \eqref{rd10_20e1}. If $E_1$ occurs, then for all initial positions $(r,\theta)\in \shs$,
\begin{equation}\label{rd09_08e1}
    \tau^A< \left(\frac{1}{\cos\alpha}+1 \right)\,\frac{r-\eta}{v} + (d_1+d_2)\frac{r}{v} =: \Gamma(r) ,\qquad \P_{r,\theta}\text{-a.s.}
\end{equation}
\end{cor}
\begin{proof}
On $E_1$, by Corollary \ref{cor_exist_taui_taua}, $M_{\tau^A}(A)\leq 1$. Consider the case where $\psi_1 < \tau^A$. By \eqref{r_leq_r0},  $\psi_1 = \nu_1$ and $\Theta^A_{\psi_1}=\tfrac{\pi}{2}+\alpha$. Thus, by Lemma \ref{Lem_exist_A_zeta} and Proposition \ref{Lem_exist_A_cst},
\begin{equation}\label{rd09_08e2}
\tau^A=\zeta+(\sigma_1-\zeta)+(\psi_1-\sigma_1)+(\tau^A-\psi_1)\leq \frac{r- \eta}{v\cos\alpha}+\frac{d_1\, r}{v}+\frac{d_2\, r}{v}+(\tau^A-\psi_1),\; \P_{r,\theta}\text{-a.s.}
\end{equation}
Furthermore, on $E_1$, as noted in Case 2 of the proof of Proposition \ref{Lem_exist_A_cst}, $\Theta^A_t\in\left[\tfrac{\pi}{2},\pi\right]$ for all $t\in[\psi_1,\tau^A]$, hence, as in \eqref{rd09_17e2}, 
$R^A_t=R^A_{\psi_1}-(t-\psi_1) v$, and setting $t = \tau^A$, we obtain
$$
   \tau^A-\psi_1=\frac{R^A_{\psi_1}-R^A_{\tau^A}}{v}=\frac{R^A_{\psi_1}-\eta}{v}<\frac{r-\eta}{v}.
$$
With \eqref{rd09_08e2}, this establishes \eqref{rd09_08e1} in the case where $\psi_1 < \tau^A$. 
The case $\psi_1 = \tau^A$ follows from \eqref{rd09_20e1}--\eqref{eq_Lem_cst_3}.
\end{proof}
\medskip

\noindent{\em Using calm spells to estimate the time to reach the target}
\medskip

\begin{prop}\label{prop_exist_p0}
There exists $p_0\in\,]0,1[\,$, depending only on the parameters $\alpha$, $\sigma$, $\eta$ and $r_0$, such that for all $(r,\theta)\in \shs$,

    (a)  $\P_{r,\theta}\left\{\psi_2=\tau^A\right\} \geq p_0$,
    
  (b) $\P_{r,\theta}\left\{\psi_i<\tau^A\right\} \leq (1-p_0)^{\lfloor \frac{i}{2} \rfloor}$ for all $i\in\N,\, i\geq2$ (where $\lfloor x \rfloor$ denotes the largest integer $n \leq x$).
\end{prop}
\begin{proof}
(a) By Corollary \ref{cor_exist_taui_taua}, $E_1\subset\left\{\psi_2=\tau^A\right\}$. Let $\Gamma(r)$ be as in \eqref{rd09_08e1}. By Corollary \ref{cor_exist_E1_tauA},
\begin{align*}
   \P_{r,\theta}\left\{\psi_2=\tau^A\right\} &\geq \P_{r,\theta}(E_1)=\P_{r,\theta}\left\{\sup_{0\leq t\leq\tau^A} \left|\sigma B_t\right|<\tfrac{\alpha}{2}\right\}\\
&\geq \P\left\{\sup_{0\leq t\leq \Gamma(r)} \left|\sigma B_t\right|<\tfrac{\alpha}{2}\right\} \geq \P\left\{\sup_{0\leq t\leq \Gamma(r_0)} \left|\sigma B_t\right|<\tfrac{\alpha}{2}\right\} =:p_0>0.
\end{align*}

(b) Notice that $\left\{\psi_i<\tau^A\right\}\subset\left\{\psi_{i-1}<\tau^A\right\}\subset\cdots\subset\left\{\psi_1<\tau^A\right\}$, so by the strong Markov property and (a), we get
\begin{align}\notag
\P_{r,\theta}\left\{\psi_i<\tau^A\right\}&=\P_{r,\theta}\left(\{\psi_i<\tau^A\}\cap\{\psi_{i-2}<\tau^A\}\right)\\ \notag
&=\E_{r,\theta}\left(\ind_{\{\psi_{i-2}<\tau^A\}}\,\P_{r,\theta}\left(\psi_i<\tau^A\mid\shf_{\psi_{i-2}}\right)\right)\\
&\leq \P_{r,\theta}\{\psi_{i-2}<\tau^A\}\, (1-p_0).
\label{ind_p}
\end{align}
If $i$ is even, then the induction is stopped at $\P_{r,\theta}\{\psi_{2}<\tau^A\}$  and the claim follows. If $i$ is odd, then the induction is stopped at $\P_{r,\theta}\{\psi_{1}<\tau^A\} \leq 1$  and the conclusion follows. 
\end{proof}

\begin{prop}\label{prop_NtauA_tauA_finite}
For all $(r,\theta)\in  \shs$,

 (a)  $\E_{r,\theta}\left(M_{\tau^A}(A)\right)<  \frac{2}{p_0} - 1 $;
 
 (b)  $\E_{r,\theta}\left(\tau^A\right)<K \frac{2}{p_0}$, where $K$ is the constant defined in \eqref{rd09_16e3}.
\end{prop}

\begin{proof}
(a) Since $t \mapsto M_t(A)$ counts the number of tacks and tacks only occur at the times $\psi_1, \psi_2, \dots$ when $\psi_i < \tau^A$, we have $\{M_{\tau^A}(A)>i\}\subset\{\psi_{i+1}<\tau^A\}$, therefore, by Proposition \ref{prop_exist_p0},
\begin{align*}
\E_{r,\theta}\left(M_{\tau^A}(A)\right)&=\sum_{i=0}^\infty\P_{r,\theta}\left\{M_{\tau^A}(A)>i\right\} \leq\sum_{i=0}^\infty\P_{r,\theta}\left\{\psi_{i+1}<\tau^A\right\} \\
&\leq 1+\sum_{i=2}^\infty\P_{r,\theta}\left\{\psi_i<\tau^A\right\} \leq 1+\sum_{i=2}^\infty(1-p_0)^{\lfloor \frac{i}{2} \rfloor}= \frac{2}{p_0} - 1.
\end{align*}

(b) We have
\begin{align*}
\E_{r,\theta}\left(\tau^A\right)&=\E_{r,\theta}\left(\sum_{i=0}^{M_{\tau^A}(A)}u_{i+1}\right)=\E_{r,\theta}\left(\sum_{k=0}^\infty\sum_{i=0}^k u_{i+1}\ind_{\left\{M_{\tau^A}(A)=k\right\}}\right)\\
&=\sum_{i=0}^\infty\E_{r,\theta}\left(\sum_{k=i}^\infty u_{i+1}\ind_{\left\{M_{\tau^A}(A)=k\right\}}\right)=\sum_{i=0}^\infty\E_{r,\theta}\left( u_{i+1}\ind_{\left\{M_{\tau^A}(A)\geq i\right\}}\right)\\
&\leq\sum_{i=0}^\infty\E_{r,\theta}\left( \E_{r,\theta}\left( u_{i+1}\mid \shf_{\psi_i}\right)\ind_{\left\{\psi_i<\tau^A\right\}}\right).
\end{align*}
By the strong Markov property, Lemma \ref{Lem_exist_A_ui} and Proposition \ref{prop_exist_p0}, we find that the right-hand side is bounded above by
$$
   K \sum_{i=0}^\infty\P_{r,\theta}\left\{\psi_i<\tau^A\right\}\leq  K\left(2+ \sum_{i=2}^\infty(1-p_0)^{\lfloor \frac{i}{2} \rfloor}\right) =K \frac{2}{p_0},
$$
and the proof of Item (b) is complete.
\end{proof}

The theorem below is the main result of the section.  As a side result, it also establishes an upper bound on the value function \eqref{Veta}  for the case $c=0$,  which will be improved upon in the next section, where a candidate optimal strategy  is discussed. We recall that $V^{c,\eta}(r,\t,a)$ (resp. $V^{\eta}(r,\t)$) is defined in \eqref{Vetac} (resp. \eqref{Veta}).

\begin{thm}\label{thm_brw_exist_V}
For every $\eta \geq 0$, $c\geq 0$ and every $(r,\t,a)\in\R_+ \times \R \times\{+1, -1\}$, $A$ belongs to $\sha^\eta_{r,\t}  \cap \sha^{c,\eta}_{r,\t,a}$.
In particular, 

 (a) $V^{c,\eta}(r,\t,a) \leq  K (1+\tfrac{2}{p_0}) + c \left( 1+ \frac{2}{p_0} \right)$;
 
   (b) $V^{\eta}(r,\t)\leq  K (1+\tfrac{2}{p_0})$,

\noindent where $K$ is the constant defined in Lemma \ref{Lem_exist_A_ui}.
\end{thm}

\begin{proof}
Let $r_0 > r$. Proposition \ref{prop_NtauA_tauA_finite} shows that the strategy $A$ constructed in \eqref{rd09_15e1}-\eqref{A_tau_def} satisfies 
 \begin{equation}
 J^{c,\eta}(r,\t, a,A) \leq K \frac{2}{p_0} + c \left( \frac{2}{p_0} - 1\right),
 \end{equation}
  for all $(r,\theta)\in \shs_{\eta, r_0, \alpha}$ and $a=1$. 
 Remark \ref{rem_favorable_initial_conditions} tells us that if the initial conditions do not belong to $\shs_{\eta, r_0, \alpha}$, that is, if $a=1$ and $(r,\t) \in C^1$,
 then after one tack, which takes an expected time bounded by $K$, the boat reaches a position in $\shs$, so in this case, $J^{c,\eta}(r,\t, a,A) \leq K (1+\tfrac{2}{p_0}) + c \left( \frac{2}{p_0} \right)$; or if $a = -1$ and $(r,\t)\in C^{-1}$, then the symmetry \eqref{eq_brw_polar_symmetry} tells us that it is equivalent to consider $a=1$ and $(r,\pim-\t)$, which is in $C^1$.
Finally, if $a= 1$ and $(r,\t)\in D^1$ (resp. $a= -1$ and $(r,\t)\in D^{-1}$), then the yacht tacks immediately, so  $J^{c,\eta}(r,\t, a,A) \leq K (1+\tfrac{2}{p_0}) + c \left( 1+ \frac{2}{p_0} \right)$.
This proves (a). Together with Lemma \ref{strong_sol}, this also shows that $A\in \sha^{c,\eta}_{r,\t,a} $, for all $\eta \geq 0$, $(r,\t) \in \R_+ \times \R$, $a\in \{+1,-1\}$.  Because $\sha^{c,\eta}_{r,\t,a} \subset \sha^\eta_{r,\t}$, it also follows that $ A \in \sha^\eta_{r,\t}$. 
Item (b)  follows from (a) after setting $c = 0$.
\end{proof}

\section{Candidate optimal strategy when $\eta> 0$ and $c=0$}\label{Sec:5}

The strategy $A$ studied in the previous section was used to show that the value functions are finite in both  cases $c >0$ and $c=0$. However, $A$ is not optimal. Indeed, consider for example the situation where $\T_t \in [0,\piq[$ 
and $A_t = 1$, which may occur with positive probability. In this situation, the boat would move towards the buoy at speed $v\sin(\T_t)$, while if $A_t$ were switched to $-1$, then the speed would become $v\cos(\T_t) > v\sin(\T_t)$. In the case $c=0$ (or $c$ small), this improvement in speed can be achieved without incurring any (or only little) additional cost.

We hypothesize that the optimal strategy in the case $c=0$ involves a feedback approach that selects, at any given position in the race field, the tack that maximizes the projection of the boat's speed along the radial direction.
We  prove that this choice leads to an admissible strategy in the sense of Definition \ref{Def:adm_con}.
 However, the coefficients of the state equation turn out not to be Lipschitz and standard results cannot be applied. Nevertheless, the coefficients do exhibit a monotonicity property, and this makes it possible  to prove strong uniqueness and to apply the classical Yamada-Watanabe argument \cite[Chapter 5, Section D]{ks}. 

\subsection{Defining the strategy}
\begin{figure}
\centering
\includegraphics[height=4cm, width=4cm,keepaspectratio]{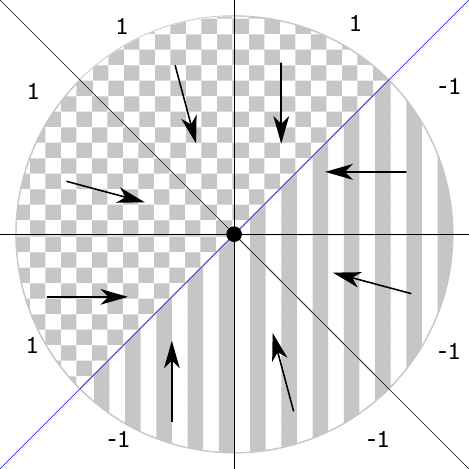}\hspace{5ex}\includegraphics[height=5cm,width=6cm,keepaspectratio]{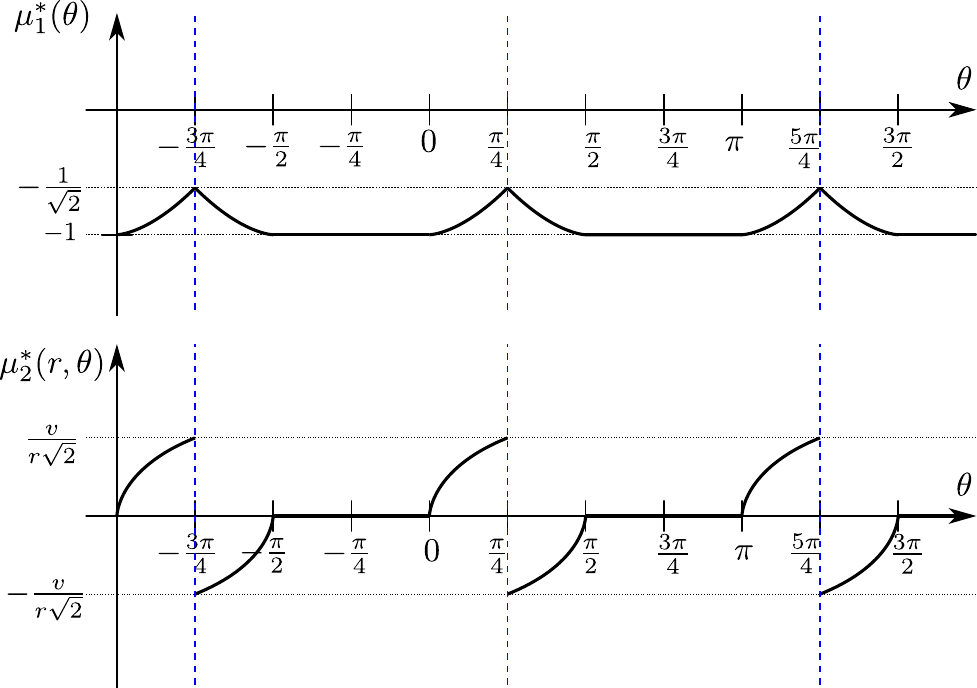}
\caption{On the left, the candidate optimal strategy $\sa$. The arrows indicate the direction of the prescribed tack.  
On the right, the drifts $\mu^\star_1$ and $\mu^\star_2(r, \cdot)$ of the processes $R^{\sa}$ and $\T^{\sa}$ extended to $\R$ by $2\pi$-periodicity. }
\label{Fig:9}
\end{figure}

Let $\eta > 0$ and fix the starting position $(r_0,\t_0)$. Our candidate optimal strategy $(A^\star_t)$ has the following form: when the boat is in the region $S = \{(r, \t): \t \text{ mod } 2\pi \in  [0, \piq[\, \cup\, ]\tfrac{5 \pi}{4}, 2\pi[ \}$, then the boat sails on starboard tack ($\sa_t = +1$), whereas when the boat is in region $P = \{(r, \t): \t \text{ mod } 2\pi \in \, ]\piq, \tfrac{5 \pi}{4} [\}$, then the boat sails on port tack ($\sa_t = -1$); see Fig.~\ref{Fig:9}. Of course, the tack is not specified on the main diagonal, but since the angle process $(\T^{A^*}_t)$ will satisfy a diffusion equation as in \eqref{state1}, the boat should spend zero (Lebesgue-measure) amount of time on this diagonal, on which a rather complex switching process will occur. This situation is comparable to the switching procedures that appear in \cite{DS,Ma87,MSV90}. 

The case $\eta = 0$ requires a somewhat different analysis: in particular, the damping procedure described in \eqref{phi}-\eqref{state} would not be applicable (however, see \cite{CC2}).

 Let $A^\star_{0-}(\theta) = a_0$ for a given $a_0 \in \{-1, +1 \}$. However, recall that since we are in the setting $c=0$,  the initial tack is not relevant. 
 
 Formally,  for $t \geq 0$, $\sa = (\sa_t)$ is defined as the feedback strategy $A_t^\star(\omega) := a^\star(\Theta_t(\omega))$, where 
\begin{equation}\label{eq50}
a^\star(\t) =
\left\{
\begin{array}{rll}
 1, &\text{if } \t \in [\frac{\pi}{4}, \frac{5\pi}{4}[,  \\[6pt]
 -1, &\text{if }\t \in [-\frac{3\pi}{4}, \frac{\pi}{4}[, 
\end{array}
\right.
\end{equation}
extended to $\R$ by $2\pi$-periodicity (see  Fig.~\ref{Fig:9}).

Recall that $\mu_1(\theta, a)$ and $\mu_2(r, \theta, a)$ are defined in \eqref{eq_drift_ra}--\eqref{eq_drift_pa}. Under the strategy $A^\star$,  the system \eqref{state1} becomes
\begin{equation*}
\left\{\begin{array}{l l}
   dR^{\sa}_t =  1_{\{R_t^{\sa} > \eta\}}\, \mu_1(\Theta^{\sa}_t, a^\star(\Theta^{\sa}_t ))dt, \\
   d\Theta^{\sa}_t =  1_{\{R_t^{\sa} > \eta\}}\left[ \mu_2(R^{\sa}_t, \Theta^{\sa}_t, a^\star(\Theta^{\sa}_t )) dt + \sigma  dB_t \right],
\end{array}
\right.
\end{equation*}
or, equivalently
\begin{equation}\label{state2}
\left\{\begin{array}{l l}
dR^{\sa}_t =  1_{\{R_t^{\sa} > \eta \}} \,\mu^\star_1(\Theta^{\sa}_t)dt, \\
d\Theta^{\sa}_t =  1_{\{R_t^{\sa} > \eta \}} \left[ \mu^\star_2(R^{\sa}_t, \Theta^{\sa}_t) dt + \sigma  dB_t\right],
\end{array}
\right.
\end{equation}
where
\begin{equation}\label{mu_1_s}
 \mu^\star_{1}(\t) := \left\{\begin{array}{l}
      \mu_1(\t,1),    \\
       \mu_1(\t,-1),
       \end{array}\right.
       \quad \mu^\star_{2}(r,\t):=\left\{\begin{array}{lll}
           \mu_2(r,\t,1), & \qquad\text{if } \t \in [\piq,\picq[, & r > 0, \\
           \mu_2(r,\t,-1), & \qquad  \text{if } \t \in[-\pitq,\frac{\pi}{4}[, & r > 0.
           \end{array}\right.
\end{equation}
 The functions $\mu_1^\star$ and $\mu_2^\star(r, \cdot)$ are then extended to all  $\R$ by $2\pi$-periodicity. They turn out to be $\pi$-periodic.

Notice that $\t \mapsto \mu_2^\star(r, \t)$ is not continuous at $\piq + k \pi$, $k \in \Z$ (for instance, $\mu_2^\star(r, \piq +) = \mu_2^\star(r, \piq) = -\frac{v}{r \sqrt{2}} \neq  \frac{v}{r \sqrt{2}} = \mu_2^\star(r, \piq -)$), and because of the diffusion term $\sigma 1_{\{R_t^{\sa} >0\}} dB_t$, if $\T^{\sa}_{t_0} = \piq$ for some $t_0 \geq 0$, then $\T^{\sa}_t$ will again hit $\piq$ infinitely many times during $[t_0, t_0 + \epsilon]$, for any $\epsilon > 0$. The existence of a strong solution to \eqref{state2} is therefore {\em not} covered by any standard theorem on systems of SDEs.

\begin{rem}
The interpretation of the strategy $A^\star$ is as follows. In a long-distance race in which the tacking cost is negligible relative to the total travel time, when sailing upwind or downwind the boat should choose the heading that takes it towards the most distant layline, so as to remain as far as possible from the laylines. When sailing in crosswinds, it should aim for the target.
\end{rem}

\subsection{Admissibility of the controlled process}

As mentioned above, the set of discontinuities of the function $\t \mapsto \mu_2^\star(r,\t)$ is
$$ 
   G := \{g_k,\, k\in \Z\}, \qquad \text{where }\quad g_k= \piq+k \pi,
 $$
 therefore  standard existence results do not apply. The coefficients do satisfy a monotonicity property (see \eqref{mon} below), but even the results in literature which rely on this property assume continuity of the coefficients in the spatial coordinates and hence are not applicable (see for example \cite[Sec.~3]{KR1} or \cite[Theorem 3.1.1]{rl}). Additionally, the diffusion is degenerate, so the classical results of \cite{veretennikov1982} cannot be used either. Our proof uses the classical Yamada-Watanabe procedure \cite[Chapter 5, Section D]{ks}, therefore we first prove pathwise uniqueness, and then establish the existence of a weak solution.

\begin{thm}\label{Solution}
Let  $\eta > 0$. For every $(r,\t)\in \, ]\eta, \infty[ \times \R$, there exists a unique strong solution to \eqref{state2} on $[0,\infty [$ with initial condition $(R^{\sa}_0, \T^{\sa}_0) = (r, \t)$.
Moreover, $\tau^{\sa}_\eta \leq \sqrt{2}\, (r - \eta)/v$, $\P_{r, \t}$-a.s. In particular,   $\sa$ is an admissible strategy in $\sha^\eta$.
\end{thm}

\begin{proof}
The coefficients in the system \eqref{state2} are not continuous functions of $r \in [0, \infty[$, and $r \mapsto \mu^\star_2(r, \t)$ diverges when $r \downarrow 0$ for many values of $\t$. 
For this reason, and since $\eta > 0$, we first  study a system in which the coefficients are damped near the origin:   Define $n_\eta:= \inf \{n \geq 1: \tfrac{1}{n} <\eta \}$, and for $n \geq n_\eta$,
\begin{equation}\label{phi}
\phi_n(r) :=\left\{\begin{array}{ll}
0, &\qquad \text{if } r \in [0,  \eta -\tfrac{1}{n}],\\
n(r -\eta)  +1, &\qquad  \text{if } r \in\, ]\eta -\tfrac{1}{n},\eta[,\\
1, &\qquad   \text{if } r \geq \eta.
\end{array}
\right.
\end{equation}
Consider the system of SDEs
\begin{equation}\label{state}
\left\{\begin{array}{l l}
dR^{n,\sa}_t =  \phi_n(R^{n,\sa}_t) \mu^\star_1(\Theta^{n,\sa}_t)dt,\\
d\Theta^{n,\sa}_t = \phi_n(R^{n,\sa}_t)\left[\mu^\star_2(R^{n,\sa}_t, \Theta^{n,\sa}_t) dt +  \sigma dB_t\right],
\end{array}
\right.
\end{equation}
with the initial condition $(R^{n,\sa}_0, \Theta^{n,\sa}_0) = (r, \t)$.
Define the stopping time $ \tau^{n,\sa}_\eta:= \inf\{t \geq 0: R^{n,\sa}_t \leq  \eta  \}.$ Since $\phi_n(r)=1$ for $r\geq \eta$, for $m > n \geq n_\eta$, it follows that $(R^{m, \sa}_t, \T^{m, \sa}_t) = (R^{n, \sa}_t, \T^{n, \sa}_t)$ when $t \leq \tau^{n,\sa}_\eta$, so $\tau^{m,\sa}_\eta = \tau^{n,\sa}_\eta$. 
Therefore, if there exists a unique strong solution $((R^{n, \sa}_t, \T^{n, \sa}_t))_{  t \in \R_+}$ to \eqref{state}  for any $n \geq n_\eta$, 
then the solution to \eqref{state2} will satisfy $\tau^{\sa}_\eta = \tau^{n,\sa}_\eta$ $\P_{r, \t}$-a.s.

First, for fixed $n \geq n_\eta$, we will  prove the existence of a strong solution $((R^{n, \sa}_t, \T^{n, \sa}_t))_{t \in \R_+}$ to \eqref{state}. We note that the coefficients of the system \eqref{state} are bounded and continuous on $[0, \infty[ \times \R$, and because $\phi_n(r) = 0$ for $0 \leq r < \eta - 1/n$, the first equation in \eqref{state} is such that any solution to \eqref{state} will remain in $[\eta - 1/n, \infty[ \times \R$ for all time.

By a standard localization procedure, this will imply the existence and uniqueness of a strong solution $((R^{\sa}_t, \T^{\sa}_t))_{ t \in \R_+}$ to  \eqref{state2}, 
since we can set $(R_t^{\sa}, T_t^{\sa}) = (R_{t \wedge \tau^{\sa}_\eta}^{n,\sa}, \T^{n,\sa}_{t \wedge \tau^{\sa}_\eta})$, $t\in \R_+$.

 Fix $n \geq n_\eta$. In the sequel, we will write $\mu_1(\t)$ for $\mu^\star_1(\t)$, $\mu_2(r,\t)$ for $\mu^\star_2(r,\t)$ and $\phi(r)$ for $\phi_n(r)$. The first step is to prove that $\mu_1$ and $\mu_2$ satisfy, for all $r_1, r_2 > 0$ and $\t_1, \t_2 \in \R$, the monotonicity condition 
\begin{align}\notag
&\left[ \begin{pmatrix} r_1 \\ \t_1 \end{pmatrix}  -\begin{pmatrix} r_2\\ \t_2 \end{pmatrix} \right] \cdot \left[ \begin{pmatrix} \phi(r_1)\mu_1(\t_1) \\ \phi(r_1)\mu_2(r_1,\t_1) \end{pmatrix}-  \begin{pmatrix} \phi(r_2)\mu_1(\t_2) \\ \phi(r_2)\mu_2(r_2,\t_2) \end{pmatrix}\right] \\
&\qquad \leq C_n [(r_1-r_2)^2 + (\t_1 - \t_2)^2],
\label{mon}
\end{align}
where $C_n$ is a constant that depends on $n$,  which may vary from one line to the other. 
In fact, \eqref{mon} can be rewritten  
\begin{align}\notag
   &(r_1-r_2)[\phi(r_1) \mu_1(\t_1)-\phi(r_2)\mu_1(\t_2)]  + (\t_1-\t_2) [\phi(r_1) \mu_2(r_1,\t_1)-\phi(r_2)\mu_2(r_2,\t_2)] \\ \notag
   &\quad= (r_1-r_2)[\phi(r_1)  -\phi(r_2)]\mu_1(\t_1)  + (r_1-r_2)\phi(r_2) [\mu_1(\t_1) - \mu_1(\t_2)] \\ \notag
     &\qquad\qquad + (\t_1-\t_2) [\phi(r_1) \mu_2(r_1,\t_1)-\phi(r_2)\mu_2(r_2,\t_1)] \\
   & \qquad\qquad + (\t_1-\t_2) \phi(r_2) [\mu_2(r_2,\t_1) - \mu_2(r_2,\t_2)].
\label{eq_5_5}
\end{align}
Using the fact that $r \mapsto \phi(r)$ and $\t \mapsto \mu_1(\t)$ are Lipschitz on $\R$, and $r \mapsto \phi(r) \mu_2(r,\t)$ is Lipschitz on $\R_+$ (with Lipschitz constant that depends on $n$ but not on $\t$), we deduce that the right-hand side of \eqref{eq_5_5} is bounded above by 
\begin{equation}\label{eq_5_6}
 C_n\, \mu_1(\t_1)(r_1-r_2)^2+ C_n\,  |r_1-r_2|| \t_1-\t_2| + \phi(r_2)(\t_1-\t_2) [ \mu_2(r_2,\t_1)-\mu_2(r_2,\t_2)].
\end{equation}

Recall the elementary inequality $ab\leq \um (a^2+b^2)$. In order to get the desired inequality \eqref{mon}, it suffices to prove that there exists a positive constant $C_n$ such that for all $r > 0$ and $\t_1, \t_2 \in \R$,
\begin{align}\label{eq:2.6a}
    \phi(r) (\t_1-\t_2) [\mu_2(r,\t_1) -\mu_2(r,\t_2) ] \leq C_n\, (\t_1-\t_2)^2  .
\end{align}

For this, let $F_n(r, \t_1,\t_2) := \phi(r) (\t_1- \t_2)[\mu_2(r,\t_1)-\mu_2(r,\t_2)]$. First observe that $F_n(r, \t_1,\t_2)=F_n(r, \t_2,\t_1)$, so it is enough to consider the case $\t_1 < \t_2$.
We claim that in fact, it is enough to prove inequality \eqref{eq:2.6a} for $\t_1$ close to $\t_2$. 
Indeed, assume that for some $\d >0$, \eqref{eq:2.6a} holds for $r > 0$ and for some constant $C_n$ when $0 < \t_2 -\t_1 \leq \d$. Let  $M_n := \sup_{r > 0,\, \t \in \R} |\phi(r) \mu_2(r, \t)|$, and set $D_n:=\frac{2M_n}{\d}\vee C_n$. This implies that for $\t_2-\t_1 \geq \d$ and $r > \eta - 1/n$,
\begin{align*}
   \phi(r) (\t_1-\t_2)(\mu_2(r,\t_1)-\mu_2(r,\t_2)) \leq |\t_1-\t_2|\, 2M_n \leq (\t_2-\t_1) D_n\,  \d \leq D_n(\t_2-\t_1)^2.
\end{align*}
Now, it suffices to replace $C_n$ by $D_n$ in \eqref{eq:2.6a}.

We  now prove that \eqref{eq:2.6a} holds for $0 < \t_2 -\t_1 \leq \piq$ (i.e.~$\delta = \piq$).
If $g_k \notin \, ]\t_1,\t_2]$, for any $k\in \Z$, then $\t \mapsto \phi(r) \mu_2(r,\t)$ is Lipschitz on $[\t_1, \t_2]$ (and the Lipschitz constant $C_n$ depends on $n$ but not on $r$) and therefore,
\begin{align*}
 \phi(r) (\t_1-\t_2)(\mu_2(r,\t_1)-\mu_2(r,\t_2)) &=  \phi(r) (\t_2-\t_1)(\mu_2(r,\t_2)-\mu_2(r,\t_1)) \\
 & \leq C_n(\t_1-\t_2)^2.
\end{align*}
Thus, \eqref{eq:2.6a} follows in this case.

Now let $\t_1,\t_2$ be such that $g_k \in\, ]\t_1,\t_2]$, for some $k \in \Z$. Since the function $\t \mapsto \mu_2(r,\t)$ is periodic, it suffices to consider angles $\t_1,\t_2$ such that  $G_0 = \piq \in\, ]\t_1,\t_2]$, and $0 < \t_2 -\t_1 \leq \piq$. 
This corresponds to $\t_1 \in [0,\piq[$ (port tack) and $\t_2 \in [\piq,\pim[$ (starboard tack).  By \eqref{eq_drift_ra} and \eqref{eq_drift_pa},  the right hand-side of \eqref{eq:2.6a} becomes 
$\phi(r) (\t_1-\t_2) v \left(\frac{\sin \t_1}{r}+\frac{\cos \t_2}{r}\right)\leq 0$  (because the second factor is negative and the first, third and fourth factors are nonnegative). This establishes inequality \eqref{eq:2.6a} and completes the proof of \eqref{mon}. 

In order to prove the existence and uniqueness of a strong solution, we first prove pathwise uniqueness and then the existence of a weak solution. Let $T>0$, $t\in [0,T]$ and  $n > n_\eta$ be fixed. 
Let $(R_t,\Theta_t)$ and $(R'_t,\Theta'_t)$  be two solutions of the system \eqref{state}.   By It\^{o}'s Lemma for continuous diffusion processes \cite[Chapter 3, Theorem 3.3]{ks},
\begin{align}\notag 
&\E \left((R_t-R'_t)^2 + (\T_t-\T'_t)^2 \right) \\ \notag
&\qquad= 2 \E \Big( \int_0^t [(R_s-R'_s) (\phi(R_s)\mu_1(\T_s)- \phi(R'_s)\mu_1(\T'_s)) \\ \notag
&\qquad\qquad\qquad\qquad  + (\T_s-\T'_s)(\phi(R_s)\mu_2(R_s,\T_s)-\phi(R'_s)\mu_2(R'_s,\T'_s))] ds \Big) \\ \notag
&\qquad\qquad\quad  + 2 \E \left( \int_0^t  (\T_s-\T'_s) (\phi(R_s) - \phi(R'_s))\s dW_s\right)  \\
&\qquad\qquad\quad + \s^2\,  \E \left( \int_0^t   (\phi(R_s) - \phi(R'_s))^2  ds \right).
\label{eq:2.6}
\end{align}
Using the monotonicity property \eqref{mon} in the first integral, the Lipschitz property of $\phi$ in the third integral and the fact that the stochastic integral is a martingale (which we check below) with mean zero, the right-hand side of \eqref{eq:2.6} is bounded above by
\begin{align}\label{eq:2.6aab}
& C \int_0^t  \E   \left((R_s-R'_s)^2 +(\T_s-\T'_s)^2 \right) ds ,
\end{align}
where $C$ is a positive constant. From \eqref{eq:2.6}, \eqref{eq:2.6aab} and Gronwall's inequality \cite[Chapter 5, Problem 2.7]{ks}, we conclude that for all $t \in [0,T]$,
\begin{align}
\E \left((R_t-R'_t)^2 + (\T_t-\T'_t)^2 \right)= 0,
\end{align}
and since the processes $(R_t, \T_t)$ and $(R'_t, \T'_t)$ are continuous, pathwise uniqueness holds.

 In order to prove the martingale property of the stochastic integral, since $\phi$ is bounded, it suffices to show that
\begin{equation}\label{391}
\E \left( \int_0^T \T_t^2 \, dt \right) < \infty.
\end{equation}
From \eqref{state}, it follows that
\begin{align*}
   \E\left(\T_t\right)^2 &\leq  3 \t_0^2 + 3 \E\left[\left(\int_{0}^{t} \phi(R_s)\mu_2(R_s,\T_s) ds\right)^2\right] + 3 \E\left[ \left(\s \int_{0}^t \phi(R_s)dW_s\right)^2\right] \\
   &\leq  3 \t_0^2 + 3 t^2 M^2 
       + 3\s^2\,  \E \left( \int_{0}^t (\phi(R_s))^2 ds \right) 
    < \infty,
\end{align*}
where $M$ is a bound for the function 
$(r, \t) \mapsto \vert \phi(r) \mu_2(r, \t)\vert$, $(r, \t) \in [0, \infty[ \times \R$, and  we have used It\^{o}'s isometry in the stochastic integral (which is possible because $r \mapsto \phi(r)$ is bounded). This shows that \eqref{391} is satisfied.

We now prove the existence of a weak solution on $[0,T]$, for any $T>0$. Let $(B_t)_{t\in [0,T]}$ be a standard Brownian motion on $(\Omega,\shf,(\shf_t),\P_{r,\t})$. Then the following system has Lipschitz and bounded coefficients: 
\begin{equation}\label{weak}
\left\{\begin{array}{l l}
dR_t = \phi(R_t)  \mu_1(\Theta_t)dt,\\
d\T_t =  \sigma dB_t.
\end{array}
\right.
\end{equation}
The existence of a unique strong solution to \eqref{weak}, for $t\in [0,T]$, follows from standard results, see for example \cite[Chap.~5, Theorem 2.9]{ks}.
Since the function $\phi(r)\mu_2(r,\t)$ is bounded  for $(r,\t)\in [0,\infty[\times \R$,  Novikov's condition \cite[Chap.~3, Corollary 5.13]{ks} implies that  the process
\begin{align*}
Z_t= \exp\left\{ -\frac{1}{2} \int_0^t \frac{[\phi(R_t) \mu_2(R_s,\Theta_s)]^2}{\sigma}  ds  +  \int_0^t \frac{\phi(R_t)\mu_2(R_s,\T_s)}{\sigma} dB_s \right\} ,\qquad t\in [0,T],
\end{align*}
is a martingale.
Let $\Q_T$ be the probability measure on $(\Omega, \shf_T)$ defined by $d\Q_T = Z_T\, d\P_{r,\t}$. By Girsanov's theorem \cite[Chap.~3, Theorem 5.1]{ks}, the process $(B'_t)_{t \in [0, T]}$ defined by $B'_0 = 0$ and
\begin{align*}
dB'_t  = dB_t -  \int_0^t \phi(R_t)\mu_2(R_s,\Theta_s)\,  ds, \qquad t\in [0,T], 
\end{align*}
is a   $(\Omega,\shf, (\shf_t),\Q_T)$ Brownian motion. The system \eqref{weak} can be written
\begin{equation}\label{weak1}
\left\{\begin{array}{l l}
dR_t =  \phi(R_t) \mu_1(\Theta_t)dt,\\
d\Theta_t = \phi(R_t) \mu_2(R_t,\Theta_t)  + \sigma dB'_t,
\end{array}
\right.
\end{equation}
so $(R_t,\T_t)$ is a weak solution to \eqref{state} in $[0,T]$.
Weak existence and pathwise uniqueness imply that there exists a unique strong solution to \eqref{state} in $[0,T]$ (see \cite[Chap.~5, Corollary 3.23]{ks}). Since $T$ is arbitrary, we obtain a strong solution $((R_t, \T_t))_{ t \in \R_+}$ of \eqref{state}.

It remains to prove  that $\tau^{\sa}_\eta<\sqrt{2}\, (r- \eta)/v$ for every $\eta  > 0$.  Indeed, the process $R_t$ satisfies the equation
\begin{align*}
R_{\tau^{\sa}_\eta} = r + \int_0^{\tau^{\sa}_\eta} \phi(R_t) \mu_1(\Theta_t) dt.
\end{align*}
If $R_t \geq \eta$, then $\phi(R_t) = 1$ and, by \eqref{mu_1_s}, \eqref{eq_drift_ra} and \eqref{eq_drift_pa}, it follows that  $\mu_1(\t) \leq  -\frac{v}{\sqrt{2}}$ for every $\theta \in \R$.
 Therefore
\begin{align*}
\eta= R_{\tau^{\sa}_\eta} = r + \int_0^{\tau^{\sa}_\eta} \mu_1(\Theta_t)\, dt \leq r- \frac{v}{\sqrt{2}}\, \tau^{\sa}_\eta,
\end{align*}
so
\begin{equation}\label{778}
\tau^{\sa}_\eta \leq \sqrt{2}\, \frac{r- \eta}{v} .
\end{equation}
This completes the proof of Proposition \ref{Solution}.
\end{proof}

\begin{rem}
For a boat starting on the main diagonal ($\theta = \piq $) at distance $r> \eta$ from the origin and traveling in a constant wind, all piecewise constant strategies will take $\sqrt{2}\, (r-\eta)/v$ to reach the target (see Fig.~\ref{Fig_points_of_sail}). Therefore, we expect that the upper bound given in Theorem \ref{Solution} should be equal to the value function in this case.  
\end{rem}

\section{Conclusion and future work}

In this paper, we have developed a stochastic control framework for modeling upwind sailboat navigation and have analyzed two complementary problems: the first involves impulse controls with a positive tacking cost, and second involves singular controls with zero tacking cost. We have established the finiteness of the corresponding value functions by identifying a strategy that guarantees reaching the target in a finite number of tacks, have obtained explicit bounds on the expected travel time under this strategy, and have identified a candidate optimal strategy in the case of zero tacking cost. This candidate strategy suggests that in a long-distance race, the boat should navigate so as to remain as far as possible from the (rotating) laylines.

Beyond its nautical motivation, the model leads to a broad class of stochastic control problems involving discontinuous drifts and degenerate diffusions, and provides a rich framework for both applied and theoretical extensions.
\bigskip

\noindent{\em Ongoing and future work}
\medskip

Two companion papers further develop the mathematical aspects of this work. 

\begin{itemize}
    \item In \cite{CC1}, we rigorously establish the \emph{optimality of the candidate strategy} introduced in Section \ref{Sec:5} for the case $c=0$ and we provide detailed information on the corresponding value function.
    \item In \cite{CC2}, we address the \emph{limiting case $\eta = 0$} with $c=0$, where the drift becomes supercritical near the origin, and we prove the existence of a strong solution to the associated stochastic differential equation.
\end{itemize}


There are many other questions that should be addressed. For instance, a complete solution to the free boundary problem described in Section \ref{appB} for the case $c>0$ remains to be obtained, along with a description of the manner in which this free boundary depends on the parameter $\sigma$. 

Modeling the wind as a Brownian motion on the unit circle without drift is the most fundamental situation. However, the wind direction may exhibit a trend, and in this case, the drift should be non-zero. This modifies the problem substantially, because it is known from experience and from numerical results obtained in \cite{DUM}, that in certain situations, the boat should navigate at an angle different from the angle $\alpha_0$ discussed in Section \ref{sec2.1}. In the case $c > 0$, the stochastic control problem is no longer an optimal switching problem, and new difficulties appear.

A substantial extension would be to model the wind as a random field that depends on both time and space. Indeed, sailors often seek to guess what the wind will be like when they have advanced some distance to one side or other of the race field, and it would be quite useful to be able to incorporate this information into a more advanced model.

Altogether, our model provides a unified framework for mathematical studies of sailboat navigation within which many interesting stochastic stochastic control problems related to sailboat navigation can be posed and solved.

\section{Appendix A: some unusual wind behaviors}\label{appA}

We briefly discuss wind behaviors that could prevent the boat from reaching the target in a finite number of tacks. These behaviors are not particularly realistic, but in a mathematical model in which all cases must be covered, they cannot be ignored. In our model with Brownian wind on a circle, these unusual behaviors have probability zero.

Imagine first that a malicious adversary controls the wind. There are many ways that this adversary can control the wind's direction over time so as to prevent the boat from reaching the origin (which is the target when $\eta = 0$) in a finite number of takings. For instance, imagine that the wind is coming from the North, the boat is navigating upwind on starboard tack and is approaching the port layline; its heading is North-West. Since the wind is variable, the boat will overshoot the layline by some safety margin before tacking onto port. Almost immediately after this tack, the adversary swiftly shifts the wind 90 degrees to the right, so that the wind is now from the East. This forces the boat, which is now on port tack, to head South-East, parallel to its initial path but in the opposite direction! This scenario can be repeated indefinitely, preventing the boat from reaching the target with a finite number of tacks.

\begin{figure}
\begin{center}
\includegraphics[height=10cm]{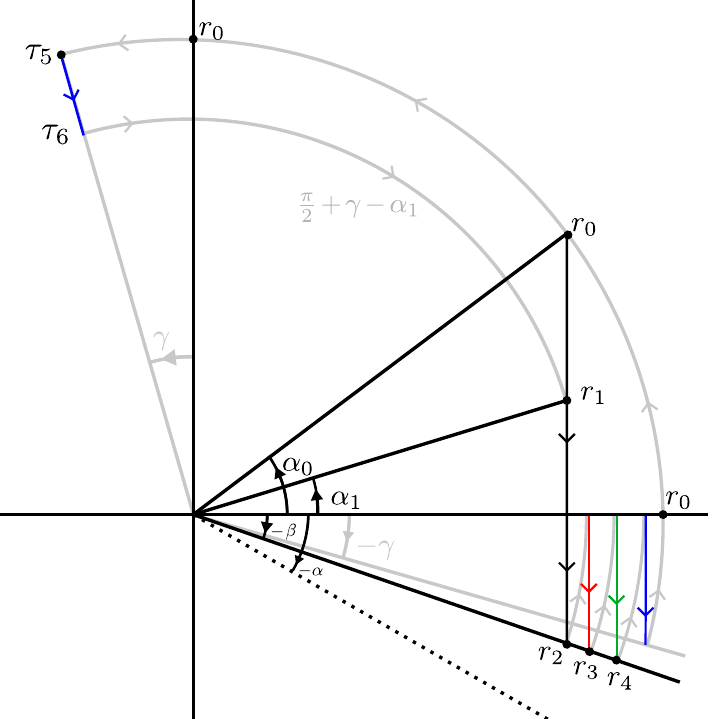} \\ (a) \\ \vskip 1cm 
\includegraphics[height=10cm]{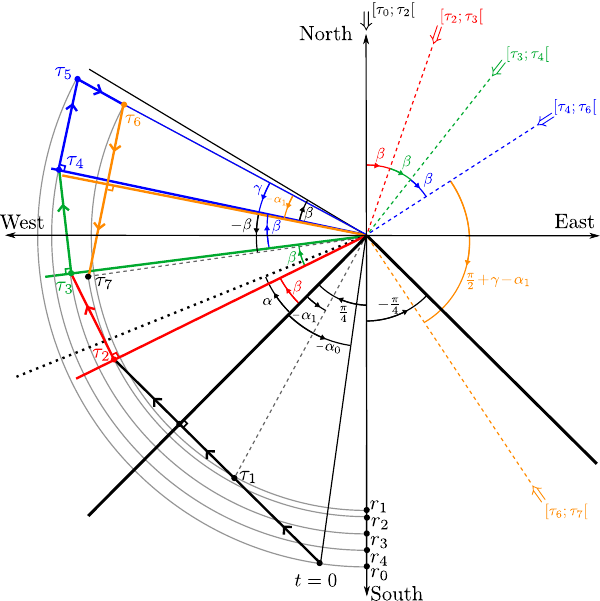} \\
(b)
\end{center}
\caption{The path of the boat when the wind angle is given by \eqref{rd11_24e1}, the boat's initial position is $(r_0, \alpha_0, 1)$ and the boat follows strategy $A$ of Section \ref{Sec:4}.  The path in picture (a) (respectively (b)) is in the rotating (respectively geographic) frame. In picture (b), the double arrows indicate the wind direction during each of the time intervals $[\tau_0,\tau_2[,\dots,[\tau_6, \tau_7[$.  \label{fig9}}
\end{figure}

We now describe a wind behavior that could prevent a boat that follows the strategy $A$ described in Section \ref{Sec:4} from reaching the target. Assume that the initial position of the boat is $(r_0, \alpha_0, 1)$, where $r_0 > \eta$ and $\alpha_0$ is close to $\piq$, that is, the boat is almost at equal distances from both laylines, at distance $r_0$ from the origin, and on starboard tack.  Fix $\beta \in \, ]0, \alpha[$ and $r_1 \in \, ]0, r_0[$ close to $r_0$. Let $\tau_1$ be the first time at which $R_t = r_1$ and suppose that $\T_t = \alpha_1 \in \, ]0, \alpha_0[$ at this time. We assume that the wind is constant up to the first time $\tau_2 > \tau_1$ that $\T_t = - \beta$. At this time, $R_t = r_2 < r_0$. Suppose that at this time, the wind turns $\beta$ radians to the right, so that the boat is on the new port layline. The boat continues sailing on starboard tack up to the first time $\tau_3 > \tau_2$ that $\T_t = -\beta$, at which time $R_t = r_3 \in\,] r_2, r_0[$. Suppose that the wind again shifts $\beta$ degrees to the right, and that this pattern is repeated at time $\tau_4$. Then, as the boat continues on starboard tack, it hits the set $D^1_2$ at time $\tau_5 > \tau_4$, and we set $-\gamma = \T_{\tau_5}$. At this time, the boat tacks to port, as prescribed by the strategy $A$. It is now back to a distance $r_0$ from the origin. The boat continues sailing on port tack until it is at distance $r_1$ from the origin, which occurs at time $\tau_6$. At this time, the wind turns $\pim + \gamma - \alpha_1$ radians to the right (in the geographic frame). This situation is equivalent to the boat's position at time $\tau_1$ (when it was on on starboard tack): the boat continues on port tack until the first time $\tau_7$ where $\T_t = \pim + \beta$, etc.,  and the process can be repeated. The boat will forever remain between distances $r_1 \cos(\alpha_1)$ and $r_0$ from the origin.

In this scenario, the wind direction $\beta_t$ (relative to the North to South direction) is given by the formula
\begin{align}\label{rd11_24e1}
   \beta_t = \beta\, \left(1_{[\tau_2, \tau_3[}(t) + 2\, 1_{[\tau_3, \tau_4[}(t) + 3\, 1_{[\tau_4, \tau_6[}(t) \right) + \left(3 \beta + \pim + \gamma - \alpha_1\right) \, 1_{[\tau_6, \tau_7[}(t).
\end{align}
The path of the boat in the rotating (respectively geographic) reference frame is shown in picture (a) (respectively (b)) of Fig.~\ref{fig9}. 

We have chosen formula \eqref{rd11_24e1} because it leads to a path for the boat that is not too difficult to describe. Obviously, the indicator functions on the right-hand side of \eqref{rd11_24e1} can be smoothed into continuous functions, and any minor variation around this formula will lead to a similar behavior. Theorem \ref{thm_brw_exist_V} tells us that even though there are (infinitely) many wind behaviors that have the same effect as the scenario described here, the set of all such behaviors has probability $0$ when the wind direction is given by a Brownian motion on a circle.

\section{Appendix B: free boundary problem and smooth fit}\label{appB}

In this appendix, we examine in more detail the HJB equations presented at the ends of Subsections \ref{Subs:1.2.1} and \ref{rd07_15ss1}. 
 
\subsection{The case $c>0$ and $\eta = 0$}

According to the HJB equation \eqref{eq:QVI}, in the ``no-tacking'', or ``continuation'' region $C$ defined in Theorem \ref{rd06_29t1}, the value function $V = V^{c,0}$ should satisfy the PDE $\mathcal{L} V(r,\t,a) + 1 = 0$. Taking into account the $2 \pi$-periodic form of the coefficients of $\mu_1$ and $\mu_2$ in the definition of $\mathcal{L}$ (see \eqref{eq_drift_pa}), we can restrict this PDE to the domain $]0, \infty[\, \times \, ]{-}\frac{\pi}{4}, \frac{7\pi}{4}[ \times \{1, -1\}$. 

We write $C = (C_1 \times \{1\}) \cup(C_{-1} \times \{-1\})$, where $C_1$ and $C_{-1}$ are subsets of $]0, \infty[\, \times \, ]{-}\frac{\pi}{4}, \frac{7\pi}{4}[$.  Using the symmetry property between tacks discussed in \eqref{eq_brw_cartesian_symmetry} and \eqref{eq_brw_polar_symmetry}, we can reduce to a single PDE
\begin{align}\label{rd07_15e1}
    \mathcal{L} V(r,\t,1) + 1 = 0, \qquad (r,\t) \in C_1,
\end{align}
with the additional equation
\begin{align*}
    V(r,\t,1) &= V(r, \t,-1) + c \\
     & = V\left(r, \tfrac{\pi}{2} - \t,1\right) + c, \qquad (r,\t) \in \left(]0, \infty[\, \times \left[{-}\tfrac{\pi}{4}, \tfrac{7\pi}{4}\right] \right) \setminus C_1,
\end{align*}
where $\tfrac{\pi}{2} - \t$ is computed modulo $2\pi$ and belongs to $]{-}\tfrac{\pi}{4}, \tfrac{7\pi}{4}[$. In this PDE, the domain $C_1$ is unknown and determining it is part of the problem of solving the HJB equation.

We conjecture that for the optimal strategy, the boundary of $C_1$ consists of the union of two curves $\Gamma_1$ and $\Gamma_2$ with respective equations $\t = \gamma(r)$ and $\t = \frac{3\pi}{2} - \gamma(r)$ (because of the symmetry \eqref{eq_brw_polar_symmetry_2}), where $\gamma: [0,\infty[\, \to \, ]{-}\frac{\pi}{4}, \frac{7\pi}{4}[$. We call these curves {\em stochastic laylines}, because when the boat is sailing in the continuation region and reaches one of these curves, then it should tack. 

\begin{figure}[ht]
\begin{center}
\includegraphics[width=0.7\textwidth]{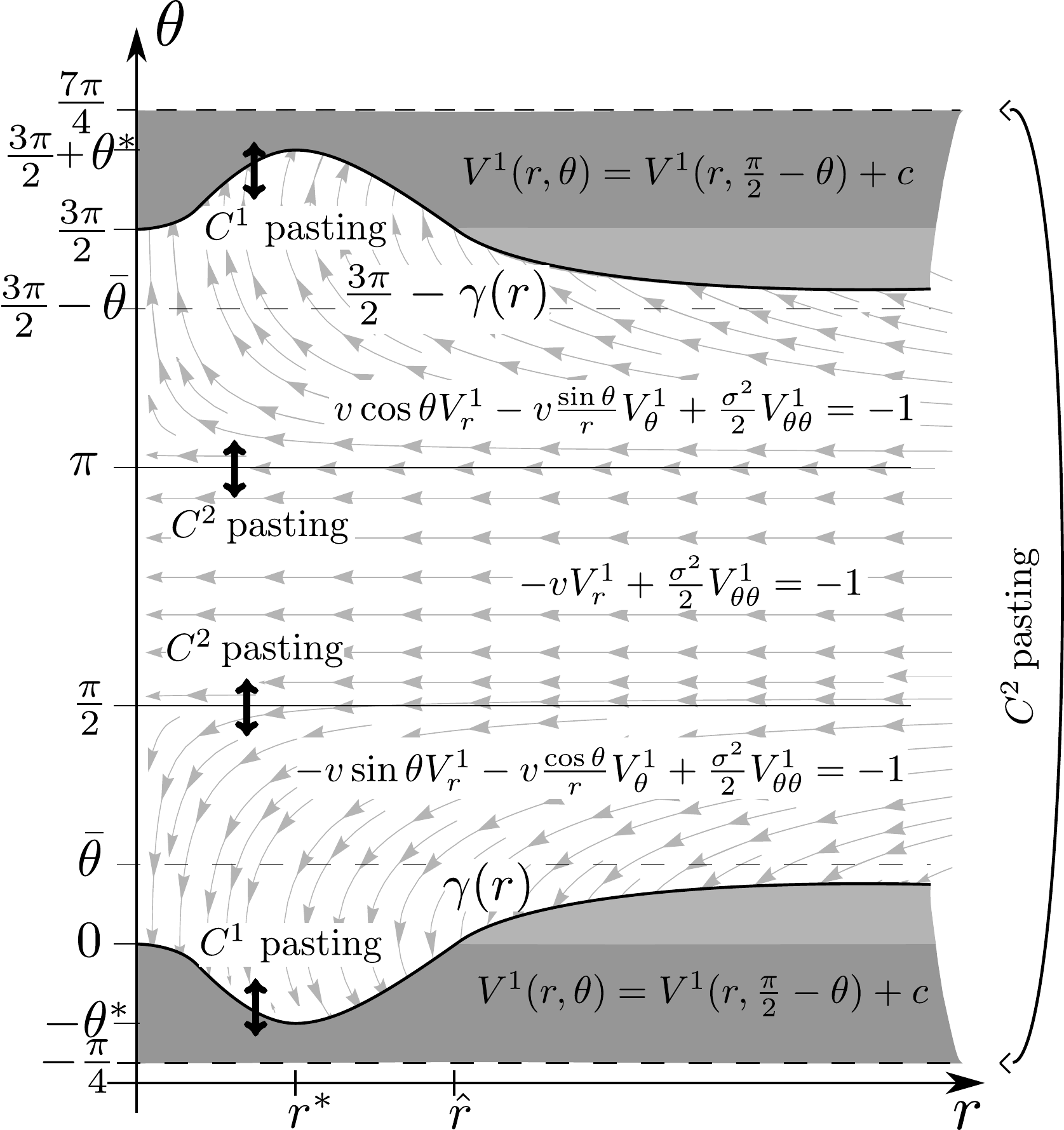}
\end{center}
\caption{Conjectured free-boundary problem in the $(r,\theta)$-plane for
the value function $V^{c,0}$ on starboard tack. The grey areas are the tacking regions,
and the white area is the continuation region. The curves
$\theta=\gamma(r)$ and $\theta=3\pi/2-\gamma(r)$ are the stochastic laylines.
The arrows represent the drift field.}
\label{fig:free-boundary-problem}
\end{figure}

The conjectured curves $\Gamma_1$ and $\Gamma_2$ are shown in Figure \ref{fig:free-boundary-problem}. Because of the form of the coefficients $\mu_1$ and $\mu_2$ in the definition of $\mathcal{L}$ (see \eqref{eq_drift_pa}), the PDE \eqref{rd07_15e1} has different expressions depending on whether $\t$ belongs to $]{-}\frac{\pi}{4}, \frac{\pi}{2}[$, $]\frac{\pi}{2}, \pi[$, or $]\pi, \frac{7\pi}{4}[$, as shown in Figure \ref{fig:free-boundary-problem}.

Because of the symmetry \eqref{eq_brw_polar_symmetry_2}, we can further reduce the domain of the PDE to $]{-}\frac{\pi}{4}, \frac{3\pi}{4}[$, and this symmetry also implies that $\frac{\partial}{\partial \t} V(r, \frac{3\pi}{4}, 1) = 0$. 

Summarizing, letting $V^1(r, \t) := V(r,\t,1)$, we seek a function $V^1: [0,\infty[\, \times [{-}\tfrac{\pi}{4}, \tfrac{3\pi}{4}] \to \R$ and $\gamma: [0,\infty[\, \to \, ]{-}\frac{\pi}{4}, \frac{7\pi}{4}[$ with $\gamma(0) = 0$ such that:
\begin{align}
\begin{cases}
  \text{(a) } v \sin(\t \wedge \tfrac{\pi}{2})\, \tfrac{\partial}{\partial r} V^1(r,\t) + v\, \tfrac{\cos\left(\t \wedge \tfrac{\pi}{2}\right)}{r}\,\tfrac{\partial}{\partial \t} V^1(r,\t) \\
  \quad\qquad\qquad\qquad\qquad\qquad - \tfrac{\sigma^2}{2} \tfrac{\partial^2}{\partial \t^2} V^1(r,\t) = 1, \qquad r > 0, \ \gamma(r) < \t < \tfrac{3\pi}{4}, \\
  \\ 
  \text{(b) }  V^1(0, \t) = \begin{cases} c, \qquad {-}\frac{\pi}{4} \leq \t < \tfrac{\pi}{2},  \\
          0, \qquad \tfrac{\pi}{2} < \t \leq \tfrac{3\pi}{4},
       \end{cases}
   \\ 
  \\
     \text{(c) } \tfrac{\partial}{\partial \t} V^1(r,\frac{3\pi}{4}) = 0,\qquad r > 0,\\
  \\
  \text{(d) }  V^1(r,\t) =  V^1\left(r, \tfrac{\pi}{2} - \t \right) + c, \qquad r > 0, \ {-}\tfrac{\pi}{4} \leq \t \leq \gamma(r).
  \end{cases}
\end{align}
We require that $V^1$ be $C^{1,2}$ in the region where the PDE (a) is satisfied. In the conditions (a)--(d), there is nothing that imposes that $V^1$ is continuous at $(r, \gamma(r))$. If $\gamma(r)$ were given, then imposing the {\em continuous pasting} condition 
$$
  \text{(e) } \quad \lim_{\t \uparrow \gamma(r)}V^1(r, \t) = \lim_{\t \downarrow \gamma(r)}V^1(r, \t),\qquad r > 0,
$$
would suffice to determine $V^1$. Since $\gamma(r)$ is not given, following \cite[Section 8]{peskshir}, we also require the {\em smooth pasting} condition
$$
  \text{(f) } \quad \lim_{\t \uparrow \gamma(r)} \tfrac{\partial}{\partial \t} V^1(r, \t) = \lim_{\t \downarrow \gamma(r)} \tfrac{\partial}{\partial \t} V^1(r, \t),\qquad r > 0.
$$
We conjecture that conditions (a)--(f) determine $V^1$ and $\gamma$ and that the unique solution $V^1$ of (a)--(f) is the value function $V^{c,0}$.

\subsection{The case $c = 0$ and $\eta > 0$}

According to the HJB equation \eqref{eq:HJB-c-zero}, in the regions $C_1$ and $C_{-1}$ defined in Theorem \ref{rd07_01t1}, the value function $V = V^{\eta}$ should satisfy the PDEs $\mathcal{L}^1 V(r,\t) + 1 = 0$ and $\mathcal{L}^{-1} V(r,\t) + 1 = 0$, respectively.

Taking into account the $2 \pi$-periodic form of the coefficients of $\mu_1$ and $\mu_2$ in the definition of $\mathcal{L}$ (see \eqref{eq_drift_pa}), we can restrict this PDE to the domain $]\eta, \infty[\, \times \, ]{-}\frac{3\pi}{4}, \frac{5\pi}{4}[$. 
 
 Looking at Figure \ref{Fig:9}, for the optimal strategy we should have $C_1 = \,]\eta, \infty[\, \times \, ]\frac{\pi}{4}, \frac{5\pi}{4}[$ and $C_{-1} = \,]\eta, \infty[\, \times \, ]{-}\frac{3\pi}{4}, \frac{\pi}{4}[$. The stochastic laylines are no longer stochastic, and in fact, $\gamma(r) \equiv \frac{\pi}{4}$. There are two symmetries present in Figure \ref{Fig:9}, with respect to the line $\t =  \frac{\pi}{4}$ and $\t = \frac{3\pi}{4}$. These imply that that we can just consider the PDE $\mathcal{L}^1 V(r,\t) + 1 = 0$ on $]\eta, \infty[\, \times \, ]\frac{\pi}{4}, \frac{3\pi}{4}[$ only, that is,
\begin{align}\nonumber
    & v \sin(\t \wedge \tfrac{\pi}{2})\, \tfrac{\partial}{\partial r} V^1(r,\t) + v\, \tfrac{\cos\left(\t \wedge \tfrac{\pi}{2}\right)}{r}\,\tfrac{\partial}{\partial \t} V^1(r,\t) \\
  &\qquad\qquad\qquad\qquad\qquad - \tfrac{\sigma^2}{2} \tfrac{\partial^2}{\partial \t^2} V^1(r,\t) = 1, \qquad r > \eta, \ \tfrac{\pi}{4} < \t < \tfrac{3\pi}{4}, 
\label{rd07_16e1}
\end{align}
with the vanishing Neumann boundary conditions
\begin{align}
     \tfrac{\partial}{\partial \t} V^1 (r,  \tfrac{\pi}{4}) = 0 = \tfrac{\partial}{\partial \t} V^1 (r,  \tfrac{3\pi}{4}), \qquad r > 0,
\label{rd07_16e2}
\end{align}
and vanishing initial conditions
\begin{align}
   V^1(\eta, \t) = 0,\qquad \tfrac{\pi}{4} \leq \t \leq \tfrac{3\pi}{4}.
\label{rd07_16e3}
\end{align}
We conjecture that the value function $V^\eta(r,\t,1)$ is the unique solution to the PDE \eqref{rd07_16e1} with the boundary conditions \eqref{rd07_16e2} and initial condition \eqref{rd07_16e3}.

\end{document}